%% file: PQR_baseclass.tex
\documentclass[10pt,a4paper]{article}

\usepackage{fancyhdr}
\usepackage{a4wide}
\usepackage{authblk}
\usepackage{amssymb}
\usepackage[figuresright]{rotating}

\usepackage{siunitx}            

\usepackage{amsmath}
\usepackage{bm}                 
\usepackage{tikz}
\usetikzlibrary{patterns}
\usetikzlibrary{arrows}
\usepackage{subfigure}
\usepackage{booktabs}          
\pgfdeclarepatternformonly{north east lines wide}{\pgfqpoint{-1pt}{-1pt}}{\pgfqpoint{10pt}{10pt}}{\pgfqpoint{9pt}{9pt}}%
{
  \pgfsetlinewidth{0.4pt}
  \pgfpathmoveto{\pgfqpoint{0pt}{0pt}}
  \pgfpathlineto{\pgfqpoint{9.1pt}{9.1pt}}
  \pgfusepath{stroke}
}

\fancypagestyle{FirstPage}{
\lfoot{\copyright~2017. This manuscript version is made available under the CC-BY-NC-ND 4.0 license http://creativecommons.org/licenses/by-nc-nd/4.0/}
}

\DeclareMathOperator{\ddiv}{div}
\DeclareMathOperator{\Span}{span}
\DeclareMathOperator{\argmin}{\arg\,\min}
\newcommand{\re}{\mathrm{Re}}
\newcommand{\AR}{\mathrm{AR}}
\newcommand{\HH}{\mathcal{H}}
\newcommand{\bs}{\boldsymbol}
\newcommand*{\de}{\mathop{}\!\mathrm{d}}
\newcommand*{\Proj}{\mathit{\Pi}}
\newcommand*{\Om}{\mathit{\Omega}}
\newcommand*{\Omh}{\widehat{\mathit{\Omega}}}
\newcommand*{\Ga}{\mathit{\Gamma}}
{\left\lbrace\begin{array}{@{}l@{}}}%
{\end{array}\right.}
\def\ss{\boldsymbol{\sigma}}
\def\mmu{\boldsymbol{\mu}}
\def\ee{\boldsymbol{\epsilon}}

\def\el {\nonumber }

\usepackage[natbib=true,backend=bibtex]{biblatex}
\bibliography{rbsissa,latexbi}

\begin{document}

\date{}

\title{Computational reduction strategies for the detection of steady bifurcations in incompressible fluid-dynamics: applications to Coanda effect in cardiology}

\author[1]{Giuseppe Pitton}
\author[2]{Annalisa Quaini}
\author[1]{Gianluigi Rozza}

\affil[1]{SISSA, International School for Advanced Studies, Mathematics Area, mathLab, Via Bonomea 265, 34136 Trieste, Italy. \\ Email: \texttt{giuseppe.pitton@sissa.it, gianluigi.rozza@sissa.it}}
\affil[2]{University of Houston, Department of Mathematics, Houston, TX, USA. Email: \texttt{quaini@math.uh.edu}}

\maketitle

\thispagestyle{FirstPage}

\begin{abstract}
We focus on reducing the computational costs associated with the hydrodynamic stability of solutions of the incompressible 
Navier-Stokes equations for a Newtonian and viscous
fluid in contraction-expansion channels. In particular, we are interested in studying 
steady bifurcations, occurring when non-unique stable solutions  appear as physical and/or geometric 
control parameters are varied.  
The formulation of the stability problem requires solving an eigenvalue problem for a partial differential operator. 
An alternative to this approach is the direct simulation of the flow to characterize
the asymptotic behavior of the solution.
Both approaches can be extremely expensive in terms
of computational time.
We propose to apply Reduced Order Modeling (ROM) techniques to reduce
the demanding computational costs associated with 
the detection of a type of steady bifurcations in fluid dynamics.
The application that motivated the present study is the onset of asymmetries (i.e., symmetry breaking bifurcation) 
in blood flow through a regurgitant mitral valve, depending on the Reynolds number
and the regurgitant mitral valve orifice shape.

\end{abstract}

\input{tex/introduction_revised.tex}
\input{tex/problem_setting_revised.tex}
\input{tex/numerical_method_revised.tex}
\input{tex/results_revised.tex}

\input{tex/conclusions.tex}

\input{tex/perspectives.tex}

\printbibliography
\end{document}

%% file: tex/introduction_revised.tex
\section{Introduction}\label{motivation}

We focus on the hydrodynamic stability of solutions of the incompressible 
Navier-Stokes equations for a Newtonian and viscous
fluid in contraction-expansion channels, with a particular concern on 
steady bifurcations.
Steady bifurcations occur when new, non-unique solution branches of the
Navier-Stokes equations appear as physical and/or geometric 
control parameters are varied.  
When the fluid domain is characterized by two or three dimensions with non-periodic
boundary conditions, the formulation of the stability problem requires solving an eigenvalue problem for a partial differential operator. See \cite{Dijkstra} for a review on numerical methods for stability analysis
based on linearized eigenvalue problems. 
An alternative to the eigenvalue problem approach is the direct simulation of the flow to characterize
the asymptotic behavior of the solution; see, e.g., \cite{goodrichg,autierip1,AQpreprint}.
Both approaches can be extremely expensive in terms
of computational time.
In this paper, we propose to apply Reduced Order Modeling (ROM) techniques to reduce
the demanding computational costs associated with flow stability analysis.

\begin{figure}
\begin{center}
    \subfigure[Location of Mitral Valve]{
        \hskip 0.2in\includegraphics[width = 0.23 \textwidth]{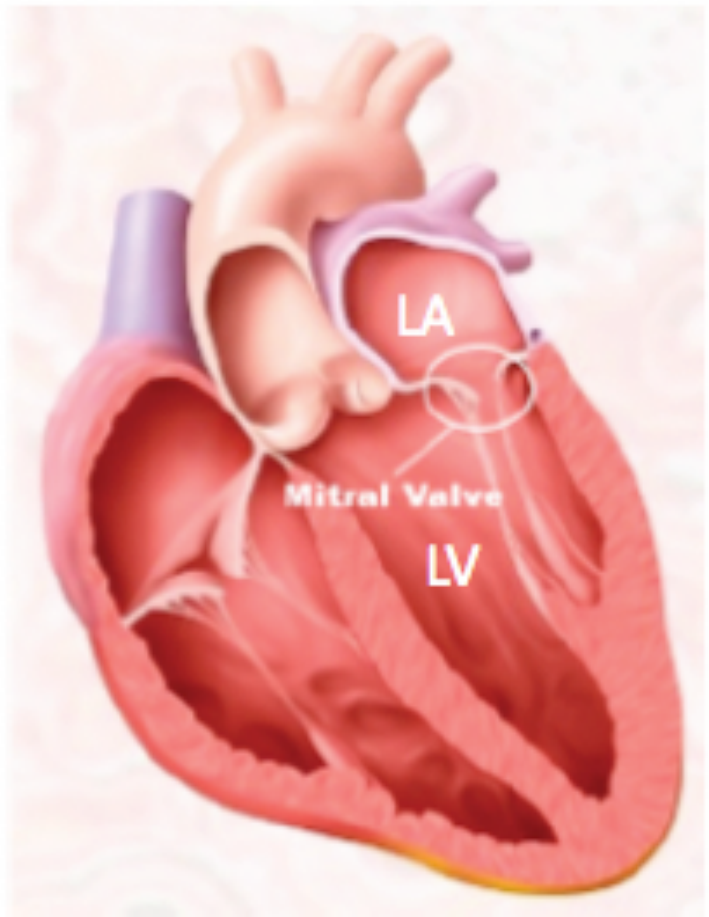}\hskip 0.2in
    }
    \subfigure[Central color Doppler jet ]{\hskip 0.2in
        \includegraphics[width = 0.23 \textwidth]{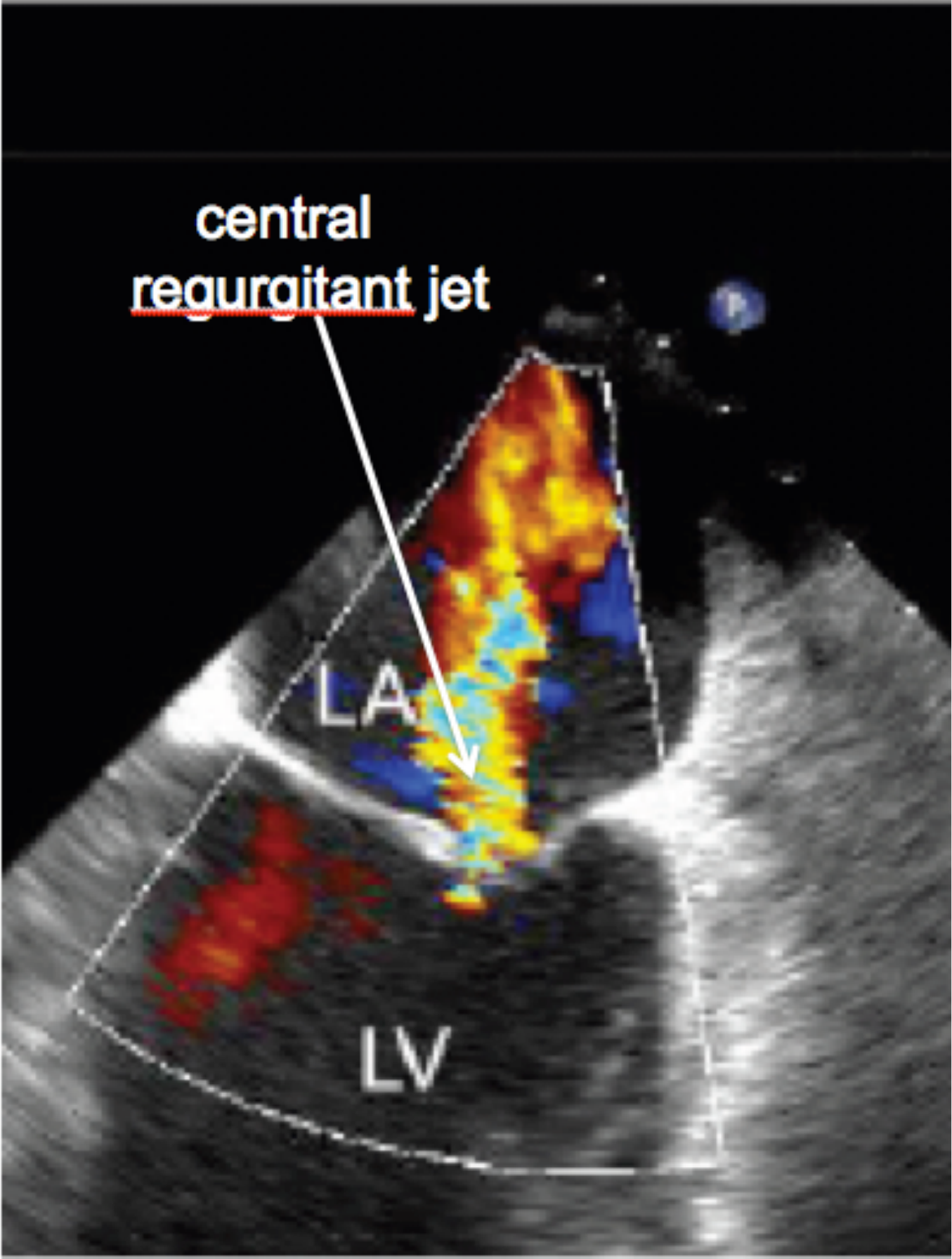}\hskip 0.2in
    }
    \subfigure[Eccentric color Doppler jet ]{\hskip 0.2in
        \includegraphics[width = 0.23 \textwidth]{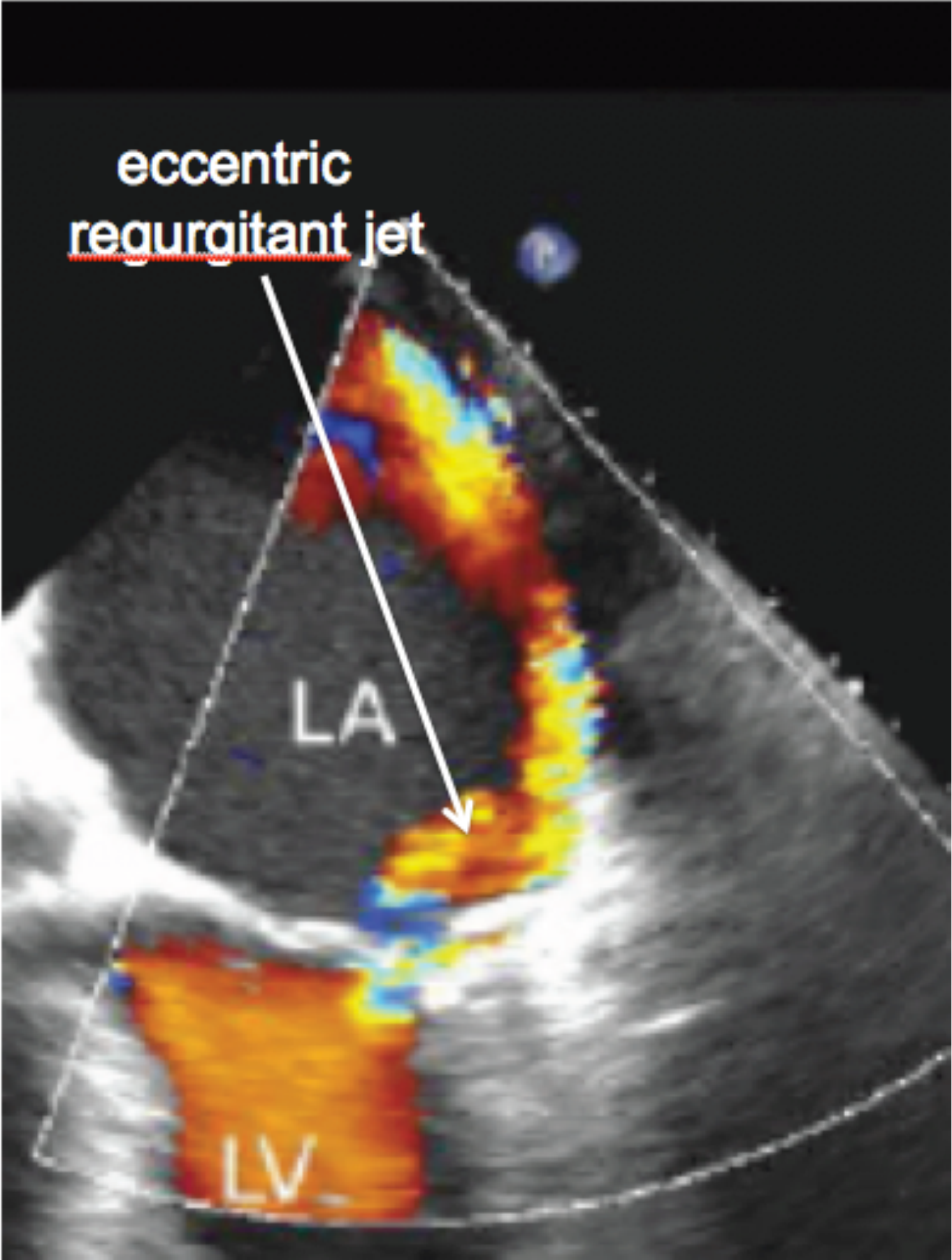}\hskip 0.2in
    }
\end{center}
        \caption{\small (a) Anatomy of the heart showing the mitral valve.
(b) Echocardiographic image of central regurgitant jet flowing from
the left ventricle (LV) to the left atrium (LA).
Colors denote different fluid velocities.
(c) Echocardiographic image of eccentric regurgitant jet, hugging the walls of the left atrium (LA)
known as the Coanda effect.}
\label{MitralValve}
\end{figure}

Practical applications of contraction-expansion channel flows
include equipments such as heat exchangers, combustion
chambers, and mixing vessel. 
An application that motivated the present study is the onset of asymmetries (i.e., symmetry breaking bifurcation) 
in blood flow through a regurgitant mitral valve, depending on the Reynolds number
and the regurgitant mitral valve orifice shape.
Mitral regurgitation is a valvular disease characterized by abnormal leaking of blood
through the mitral valve from the left ventricle into the left atrium of the heart. See Figure~\ref{MitralValve}.
In certain cases  the regurgitant jet ``hugs" the wall of the heart's atrium as shown in  Figure~\ref{MitralValve}(c).
These eccentric, wall-hugging, non-symmetric regurgitant jets have been observed at low
Reynolds numbers \cite{regurgRe1000,regurgRe50} and are said to undergo the Coanda effect
\cite{CoandaBook,willef1}.
This effect, described as the tendency of a fluid jet to be attracted to a nearby surface, 
owes its name to Romanian aerodynamics pioneer Henri Coanda.
The primary tool to assess the severity of mitral regurgitation is echocardiography \cite{Zoghbi_et_al_2003}. 
One of the biggest challenges in echocardiographic assessment of mitral regurgitation is
the Coanda effect: the wall-hugging jets appear smaller in the color Doppler image of regurgitant flow,
leading to a gross under-estimation of regurgitant volume by inexperienced observers
\cite{ginghina1,chaom1}. As a result, patients requiring treatment may not be recognized.

Despite the large cardiovascular and bioengineering literature reporting on the Coanda effect in echocardiographic 
assessment of mitral regurgitation, there is very little connection with the fluid dynamics literature that could help
identify and understand the main features of the corresponding flow conditions. In this paper, our goal is to understand 
what triggers the Coanda effect in a simplified setting. 
A contraction-expansion channel is a simplified setting which has the same geometric features of mitral regurgitation. 
In fact, a mitral regurgitant jet flows from the left ventricle through the contraction between the
mitral leaflet, called regurgitant orifice, into the left atrium. 
First, we focus on planar contraction-expansion channels (see Fig.~\ref{fig:scheme_2d})
and investigate the influence of the Reynolds number and the contraction width $w_c$ (i.e., the orifice height)
on the flow.
Then, we consider the 3D geometry reported in Fig.~\ref{fig:scheme_3d} to understand the role played by
the channel depth $h$ (i.e., the orifice length). 
Eccentric regurgitant jets typically occur in prolapsed mitral valves, i.e. when two valve flaps of the mitral valve do not close evenly.
Thus, another parameter of interest, although not considered in this work, could be the orifice depth. 
Moreover, for a more realistic setting one would have to account for the pulsatility of the flow
and include the Strouhal number among the parameters. 

We remark that the focus of this paper is to investigate the cause of the Coanda effect in simplified settings. 
Nonetheless, it is thanks to the results reported here that our medical collaborators at the 
Houston Methodist DeBakey Heart \& Vascular Center 
were able to reproduce the Coanda effect in a mock heart chamber (see Sec.~\ref{3d_results}). A comparison between the experiments
in vitro and corresponding 3D simulations is presented in 
\cite{WangQuaini}.

The incompressible fluid dynamics in a planar contraction-expansion channel 
has been widely studied from both theoretical and practical perspectives; see, e.g., 
\cite{moffatt1,Drikakis,sobeyd1,fearnm1,hawar1,mishraj1} and references therein.
In the two-dimensional geometry reported in Fig.~\ref{fig:scheme_2d}, the wall-hugging effect 
happens only above a critical Reynolds number (\ref{eq:re_2d}), which depends on the 
expansion ratio $\lambda$ defined in (\ref{eq:lambda}). 
Compare Fig.~\ref{fig:snapshots_2d}(e) and Fig.~\ref{fig:snapshots_2d}(b), which correspond 
to a Reynolds number above and below the critical value, respectively.
The asymmetric, wall-hugging solution remains stable for a certain range of Reynolds number 
and asymmetries become stronger with the increasing Reynolds number, as shown in \cite{mishraj1}. 
The formation of stable asymmetric vortices in 2D planar expansion is attributed to an increase in velocity near one wall
that leads to a decrease in pressure near that wall \cite{willef1}. Once a pressure difference is established across 
the channel, it will maintain the asymmetry of the flow. 
The critical value of the Reynolds number has been identified for different expansion ratios $\lambda$. 
In particular, it was found that such critical value decreases with increasing value of $\lambda$ (see \cite{Drikakis,revuelta1}).

In the three-dimensional geometry reported in Fig.~\ref{fig:scheme_3d}, the critical Reynolds number
for the symmetry-breaking (i.e., the wall-hugging)
varies with the expansion ratio and the aspect ratio defined in (\ref{eq:aspect_ratio}), as shown in \cite{cherdrond1,Chiang2000,Oliveira}.
When the expansion ratio is fixed and the aspect ratio decreases, the endwall influence 
becomes more important: the critical Reynolds number 
increases \cite{Chiang2000,Oliveira}.
For moderate aspect ratios, the flow is steady in time but highly three-dimensional, and complex
spiraling structures are observed, which are not closed recirculating
cells as in the case of 2D flows. See Fig.~\ref{fig:cmp_vortices}.
The numerical studies in \cite{Tsai2006} found that
the flow only resembles a 2D flow for very large aspect
ratios.
The theoretical study of Lauga et al. \cite{Lauga2004}
shows that for low aspect ratios the flow is highly three-dimensional. 
The numerical and experimental
studies in \cite{Oliveira} show that the strong three-dimensional effects
appearing for low aspect ratios inhibit the wall-hugging effect observed 
in geometries with high aspect ratios at the same Reynolds number.
This suggests that the eccentric regurgitant jets, such as the one in Figure~\ref{MitralValve}(c),
occur when the regurgitant orifice is long (large aspect ratio) and narrow (large expansion ratio).

Given the relatively fast decay of energy spectrum for flows at sufficiently low Reynolds
numbers, a ROM technique is expected to be an efficient tool to
reduce the prohibitive computational costs associated to identifying
the flow conditions and geometries that trigger asymmetries.
Recent developments of ROM techniques have focused on the reduction of
computational time for a wide range of differential problems \cite{RozzaEncyclopedia,QuarteroniRozza2013}, while maintaining a prescribed tolerance on
error bounds \cite{Rozza:ARCME,hesthaven2015certified,QuarteroniManzoniNegri}.
Terragni and Vega \cite{Terragni:2012} showed that a 
Proper Orthogonal Decomposition (POD) approach allows for considerable
computational time savings for the analysis of bifurcations in some nonlinear dissipative systems. 
Herrero, Maday and Pla~\cite{Maday:RB2} have used a Reduced Basis (RB) method to speed up the computations of different solution branches of a two-dimensional natural convection problem (Rayleigh-B\'{e}nard), achieving a good accuracy but without investigating the approximation of bifurcation points.
For each fixed aspect ratio, multiple steady solutions for the Rayleigh-B\'enard
problem can be found for different Rayleigh numbers and stable solutions
coexist at the same values of external physical parameters.
In \cite{Maday:RB2}, it is shown that stable and unstable solutions are correctly identified by the RB method. 
Yano and Patera \cite{Yano:2013} introduced a RB
method for the stability of flows under perturbations in the forcing term or in the boundary conditions,
which is based on a space-time framework that allows for particularly sharp error estimates. Furthermore, in~\cite{Yano:2013} 
it is shown how a space-time inf-sup constant approaches zero as the computed solutions get close to a bifurcating value.

In a previous work \cite{PittonRozza} we have investigated steady and Hopf bifurcations in a 
natural convection problem dealing with a geometrical parameter (the cavity length), 
and a physical parameter (the Grashof number). 
This work is an extension of that study and provides a proof of concept of the applicability of reduced order methods 
to investigate stability and bifurcations in complex fluid dynamic problems
at a reasonable computational cost.
The proposed framework allows the use of a black-box input-output toolbox to be managed also 
by non-expert scientists in computational sciences. The offline-online splitting of the computational procedure 
is crucial in view of the use of a High Performance Computing (HPC) infrastructure for the offline computational step 
(expensive and time consuming) and a light modern device, such as tablet or smart phone, for online calculations. 
The idea is to use different platforms (and methodologies) for a strategic computational collaboration 
between high order and reduced order methods, with competitive computational costs for complex simulations. 
This could still be considered a research frontier in computational fluid dynamics, especially in view of real life applications.

The outline of the paper is as follows. The general problem setting is described in Section 2, while the numerical methods 
are treated in Section 3. Numerical results in 2D and 3D contraction-expansion channels 
under different parametrizations are reported in Section 4. Conclusions and perspectives follow in Section 5.

%% file: tex/problem_setting_revised.tex
\section{Problem setting}
\label{setting}

Let $\Om\subset\mathbb{R}^d$, $d=2,3$, be the computational domain. The motion of an incompressible, viscous fluid in a spatial domain $\Om$ over a time interval of interest $(0, T)$ is governed by the incompressible Navier-Stokes equations
\begin{align}
\frac{\partial \bs{u}}{\partial t}+(\bs{u} \cdot\nabla)\bs{u} - \nabla \cdot \ss =\bs{0} \qquad \text{in}~ \Om \times (0,T), \label{eq:uns_strong1} \\
\nabla \cdot \bs{u} =0 \qquad \text{in}~ \Om\times (0,T), \label{eq:uns_strong2}
\end{align}
where $\bs{u}$ is the velocity and $\ss$ is  the Cauchy stress tensor. In large arteries and inside the heart, it is widely accepted to model blood as a Newtonian fluid. See, e.g., \cite{B-haemodel} and references therein. For such fluids, 
$\boldsymbol{\sigma} (\bs{u}, p) = -p \bold{I} +2\nu \boldsymbol{\epsilon}(\bs{u})$, 
where $p$ is the pressure, $\nu$ is the fluid kinematic viscosity, and 
$\boldsymbol{\epsilon}(\bs{u}) = (\nabla \bs{u} + (\nabla\bs{u})^T)/2$
is the strain rate tensor.
Eq. (\ref{eq:uns_strong1}) represents the conservation of the linear momentum, while eq. (\ref{eq:uns_strong2}) represents the conservation of the mass.
Eq. (\ref{eq:uns_strong1})-(\ref{eq:uns_strong2}) need to be endowed with boundary and initial conditions, e.g.:
\begin{align}
\bs{u} &= \bs{d} \quad \text{on} ~ \Ga_D \times (0,T), \label{eq:dir_bc} \\
\ss \cdot \bs{n} &= \bs{g}\quad \text{on}  ~ \Ga_N \times (0,T), \label{eq:neu_bc} \\
\bs{u} &= \bs{u}_0 \quad \text{in}~ \Om \times \{ 0 \}, \el
\end{align}
where $\overline{\Ga}_D \cup \overline{\Ga}_N = \partial \Om$ and $\Ga_D \cap \Ga_N = \emptyset$. 
Here, $\bs{d}$, $\bs{g}$, and $\bs{u}_0$ are given. We assume that there is no external body force acting on the fluid 
(see eq. (\ref{eq:uns_strong1})) and the motion is driven by the boundary conditions (\ref{eq:dir_bc})-(\ref{eq:neu_bc}). 

When the fluid acceleration is negligible (i.e., the system has evolved towards a steady state), 
eq. (\ref{eq:uns_strong1})-(\ref{eq:uns_strong2}) can be replaced by: 
\begin{align}
(\bs{u} \cdot\nabla)\bs{u} - \nabla \cdot \ss =\bs{0} \qquad \text{in}~ \Om, \label{eq:ns_strong1} \\
\nabla \cdot \bs{u} =0 \qquad \text{in}~ \Om, \label{eq:ns_strong2}
\end{align}

To characterize the flow regime under consideration, we define the Reynolds number as
\begin{equation}\label{eq:re}
\re = \frac{U L}{\nu},
\end{equation}
where $U$ and $L$ are characteristic macroscopic velocity and length respectively. We will characterize $U$ and $L$ for the specific cases under consideration in Sec. \ref{2dcase} and \ref{3dcase}. 

For the variational formulation of the fluid problem (\ref{eq:uns_strong1})-(\ref{eq:uns_strong2}), we indicate with $L^2(\Om)$ the space of square integrable functions in $\Om$ and with $H^1(\Om)$ the space of functions in $L^2(\Om)$ with first derivatives in $L^2(\Om)$. We use $(\cdot,\cdot)_\Om$ and $\langle \cdot , \cdot \rangle_\Om$ to denote the $L^2$ product and the duality pairing between $H^{1/2}(\Om)$ and $H^{-1/2}(\Om)$, respectively. Moreover, let us define:
\begin{equation*}
\begin{aligned}
\bs{V_D}=H^1_D(\Om):=\{\bs{v}\in H^1(\Om)\text{ s.t. } \bs{v}=\bs{d}\text{ on }\Ga_D \}, \quad \bs{V}&=H^1_0(\Om):=\{\bs{v}\in H^1(\Om)\text{ s.t. }\bs{v}=\bs{0}\text{ on }\Ga_D \}, \quad
Q&=L^2(\Om).
\end{aligned}
\label{eq:whoare_V_Q}
\end{equation*}

The variational formulation of problem (\ref{eq:uns_strong1})-(\ref{eq:uns_strong2}) with boundary conditions (\ref{eq:dir_bc})-(\ref{eq:neu_bc}) reads: find $(\bs{u},p)\in \bs{V_D}\times Q$ such that
\begin{align}
\left( \frac{\partial \bs{u}}{\partial t}, \bs{v} \right)_\Om + 
(\bs{v},(\bs{u}\cdot\nabla)\bs{u})_\Om+\nu(\ee(\bs{v}),\ee(\bs{u}))_\Om-(\nabla \cdot \bs{v},p)_\Om&= \langle \bs{g}, \bs{v}\rangle_{\Ga_N} \qquad \forall\bs{v}\in\bs{V},\label{eq:weak_ns1} \\
(q,\ddiv \bs{u})_\Om&=0 \hspace{1.6cm}  \forall q\in Q. \label{eq:weak_ns2} 
\end{align}
If eq.~(\ref{eq:ns_strong1})-(\ref{eq:ns_strong2}) are used to model the fluid dynamics, then the first term 
in eq.~(\ref{eq:weak_ns1}) is disregarded.


The nonlinearity in problem (\ref{eq:uns_strong1})-(\ref{eq:uns_strong2}) can produce a loss of uniqueness for 
the solution, with multiple solutions branching from a known solution
at a bifurcation point. As will be explained in Sec.~\ref{2dcase}, the Coanda effect
is associated with a steady bifurcation point.
To detect numerically the presence of a steady bifurcation point,
we will rely on the spectrum analysis of a linearized operator. See, e.g., \cite{Prodi} for a theoretical introduction
to bifurcation theory, and~\cite{Cliffe_ST} and~\cite{Dijkstra} for applications to numerical analysis.
Following~\citep[par.~7.3]{Cliffe_ST}, we introduce the linearization $\mathcal{L}:\bs{V_D}\times\bs{V}\to\bs{V}$ 
of the convection 
operator in eq.~(\ref{eq:ns_strong1}), obtained by Fr\'{e}chet differentiation about a base point $\bs{u}^*$ of the term $\bs{u}\cdot\nabla\bs{u}$:
\begin{equation}
\mathcal{L}(\bs{u}^*)[\bs{v}]=\bs{u}^*\cdot\nabla\bs{v}+\bs{v}\cdot\nabla\bs{u}^*.
\label{eq:def_L_lin}
\end{equation}
At a symmetry breaking bifurcation point, a simple eigenvalue of $\mathcal{L}$ changes sign. 

We remark that for the kind of bifurcations we are interested in (namely supercritical pitchfork bifurcations) it is sufficient to study the spectrum of the antisymmetric part of the linearized operator~\cite{Cliffe_ST}. Other techniques would have to be used 
for different kind of singular points (e.g., Hopf bifurcations, fold points) occurring at higher Reynolds numbers or in different settings.

\subsection{2D case}
\label{2dcase}
\begin{figure}
\centering
\includegraphics[width=0.5\textwidth]{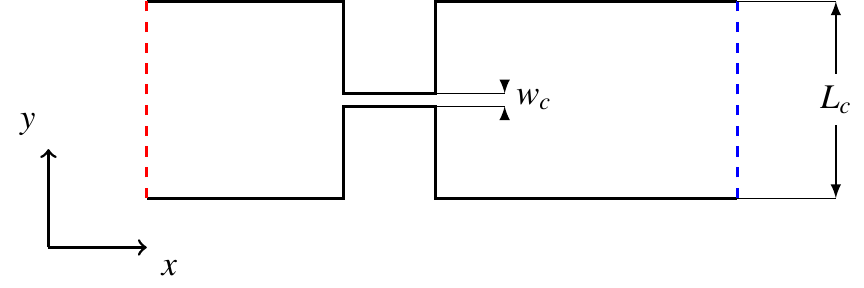}
\caption{Scheme of the planar contraction-expansion channel considered in this work. The notation is the same used in~\cite{Oliveira}.}
\label{fig:scheme_2d}
\end{figure}

The domain under consideration for the two dimensional case is shown in figure~\ref{fig:scheme_2d}. The following boundary conditions are imposed in this case: homogeneous Dirichlet (no-slip) boundary condition on the sides drawn with a continuous line, non-homogeneous Dirichlet boundary condition (time-independent parabolic velocity profile) on the red dashed line which 
corresponds to the inlet, and homogeneous Neumann (stress-free) boundary condition on the blue dashed line
which corresponds to the outlet. 
For space limitation, the channel depicted in Fig.~\ref{fig:scheme_2d} is shorter than the actual one. The actual domain length past the expansion is 6 times the channel height $L_c$. The length of the contraction channel is equal to $(L_c-w_c)/2$, and the distance between the inflow section and the contraction is equal to $L_c$.

For the characterization of the flow in the 2D case, we introduce the following quantities: 
\begin{itemize}
\item[-] {\bf expansion ratio}: 
\begin{align}\label{eq:lambda}
\displaystyle{\lambda=\frac{L_c}{w_c}};
\end{align}
\item[-] {\bf average horizontal velocity}: $\displaystyle{\langle v_x\rangle=\frac{Q}{w_c}}$, where $Q$ is the flow rate;
\item[-] {\bf Reynolds number}: 
\begin{align}\label{eq:re_2d}
\displaystyle{\re_{2D}= 2\frac{\langle v_x\rangle w_c}{\nu}}.
 \end{align}
\end{itemize}

Notice that the above definition of $\re_{2D}$ does not coincide with taking $L = w_c$ and $U = \langle v_x\rangle$ in (\ref{eq:re}). 
The reason for the extra factor 2 will be explained in the next subsection. In the numerical simulations, 
we vary $\re_{2D}$ by changing the value of the viscosity $\nu$.

In the geometry considered by~\cite{Oliveira} the expansion ratio $\lambda$ is equal to $15.4$. 
In order to reproduce the results in \cite{Oliveira} for validation purposes, in Sec.~\ref{2d_results_1param} we set 
$\lambda =15.4$ and we focus on the interval $\re\in[0.01,90]$.
In this interval, the flow configuration evolves as follows when the Reynolds numbers increases:
\begin{itemize}
\item[-] {\bf Creeping flow}: for very low Reynolds numbers the velocity field presents a double symmetry, with respect to both the horizontal and vertical symmetry axes of the domain geometry. See Fig.~\ref{fig:snapshots_2d}(a).
\item[-] {\bf Symmetric jet}: for slightly larger Reynolds numbers there is a breaking of the vertical symmetry. The flow is still symmetric with respect to the horizontal axis, but the two vortices downstream of the expansion are larger than the vortices upstream. See Fig.~\ref{fig:snapshots_2d}(b). At a further increase of the Reynolds number, the vertical asymmetry of the flow becomes increasingly evident,
 yet the horizontal symmetry is maintained.
\item[-] {\bf Asymmetric jets}: when the Reynolds number is sufficiently large, the configuration with a symmetric jet is still possible, 
but unstable \cite{sobeyd1}. In fact, small perturbations \footnote{That can be realized in several ways, e.g. with a slight variation of the boundary conditions, forcing terms, or superimposing a small random field to an established flow field and using this as a new initial condition.} expand one recirculation zone and shrink the other, causing a drastic variation in the flow.
See Fig.~\ref{fig:snapshots_2d}(c) and (d). This horizontally asymmetric solution remains stable for a certain range of $Re_{2D}$ and asymmetries become stronger with the increasing Reynolds number, as shown in \cite{mishraj1}. See also Fig.~\ref{fig:snapshots_2d}(e): the upper recirculation has enlarged and pushed the high velocity jet to the upper wall. Notice that the flow could have evolved to its mirrored image configuration with respect to the domain symmetry axis.

\end{itemize}

At the minimum value of the Reynolds number for which the asymmetric jet configuration exists there is a \emph{symmetry breaking bifurcation point} or \emph{steady-state bifurcation point}. In Sec.~\ref{2d_results_2param} we will show how the critical value of the
Reynolds number for the symmetry breaking bifurcation changes as $\lambda$ varies. 

\subsection{3D case}
\label{3dcase}
\begin{figure}
\centering
\includegraphics[width=0.5\textwidth]{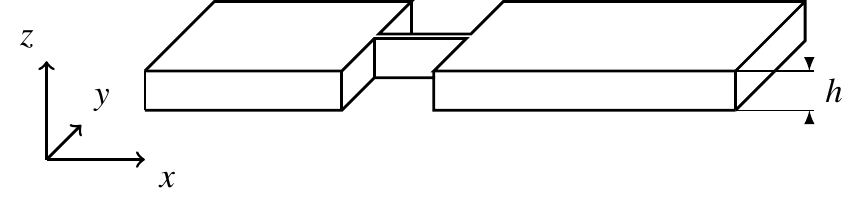}
\caption{Scheme for the 3D geometry.}
\label{fig:scheme_3d}
\end{figure}

A 3D channel is obtained by extruding the 2D geometry considered in the previous section in the direction orthogonal to the flow plane.
The goal of this test case is to study the influence of the channel depth on the flow pattern, and in particular  
on the symmetry breaking bifurcation point. Once $\lambda$ is fixed, the 3D problem depends on two parameters: 
the Reynolds number and the channel depth $h$. See Figure~\ref{fig:scheme_3d}. 
Note that also in this case the geometry reported in the figure has been cropped due to space limitation. 
In reality, we considered a channel length past the expansion equal to 6 times the channel height $L_c$.

The boundary conditions are the same as in the 2D case: we impose a parabolic velocity profile at the inlet, a stress-free condition at the outlet, and a no-slip condition everywhere else.

We introduce the following quantities, which are useful in the characterization of the numerical simulation~\cite{Oliveira}:
\begin{itemize}
\item[-] {\bf aspect ratio}: 
\begin{align}\label{eq:aspect_ratio}
\AR=\frac{h}{w_c};
\end{align}
\item[-] {\bf normalized channel depth}: $\displaystyle{\HH=\frac{h}{h+w_c}=\frac{\AR}{\AR+1}}$;
\item[-] {\bf average horizontal velocity}: $\displaystyle{\langle v_x\rangle=\frac{Q}{w_c h}}$, with $Q$ flow rate;
\item[-] {\bf Reynolds number}: $\displaystyle{\re_{3D}=\frac{\langle v_x\rangle}{\nu} \frac{2 w_c h}{h + w_c} = \frac{\langle v_x \rangle w_c}{\nu}\frac{2\AR}{\AR+1}}$. 
\end{itemize}

Note that $\HH=1$ is the limit case of infinite channel depth, which corresponds to the 2D configuration.
The definition of $\re_{3D}$ has been obtained by setting $U = \langle v_x\rangle$ and $L$ 
equal to the hydraulic diameter of the contraction channel in (\ref{eq:re}). 
We remark that the 2D case can be seen as a limit of the 3D case for $\AR\to\infty$. 
This justifies the factor 2 in the definition of $\re_{2D}$ in (\ref{eq:re_2d}).
Another reason to define $\re_{2D}$ in (\ref{eq:re_2d}) is to compare
our results with \cite{Oliveira} (see Sec.~\ref{2d_results_1param}).

For our 3D tests in Sec.~\ref{3d_results}, we set $\lambda = 15.4$ and let the normalized channel depth
$\HH$ span interval $[0.2,0.95]$. This corresponds to a wide range of aspect ratios: $\AR\in[0.2635,19.71]$.
As for the Reynolds number, we consider the same interval of interest   
used for the 2D case, namely $[0.01,90]$.

%% file: tex/numerical_method_revised.tex
\section{Numerical method}
\label{method}

We are interested in adopting a Reduced Order Model (ROM) for
the numerical solution of the problems presented in the previous section.
Reduced Order Models have been introduced for parametrized problems requiring real-time capabilities due to a many-query setting. 
The goal is to compute reliable results at a fraction of the cost of a conventional (full order) method. A practical way to realize 
this is to organize the computation in two steps:
\begin{itemize}
\item[-] An {\bf offline phase}: full order approximation solutions corresponding to selected representative 
parameters values/system configurations are computed and stored, together with other information 
concerning the parametrized problem.
This is a computationally expensive step usually performed on high performance computing facilities.
\item[-] An {\bf online phase}: the information obtained during the offline phase is used to compute the solution for a newly specified value of the parameters in a short amount of time (ideally in real time), even on a relatively low power device such as a laptop or a smartphone.
\end{itemize}
These split computational procedures are built in such a way that new parameter dependent quantities are easily and quickly computed online, while representative basis functions for selected parameter values and more demanding quantities are pre-computed offline.
We refer to \cite{LMQR:2014} for a review of ROM in Fluid Mechanics.

The problem under consideration might depend on several parameters, with each parameter varying in a certain range.  
We introduce a parameter vector $\mmu$ that contains all parameters. If the problem depends on two parameters, 
we have $\mmu= (\mu_1, \mu_2) \in \mathcal{D}_1 \times \mathcal{D}_2 = \mathcal{D} $; 
for example $\mmu = (\re,\HH)\in \mathcal{D}=[0.01,90]\times[0.2,0.95]$. 
We consider both physical parameters (e.g., $\re$) and geometric parameters (e.g., $\lambda$ and $\HH$). 
To stress the solution dependence on the parameter(s), we will use the notation $\bs{u}=\bs{u}(\mmu)$ and $p = p(\mmu)$, 
without implying that there is a one-to-one correspondence between $\mmu$ and $\bs{u}$ or $p$. 
The maps $\mmu\mapsto\bs{u}(\mmu)$ and $\mmu\mapsto p(\mmu)$ are one-to-one only on the region of $\mathcal{D}$ where there exists a unique solution of problem (\ref{eq:uns_strong1})-(\ref{eq:uns_strong2}).

The treatment of the geometric parametrization deserves further explanation. Let $\xi \in \mathcal{D}$ be a geometric parameter the problem depends on. We select a reference domain $\Omh$ that is mapped to the parametrized domain $\Om(\xi)$ through a one-to-one, orientation preserving transformation $\mathcal{T}:\mathcal{D}\times\Omh\to\Om(\xi)$. 
Using this map, we can cast eq. (\ref{eq:weak_ns1})-(\ref{eq:weak_ns2}) into the reference domain. For instance, eq. (\ref{eq:weak_ns2}) becomes:
\begin{equation}
(q,\ddiv\bs{u})_{\Om(\xi)}=\int_{\Om(\xi)} q\ddiv \bs{u}\de \bs{x}=
\int_{\Omh} q\,\mathcal{F}^{-T}|J(\xi)|\ddiv \bs{u} \de\widehat{\bs{x}},
\label{eq:mapping_div}
\end{equation}
where $\mathcal{F}(\xi)$ is the Jacobian matrix of transformation $\mathcal{T}(\xi)$ and $J(\xi)$ its determinant.

Among many Reduced Order Models available in the literature, we choose a Reduced Basis (RB) method. We will briefly recall the main features of RB methods in the following sections. For a general review on the RB method we refer to, e.g.,~\cite{Rozza:ARCME,hesthaven2015certified,QuarteroniManzoniNegri}. 

\subsection{Full order approximation}
\label{full_order}
As full order approximation scheme for eq. (\ref{eq:weak_ns1})-(\ref{eq:weak_ns2}) to be used in the offline phase, 
we choose the Spectral Element Method (SEM). See, e.g., \cite{Fischer:02,CHQZ1,CHQZ2} for a general review of SEM 
and application to fluid mechanics. We adopt the SEM implementation available in open source software \texttt{Nek5000}~\cite{nek5000}, where the basis functions for each element are the Lagrange interpolants on a Gauss-Lobatto-Legendre tensor product grid. We refer to~\cite{Fischer:02} for an introduction to efficient SEM implementation.

Let $\bs{V}^\mathcal{N}$, $Q^\mathcal{N}$, and $\bs{V}^\mathcal{N}_0$ be the Spectral Element spaces, which are finite dimensional 
approximations of the infinite dimensional spaces $\bs{V}_D$, $Q$, and $\bs{V}$, respectively.
The full order approximation problem reads: for a given $\mmu\in\mathcal{D}$, find $(\bs{u}^\mathcal{N}(\mmu),p^\mathcal{N}(\mmu)) \in \bs{V}^\mathcal{N} \times Q^\mathcal{N}$ such that
\begin{align}
\left( \frac{\partial \bs{u}^\mathcal{N}}{\partial t}, \bs{v} \right)_\Om +
(\bs{v},(\bs{u}^\mathcal{N}\cdot\nabla)\bs{u}^\mathcal{N})_\Om+\nu(\ee(\bs{v}),\ee(\bs{u}^\mathcal{N}))_\Om-(\ddiv\bs{v},p^\mathcal{N})_\Om&=0,\qquad \forall\bs{v}\in\bs{V}^\mathcal{N}_0, \label{eq:NS_full_order1}\\
(q,\ddiv\bs{u}^\mathcal{N})_\Om&=0, \qquad \forall q\in Q^\mathcal{N}  \label{eq:NS_full_order2}.
\end{align}
Notice that in eq. (\ref{eq:NS_full_order1})-(\ref{eq:NS_full_order2}) we have already accounted for the fact that $\bs{g} = \bs{0}$ in (\ref{eq:neu_bc}) for both the 2D and 3D case.

For the computations in Sec.~\ref{sec:res}, we choose the stable $\mathbb{P}_{11}-\mathbb{P}_{9}$ couple for velocity and pressure approximation. In the \texttt{Nek5000} solver, the aliasing errors associated with the choice of high order polynomials for approximating the nonlinearity are dealt with the \emph{3/2 rule} (also called \emph{zero-padding rule}, see~\cite{Boyd}). This rule consists in evaluating the integrals of the nonlinear term (to be liearized) using a quadrature formula with $3/2$ times the quadrature points of the other terms, so that the aliasing errors contribute only for those wavelengths that are filtered out by the grid size.  
For the time discretization of eq. (\ref{eq:NS_full_order1})-(\ref{eq:NS_full_order2}) we adopt 
a Backward Differentiation Formula of order 3 (BDF3; see, e.g., \cite{Fischer:02}). 
The convective term is treated explicitly, with a third order extrapolation formula
as explained in \texttt{Nek5000} documentation~\cite{nek5000doc}. 
Such a treatment of the convective term does not guarantee the unconditional stability in time of the linearized numerical scheme.
A CFL condition 
has to be verified at every collocation point.



Given an initial solution, we consider the system to be close enough to the steady state
when the following stopping condition is satisfied:
\begin{equation}
\frac{\|\bs{u}^\mathcal{N}_n-\bs{u}^\mathcal{N}_{n-1}\|_{L^2(\Om)}}{\|\bs{u}^\mathcal{N}_n\|_{L^2(\Om)}}<\mathtt{tol},
\label{eq:def_increment}
\end{equation}
with tolerance $\mathtt{tol}=10^{-8}$. When the stopping criterion~\eqref{eq:def_increment} is met, the simulation is interrupted.

\subsection{Sampling}

The sampling process consist in selecting $N$ parameters $\{\mmu^i\}$, with $i=1, \dots, N$, in the parameter space $\mathcal{D}$, whose corresponding solutions $\{\bs{u}^\mathcal{N}(\mmu^i)\}\subset\bs{V}^\mathcal{N}$ and $\{p^\mathcal{N}(\mmu^i)\}\subset Q^\mathcal{N}$ will be used to construct the  \emph{Reduced Basis spaces} for velocity and pressure, respectively. 
Solutions $\bs{u}^\mathcal{N}(\mmu^i)$ and $p^\mathcal{N}(\mmu^i)$, with $i=1, \dots, N$, are called \emph{snapshots}. In order to simplify the notation, we will denote $\bs{u}^\mathcal{N}(\mmu^i)$ by $\bs{u}(\mmu^i)$ and $p^\mathcal{N}(\mmu^i)$ by $p(\mmu^i)$.
In this section, we are going to explain how to select the velocity snapshots. The same procedure can be 
applied to obtain the pressure snapshots.

\begin{figure}[h!]
\centering
\includegraphics[width=0.5\textwidth]{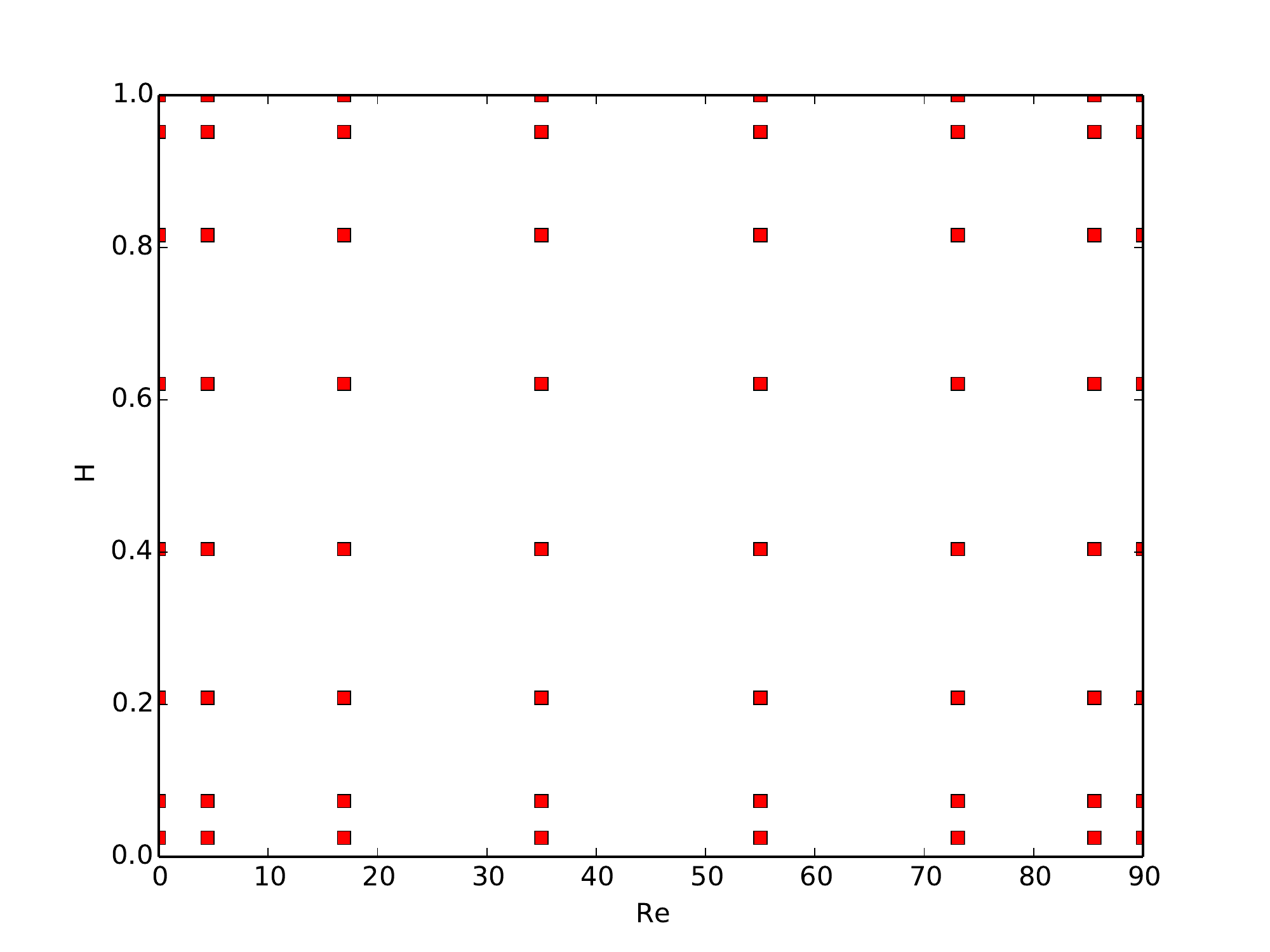}
\caption{Sample parameters for the 3D case with $\mmu = (\re,\HH)\in =[0.01,90]\times[0.025,1]$}
\label{fig:plot_cheb}
\end{figure}

Let $\mu_k$ be the $k$-th component of parameter vector $\mmu$ and let $\mathcal{D}_k$ be the interval of interest for such component. The sampling procedure described below will select $N_k$ values of $\mu_k$ in $\mathcal{D}_k$ and the total number of sample parameter vectors is $N = \prod_k N_k$:
\begin{equation}
\{\mmu^i\}_{i=1}^N=\otimes_{k}\{\mu_k^j\}_{j=1}^{N_k}.
\label{eq:Cheb_Tensor}
\end{equation}
For each component $\mu_{k}$ of the parameter vector $\mmu$ we choose as $\mu_k^j$ the Chebyshev points:
\begin{equation}
\mu_k^j=\mu_{k,\mathrm{min}}+\frac{\mu_{k,\mathrm{max}}-\mu_{k,\mathrm{min}}}{2} \cos((j-1)\pi/(N_k-1)), \quad j = 1, \dots, N_k,
\label{eq:Cheb_1D}
\end{equation}
where
\begin{equation}
\mu_{k,\mathrm{min}}=\min_{\mu_k \in\mathcal{D}_k}\mu_k \qquad \mu_{k,\mathrm{max}}= \max_{\mu_k \in \mathcal{D}_k}\mu_k.
\label{eq:who_are_min_max}
\end{equation}
This procedure to sample sample points in $\mathcal{D}$ is called Gauss-Lobatto-Chebyshev (GLC) tensor product collocation strategy~\cite{XiuHesthaven}. 
For example, Fig.~\ref{fig:plot_cheb} shows the sample parameters considered for the 3D case where $\mmu = (\re,\HH)\in =[0.01,90]\times[0.025,1]$. 
The corresponding values of the Reynolds number and $\HH$ are reported in tables~\ref{tab:which_reynolds} and~\ref{tab:sampling_data}, respectively.
Note that the tensor product collocation allows to choose a different number of sampling points $N_k$ for each component of the parameter space. 

The GLC collocation points are a practical choice, since in this case many sampling points share the same geometric
parameter, requiring to start the continuation method very few times \cite{XiuHesthaven}. Other sampling methods such as Greedy \cite{Rozza:ARCME} or CVT  \cite{PittonRozza} would require to start the continuation method for each new sampling point.
A disadvantage of the GLC collocation strategy is that the approximation spaces are not hierarchical, meaning that the RB spaces for a certain value of $N$ are not in general subspaces of the RB spaces obtained for a higher value of $N$. This could increase
the offline computational cost in case we need to enrich the RB spaces.

Our conjecture is that for bifurcation problems it may be useful to cluster the sampling points close to the bifurcation points. However, to the best of our knowledge, there are no error estimates for different sampling methods for steady state Navier-Stokes equations involving bifurcation points. In any case we are enriching our investigation with an eigenvalue analysis to detect the bifurcation point at the reduced level as well, as we will introduce in Sec.~\ref{sec:bif_detection}.

\begin{table}[h!]
\caption{2D and 3D case: values of the Reynolds number used for the Chebyshev collocation sampling.}
\label{tab:which_reynolds}
\centering
\begin{tabular}[tbp]{lSSSSSSSSS}
\toprule
        &  1  &  2  &  3  &  4  &  5  &  6  &  7  &  8 & 9 \\
\midrule
$\re$   & 0.010& 4.466&13.19&27.79&45.01&62.22&76.82&86.58&90.00 \\
\bottomrule
\end{tabular}
\end{table}

\begin{table}[h!]
\caption{3D case: geometric parameters used for the Chebyshev collocation sampling. Notice that sample 8 corresponds to the 2D case.}
\label{tab:sampling_data}
\centering
\begin{tabular}[tbp]{lSSSSSSSS}
\toprule
        &  1  &  2  &  3  &  4  &  5  &  6  &  7  &  8  \\
\midrule
$\HH$   & 0.025 & 0.0733 & 0.2085 & 0.4040 & 0.6210 & 0.8165 & 0.9517 & 1 \\
$\AR$   & 0.026 & 0.0791 & 0.2634 & 0.6779 & 1.6389 & 4.4550 & 19.704 & $\infty$ \\
\bottomrule
\end{tabular}
\end{table}

\subsection{Construction of the RB spaces}\label{sec:rb_spaces}
For every selected sample $\mmu^i$, we solve the full order approximation problem (\ref{eq:NS_full_order1})-(\ref{eq:NS_full_order2})
until stopping criterion \eqref{eq:def_increment} is satisfied to get $\bs{u}(\mmu^i)$ and $p(\mmu^i)$.
After the sampling is complete, we have two sets of snapshots $\{\bs{u}(\mmu^i)\}_{i=1}^N$ and $\{p(\mmu^i)\}_{i=1}^N$ which
generate the finite dimensional subspaces $\bs{V}^N$ and  $Q^N$, called Reduced Basis spaces.
The key feature of a correct ROM is that the dimension of the reduced order space is much lower that the dimension of the full order space:
\begin{equation}
        N=\dim \bs{V}^N\ll\mathcal{N}=\dim \bs{V}^\mathcal{N}.
        \label{eq:N_vs_Ncal}
\end{equation}
In this way, all the computations required by the online phase (see Sec. \ref{sec:online}) 
will be much less expensive that the computations required during the offline phase.

In this section, we are going to focus on how to construct the Reduced Basis $\{\bs{\phi}_i\}_{i=1}^N$ for $\bs{V}^N$. 
The same procedure can be applied to obtain the Reduced Basis $\{\sigma_i\}_{i=1}^N$ 
for $Q^N$.
To actually build the space $\bs{V}^N$, it is usually preferred not to evaluate directly the inner products in eq.~(\ref{eq:NS_full_order1})-(\ref{eq:NS_full_order2}) for all the snapshots $\{\bs{u}(\mmu^i)\}$, for two reasons:
\begin{itemize}
        \item[-] the snapshots may contain some redundant information, due to the sampling procedure, thus 
        leading to linear dependence and ill-conditioning during the matrix assembling.
      \item[-] in general the snapshots will not be orthogonal to each other, and consequently they will generate a full mass matrix, increasing the storage requirement and operations count.
\end{itemize}
Thus, it is preferred to compute an orthonormal generating set $\{\bs{\phi}_i\}_{i=1}^N$ for $\bs{V}^N$ so that the resulting linearized problem is well conditioned.
Two of the most popular techniques to compute othogonal basis functions are the Proper Orthogonal decomposition (POD) \cite {Volkwein} and the Gram-Schmidt orthogonalization (GS) with its variants \cite{RozzaVeroy}. 

One of the ways to build a POD is to compute the correlation matrix for the set of snapshots, defined as as
\begin{equation}
  \mathbb{C}_{ij}=(\bs{u}^{\mathcal{N}}(\mmu^i),\bs{u}^{\mathcal{N}}(\mmu^j)).
  \label{eq:def_Cij}
\end{equation}
The eigenvalues $\lambda_i$ and eigenvectors $\boldsymbol{\psi}_i$ of $\mathbb{C}$ are computed, and each POD basis vector is defined as
\begin{equation}
  \bs{\phi}_i=\sum_{k=1}^N \psi_{i,k} \bs{u}^{\mathcal{N}}(\mmu^k)
  \label{eq:def_POD_basis}
\end{equation}
where $\psi_{i,k}$ denotes the $k$-th component of the $i$-th eigenvector. The eigenvalue $\lambda_i$ associated to each POD mode is related to the fraction of energy stored in the corresponding mode.
The POD modes are automatically orthogonal in the $L^2$ inner product, but not normal in general. 

It can be shown~\cite{Volkwein,gunzburger_pfco} that the space generated by the POD, denoted by $\bs{V}^N_{\mathrm{POD}}$, minimizes the projection error in the $L^2$ norm:
\begin{equation}
  \bs{V}^N_{\mathrm{POD}}=\argmin_{\bs{V}^N\subset \bs{V}, \dim \bs{V}^N=N} \sum_{i=1}^N\left\|\bs{u}(\mmu^i)-\Proj_{\bs{V}^N}\bs{u}(\mmu^i)\right\|.
  \label{eq:minPOD}
\end{equation}
where $\Proj_{\bs{V}^N}:\bs{V}\to \bs{V}^N$ is the projection operator on the space $\bs{V}^N$ generated by the POD modes.

For the results reported in Sec.~\ref{sec:res}, we have used POD to compute the basis functions. 
However, we would like to remark that POD has one major drawback when the number of snapshots is large: the number of operations required to compute the correlation matrix (\ref{eq:def_Cij}) and its eigenpairs becomes prohibitive. 
If this is the case, a Gram-Schmidt orthogonalization is usually
preferred over computing the basis functions by means of a POD \cite{RozzaVeroy}:
\begin{equation}
  \bs{\phi}_i =\bs{u}^{\mathcal{N}}(\mmu^i)-\sum_{j=1}^i\bs{\phi}_j(\bs{u}^{\mathcal{N}}(\mmu^i),\bs{\phi}_j)\qquad \bs{\phi}_i\leftarrow\frac{\bs{\phi}_i}{\|\bs{\phi}_i\|},
\label{eq:Gram_Schmidt}
\end{equation}
where the normalization step is not always adopted, since it may lead to an ill-conditioned linear system.

As mentioned above, the same orthonomalization method can be applied to obtain the pressure basis $\{\sigma_i\}_{i=1}^N$. 
The reduced basis spaces $\bs{V}^N$ and $Q^N$ are defined as:
\begin{equation}\label{eq:VN_QN}
  \bs{V}^N=\Span\{\bs{\phi}_i\}_{i=1}^{N}\quad \text{and} \quad
  Q^N=\Span\{\sigma_k\}_{k=1}^{N}.
\end{equation}

\subsection{Online phase computation}\label{sec:online}
After the construction of the Reduced Basis spaces, an \emph{online approximation} 
of the solution can be computed by applying the Galerkin projection to the spaces $\bs{V}^N$ 
and $Q^N$. Namely, given a target parameter $\mmu\in\mathcal{D}$
we search for $(\bs{u}^N(\mmu),p^N(\mmu))\in\bs{V}^N\times\bs{Q}^N$ such that

\begin{align}
\left( \frac{\partial \bs{u}^N(\mmu)}{\partial t}, \bs{v} \right)_\Om + 
\left(\bs{v},(\bs{u}^N(\mmu)\cdot\nabla)\bs{u}^N(\mmu)\right)_\Om+\nu\left(\ee(\bs{v}),\ee(\bs{u}^N(\mmu))\right)_\Om-\left(\nabla \cdot \bs{v},p^N(\mmu)\right)_\Om&= 0 \qquad \forall\bs{v}\in\bs{V}^N,\label{eq:NS_variational_RB1} \\
\left(q,\ddiv \bs{u}^N(\mmu)\right)_\Om&=0 \hspace{.7cm}  \forall q\in Q^N. \label{eq:NS_variational_RB2}
\end{align}
For convenience, 
for the rest of this section we are going to assume that the first term in eq.~\eqref{eq:NS_variational_RB1}
is negligible, as if the flow was modeled by eq.~(\ref{eq:ns_strong1})-(\ref{eq:ns_strong2}).

The convective term is linearized with a fixed point scheme. 
Suppose that an initial tentative solution $\bs{u}^N_0(\mmu)$ is known. Given $\bs{u}^N_{k-1}(\mmu)$, at the $k$-th iteration of the 
fixed point method we solve problem:

\begin{align}
 \left(\bs{v},(\bs{u}^N_{k-1}(\mmu)\cdot\nabla)\bs{u}^N_k(\mmu)\right)_\Om+\nu\left(\ee(\bs{v}),\ee(\bs{u}^N_k(\mmu))\right)_\Om-\left(\nabla \cdot \bs{v},p^N_k(\mmu)\right)_\Om&= 0 \qquad \forall\bs{v}\in\bs{V}^N,\label{eq:NS_RB_fixed_point1} \\
\left(q,\ddiv \bs{u}^N_k(\mmu)\right)_\Om&=0 \hspace{.7cm}  \forall q\in Q^N. \label{eq:NS_RB_fixed_point2}
\end{align}
The iterative scheme can be stopped for instance when an increment-based residual:
  \begin{equation}
  \texttt{res}_k=\frac{\|\bs{u}^N_k(\mmu)-\bs{u}^N_{k-1}(\mmu)\|}{\|\bs{u}^N_k(\mmu)\|}
  \label{eq:def_increment_res}
  \end{equation}
is below a given tolerance.

The solution scheme described so far requires the $\bs{V}^N-Q^N$ pair to satisfy a stability condition 
called \emph{inf-sup condition} or \emph{Ladyzhenskaya-Brezzi-Babu\v{s}ka (LBB) condition}: 
\begin{equation}
\inf_{q\in Q^N}\sup_{\bs{v}\in\bs{V}^N}\frac{(q,\ddiv\bs{v})}{\|q\|_Q\|\bs{v}\|_{\bs{V}}}=\beta^N>0.
\label{eq:inf_sup}
\end{equation}
See, e.g., \cite{B-ladyzhenskaya,babuska1,brezzi1,Brezzi:mixed,Ern_Guermond}.
Spaces $\bs{V}^N$ and $Q^N$ in (\ref{eq:VN_QN}) computed using the POD or GS modes 
as explained in Sec.~\ref{sec:rb_spaces} are not guaranteed to fulfill condition~\eqref{eq:inf_sup}. 
There are two options for circumventing this issue: casting the problem into a divergence-free space 
(see, e.g., \cite{PittonRozza}) and enforcing approximation stability properties for the 
$\bs{V}^N-Q^N$ pair (see, e.g., \cite{RozzaVeroy,RZinfsup,BallarinManzoniQuarteroniRozza2014}). 
Here, we choose the former approach.
This means that we require $\bs{V}^N$ to be a subset of $H^1_{\ddiv}(\Om)$:
\begin{equation}
H^1_{\ddiv}(\Om):=\{\bs{v} \in H^1(\Om)\text{ s.t. } (q,\ddiv \bs{v})=0\,\,\forall q\in L^2(\Om)\},
\label{eq:whois_H1div}
\end{equation}
which is a subspace of $H^1(\Om)$.
If the basis functions for $\bs{V}^N$ are divergence-free, eq. (\ref{eq:NS_RB_fixed_point2}) is no longer needed. 
Thus, the pressure disappears from the variational formulation and we do not need to build the space $Q^N$.
See, e.g.,~\cite{FoiasTemam}.

With a divergence-free basis set for $\bs{V}^N$, at every fixed-point iteration we have to solve the following linear system:
\begin{equation}
\mathtt{A}^k(\mmu)\mathtt{u}_k=\mathtt{b}_k
\label{eq:fixed_point_algebra}
\end{equation}
where $\mathtt{u}_k$ is the vector containing the projection coefficients of $\bs{u}^N_{k-1}$ onto the space 
$\bs{V}^N$, $\mathtt{b}_k\in\mathbb{R}^N$ depends from the specified boundary conditions, and $\mathtt{A}_k(\mmu)\in\mathbb{R}^{N\times N}$ is given by:
\begin{equation}
\mathtt{A}^k_{lj}(\mmu)= \left(\bs{\phi}_l,\bs{u}^N_{k-1}(\mmu)\cdot\nabla\bs{\phi}_j\right)_\Om
+\nu \left(\nabla\bs{\phi}_l,\nabla\bs{\phi}_j \right)_\Om.
\label{eq:whois_A}
\end{equation}
Once the velocity $\bs{u}^N(\mmu) \in H^1_{\ddiv}(\Om)$ has been computed, the pressure can be recovered, for example, 
by solving a Poisson problem online:
\begin{align}
\Delta p^N(\mmu) = - \ddiv \left(\bs{u}^N(\mmu) \cdot \nabla \bs{u}^N(\mmu) \right). \el
\end{align}
We refer to,e.g., \cite{Caiazzo} for an analysis of velocity-pressure reduced order models.

In equation~\eqref{eq:whois_A}, we wrote explicitly the dependence of matrix $\mathtt{A}$ on the parameter vector $\mmu$. 
Such dependence is more or less evident for the different type of parameters. For instance, 
if the Reynolds number is the only parameter, i.e. $\mmu = \mu_1 = \re $, 
from (\ref{eq:re}) we have $\nu = UL/\re$ and
matrix $\mathtt{A}$ can be written as:
\begin{align}
\mathtt{A}^k_{lj}(\re)= \left(\bs{\phi}_l,\bs{u}^N_{k-1}(\re)\cdot\nabla\bs{\phi}_j\right)_\Om
+ \frac{LU}{\re} \left(\nabla\bs{\phi}_l,\nabla\bs{\phi}_j \right)_\Om, \el
\end{align}
with a linear dependence on $\re^{-1}$.
On the other hand, if $\HH$ is the only parameter,  i.e. $\mmu = \mu_1 = \HH $, the 
dependence of $\mathtt{A}$ on it is hidden 
in the inner products and differential operators. This holds true in general for
geometric parameters. 
Let $\xi \in \mathcal{D}$ be a geometric parameter.  If 
the geometric transformation $\mathcal{T}:\mathcal{D} \times \Omh\to\Om$ is affine, it is possible to express the inner 
products as a linear combination of the inner products on the reference domain:
\begin{equation}
\mathtt{A}(\xi)=\sum_{i=1}^{\dim\mathcal{D}}\Theta^i(\xi)\mathtt{A}^i.
\label{eq:affine_decomposition_A}
\end{equation}
Only functions $\Theta^i$ depend on $\xi$ and need to be evaluated online. 
Matrices $\mathtt{A}^i$ are assembled offline since they do not depend on $\xi$. 
Thus, the affine decomposition (\ref{eq:affine_decomposition_A}) allows for important 
computational time savings. 
In this work, we will consider only affine decompositions. 
If $\mathtt{A}$ depends nonlinearly on $\xi$, 
it has to be computed from scratch for each value of $\xi$.
The efficient assembling of $\mathtt{A}$ when the geometric
transformation is non-affine is still an active research area, 
one of the most popular techniques being 
the \emph{Empirical Interpolation Method}~\cite{EIM}. 

The construction of a divergence-free basis set for $\bs{V}^N$ when
geometric parameters are considered is less trivial than in the case
of physical parameters only. Thus, it requires further explanation.
The \emph{Piola transformation} $\mathcal{P}$ can be seen as the composition of the map $\mathcal{T}$ in eq.~\eqref{eq:mapping_div} with any function $f$ defined on the image (or preimage) of $\mathcal{T}$.
For example, if $f:\widehat{\Om}\to\mathbb{R}$, a new function $g:\Om\to\mathbb{R}$ can be obtained by considering $g(\bs{x})=f(\mathcal{T}^{-1}\bs{x})$ for $\bs{x}\in\Om$.
The Piola transformation $\mathcal{P}$ acts as a map between finite dimensional Hilbert spaces 
$\mathcal{D}\times\bs{V}^N(\Omh)$ and $\bs{V}^N(\Om)$: 
$$
\mathcal{P}:\mathcal{D}\times\bs{V}^N(\Omh)\to\bs{V}^N(\Om).
$$
Its use in an offline-online setting is as follows:
\begin{enumerate}
\item The snapshots $\{\bs{u}^\mathcal{N}(\mmu^i)\}_{i=1}^N$ are divergence-free on the original domain $\Om$. By pulling back the divergence operator to the reference domain $\Omh$ through the Piola map, we obtain a set of snapshots that are divergence free on the reference domain.
\item Perform the POD or GS orthogonalization for the divergence-free snapshots on the reference domain to obtain a basis for $\bs{V}^N(\Omh)$. These basis functions are divergence free on $\Omh$, but not on $\Om$ unless mapped with the Piola transformation.
\item Compute the matrices $\mathtt{A}^i$ in (\ref{eq:affine_decomposition_A}) on the reference domain and with the orthogonal divergence-free basis set.
\item During the online phase, apply the Piola transformation to the matrices $\mathtt{A}^i$ computed at step 3 
so that their entries coincides with the Piola-transformed divergence-free basis functions computed on $\Om$.
\end{enumerate}
We refer to~\cite{Brezzi:mixed} for details on the Piola transformation, and to~\cite{Lovgren} for an application to RB methods in incompressible fluid mechanics in laminar regime. For the application of 
RB methods to moderately turbulent flows we refer for example to \cite{LCLR}, and references therein.

Regarding the boundary conditions, the global support of the RB modes does not allow to impose pointwise values for the non-homogeneous Dirichlet condition. An equivalent way to impose the desired flow conditions is to impose the mass flow rate, instead of the inflow velocity profile. The physically correct definition of mass flow rate is: 
\begin{equation}
  \dot{v}_x=\frac{\int_{\Om} u_x\de\bs{x}}{\int_{\Om}\de\bs{x}},
\label{eq:def_average_mfr}
\end{equation}
where $u_x$ is the $x$-component of the velocity.
Notice that due to incompressibility and the prescribed boundary conditions, (\ref{eq:def_average_mfr}) is equivalent to:
\begin{equation*}
  \dot{v}_x=\frac{\int_{\Ga_{\text{inlet}}}u_x\de\bs{x}}{\int_{\Ga_{\text{inlet}}}\de\bs{x}}.
  \label{eq:def_inlet_mass_flow_rate}
\end{equation*}

From the implementation point of view, it is more convenient to impose the integrated mass flow rate $\dot{w}_x$, that for a given inlet velocity profile is defined as:
\begin{equation}
\dot{w}_x=\int_{\Om}u_x\de\bs{x}.
\label{eq:def_mass_flow_rate}
\end{equation}
We impose the average mass flow condition for the RB simulation through a Lagrange multiplier approach
as follows.
We compute the integrated mass flow rate for each of the RB functions:
\begin{equation}
\dot{c_i}=\int_{\Om} \phi_{i,x}\de\bs{x},
\label{eq:def_average_phi_x}
\end{equation}
and collect all the $\dot{c_i}$ in a vector $\mathtt{C}\in\mathbb{R}^N$.  Let $\alpha_k \in \mathbb{R}$ be the
Lagrange multiplier associated with the mass flow rate constraint at the $k$-th iteration of the fixed point method
described above. Notice that this is a new unknown in the problem.
Then, instead of solving system (\ref{eq:fixed_point_algebra}), at each fixed-point iteration we solve the 
following linear system:
\begin{equation}
\begin{bmatrix}
\mathtt{A}^k(\mmu)  & \mathtt{C}^T \\
\mathtt{C} & \mathtt{0}
\end{bmatrix}
\begin{pmatrix}
\mathtt{u}_k \\
\alpha_k
\end{pmatrix}
=
\begin{pmatrix}
\mathtt{b}_k \\
\dot{w}_x
\end{pmatrix}.
\label{eq:def_lagrange_matrix}
\end{equation}
We remark that imposing a constrainted condition by a Lagrange multiplier is fairly common
in the Reduced Basis context, see, e.g.,\cite{NVP}.

\subsection{Bifurcation detection}\label{sec:bif_detection}

In the configuration described in Section~\ref{setting}, the first pitchfork bifurcation point is determined by a classical modal stability analysis, that can be set up as follows.
Let us consider the 3D case, for which the parameters are $\re_{3D}$ and $\HH$.
Suppose that an initial RB solution $\bs{u}^N(\mmu_i)$ of the steady state problem is known 
for a given value of the parameter $\mmu_i = (\re_{3D,i}, \HH_i)$, 
characterized by a sufficiently small Reynolds number $\re_{3D,i}$ so that the solution is surely unique. 
We proceed as follows: set $s = 1$ and $\re_s = \re_{3D,i}$, then:
\begin{enumerate}
  \item Keeping fixed the value of the geometric parameter $\HH_i$, increase the value of the Reynolds number by a sufficiently small increment $\Delta\re$ (i.e., small enough so that the corrector step will converge to a solution in the desired branch) and set $\mmu_{s+1} = (\re_{s}+\Delta\re, \HH_i)$. 
\item Compute the RB solution $\bs{u}^N(\mmu_{s+1})$ of the steady state problem for the new parameter value $\mmu_{s+1}$. 
\item Compute the Galerkin projection of operator $\mathcal{L}$ defined in (\ref{eq:def_L_lin}) on the RB space $\bs{V}^N$ to form the matrix $\mathtt{L}(\mmu_{s+1})$:
\begin{equation}
\mathtt{L}_{kl}(\mmu_{s+1})=(\bs{\phi}_k,\mathcal{L}(\bs{u}^N(\mmu_{s+1}))[\bs{\phi}_l]).
\label{eq:def_L_matrix}
\end{equation}
\item Compute the eigenvalues 
of $\mathtt{L}(\mmu_{s+1})$ and check if there is one eigenvalue that has changed sign with respect to the previous iteration. If not, set $s = s +1$ and go back to step 1.
\end{enumerate}

We remark that the above algorithm may be unstable in the sense that in a neighborhood of the bifurcation point it may abruptly switch the approximated solution branch, or fail to converge. To make sure that the approximation is always laying on the correct branch a continuation method may be used.

Continuation methods rely on a predictor-corrector iteration to compute solutions lying on the same branch. Suppose that the $\bs{u}(\mmu_s)$ is a solution of equations~\eqref{eq:weak_ns1}-\eqref{eq:weak_ns2}, and is known to lie on a certain branch of interest. The predictor step consists in the computation of an initial guess for the velocity increment $\mathbf{\Delta u}$ due to 
an increase of a single parameter, denoted by $\bs{\Delta \mu}$, by solving the linearized Navier-Stokes equations.
Then, starting from the prediction $\tilde{\bs{u}}=\bs{u}(\mmu_s)+\bs{\Delta u}$, a Newton iteration is set up to impose that the new solution $\bs{u}(\mmu_s+\bs{\Delta\mu})$ solves the original problem~\eqref{eq:weak_ns1}-\eqref{eq:weak_ns2}, under the constraint that the solution be orthogonal to the tangent plane at the point $(\mmu_s,\bs{u}(\mmu_s))$ in the parameter-solution space. We refer to e.g.~\cite{Dijkstra} for an introduction to continuation methods in fluid mechanics.

The continuation method is computationally quite expensive.
Thanks to the fact the the GLC collocation strategy keeps the number of sample values small,
the number of times the continuation method has to be restarted is reduced, allowing for important
computational time savings. 

Note that the matrix $\mathtt{L}$ is dense but has rank equal to $N$, with $N$ of the order of a few tenths at most. Hence all the eigenvalues can be computed inexpensively with QR iterations~\cite{Golub}, for instance. If the spectrum analysis had to be carried out on the full-order model, only a few of the eigenvalues closer to zero could be computed. Moreover, the computations would be much more expensive, requiring Krylov subspace methods~\cite{SaadEig} and most likely a supercomputer.

For the 2D case, we use an analogous algorithm, the only difference being that the parameters are $\re_{2D}$ and $\lambda$.


Lately, increasing attention has been devoted to eigenvalue calculation (as bifurcation detector tool) at the reduced order level \cite{MehrmannSchroder,MehrmannSchroder2}. 
We refer to~\cite{Galdi:NSanalysis} for a theoretical analysis of bifurcation detection techniques in Navier-Stokes equations 
and to~\cite{cliffe} for a bifurcation detection method in a similar geometry.

%% file: tex/results_revised.tex
\section{Results}\label{sec:res}


In this section the method described in section~\ref{method} will be validated against benchmark problems reported in
\cite{Oliveira,Drikakis}. We start with the test cases in two dimensions and then consider problems in three dimensions.
We show that our RB method successfully captures the bifurcation points reported in \cite{Oliveira,Drikakis}. 
We compare our results with full order solutions and provide
an estimate of the computational savings.
Moreover, we carry out an extensive set of simulations that will allow us to 
confirm that  the eccentric mitral regurgitant jets occur when the regurgitant orifice is long 
(large aspect ratio) and narrow (large expansion ratio).

\subsection{2D case: one parameter study}\label{2d_results_1param}

We start with the validation of the bifurcation detection method presented in section~\ref{sec:bif_detection} for the 2D test case
with the Reynolds number as the only varying parameter. For the moment, the geometry is kept fixed. We set the expansion ratio $\lambda$ to $15.4$ in order to compare our results with those reported in reference \cite{Oliveira}. In this case we choose a mesh with 308 spectral elements of order 11, with careful refinement near the re-entrant corners of the domain, where we can expect a loss of regularity for the solution.

As shown in table~\ref{tab:which_reynolds}, we sample nine values for the Reynolds number in the interval $\re_{2D}\in[0.01,90]$.
For the first four values of $\re_{2D}$ in table~\ref{tab:which_reynolds} 
the offline solver returned only the symmetric solution, as expected. 
For the remaining five values, the solver returned two snapshots: one for the symmetric solution (unstable \cite{sobeyd1}) and one for the asymmetric solution (stable). 
As mentioned in Sec.~\ref{2dcase}, at a Reynolds number higher than the critical value
for the symmetry breaking two stable solutions co-exist, which are
one the mirrored image of the other with respect to the horizontal axis (see, e.g., \cite{battagliat1}).
 Bifurcation theory allows to clarify the nature of the multiplicity of possible flows, whereas a (numerical or laboratory) experiment will give one or the other of the stable symmetric solutions.
 Thus, for the multi-parameter case we will disregard the symmetric unstable solution and retain only the stable solutions.

\begin{figure}[htbp]
  \centering
  \subfigure[$\re_{2D} = 0.01$]{
    \includegraphics[width=0.33\columnwidth]{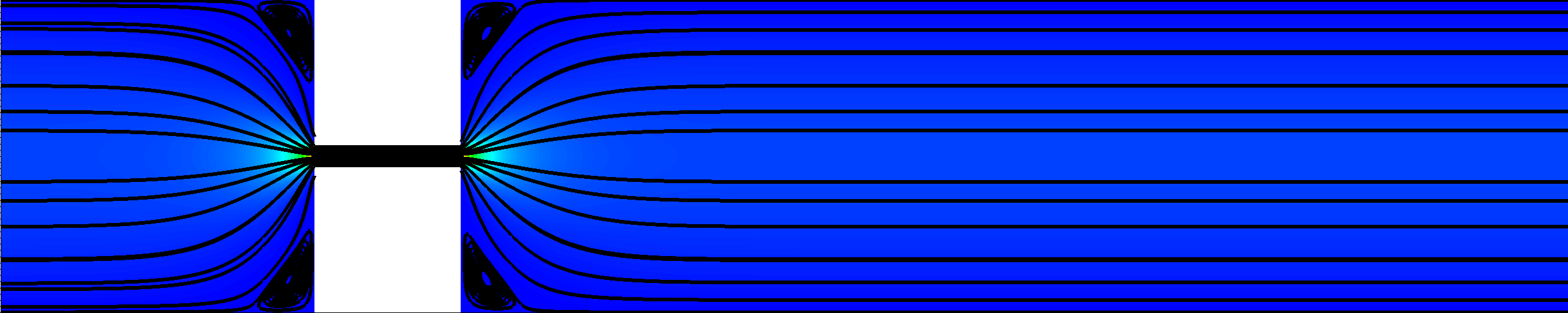}
    }
     \subfigure[$\re_{2D} = 13.19$]{
    \includegraphics[width=0.33\columnwidth]{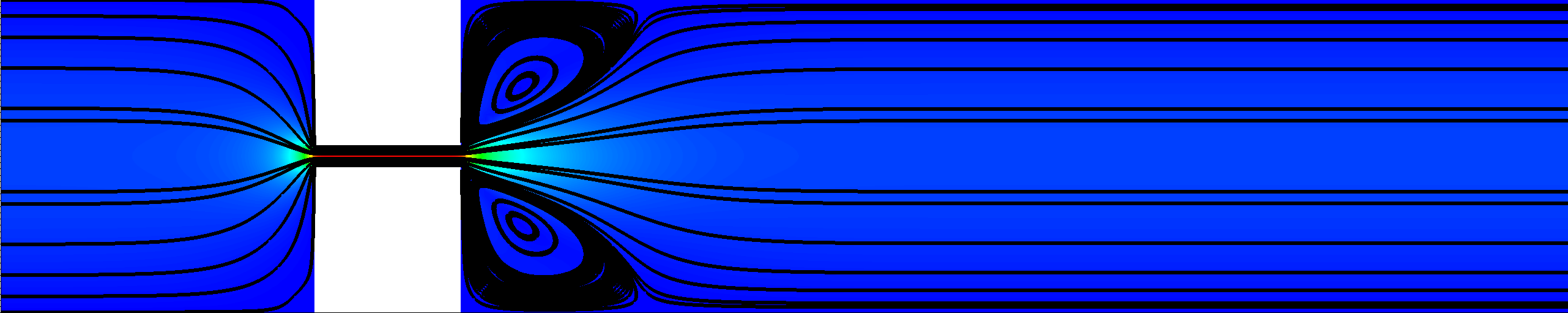}
    }
     \subfigure[$\re_{2D} = 27.7$, unstable solution]{
    \includegraphics[width=0.33\columnwidth]{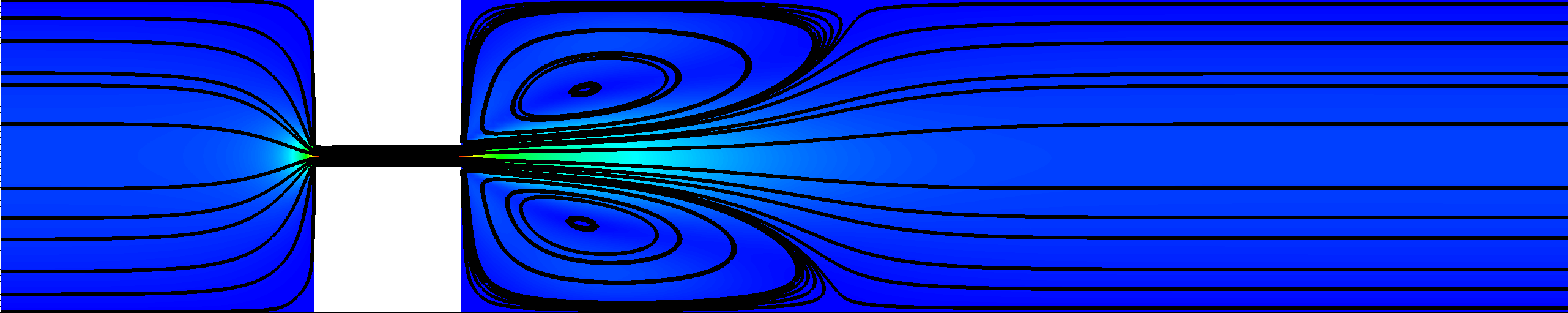}
    }
     \subfigure[$\re_{2D} = 27.7$, stable solution]{
    \includegraphics[width=0.33\columnwidth]{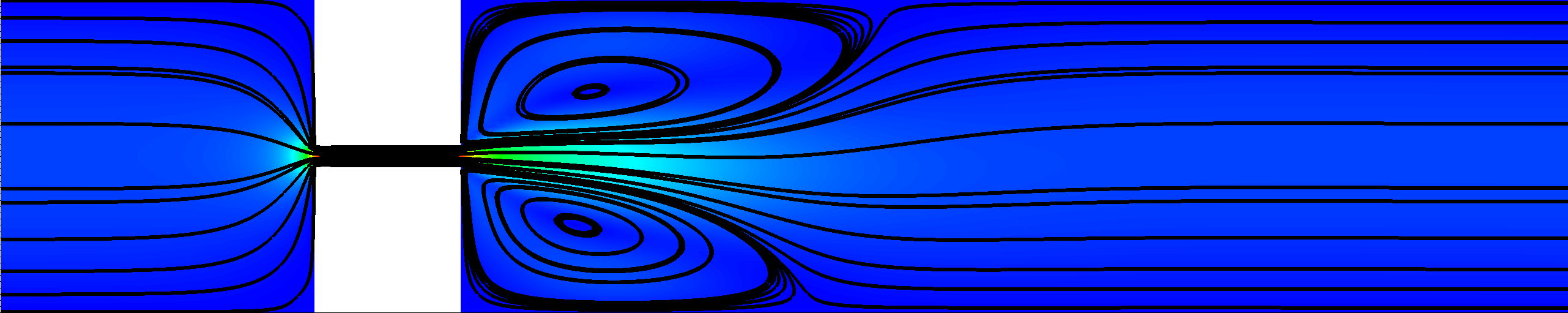}
    }
     \subfigure[$\re_{2D} = 62.22$]{
    \includegraphics[width=0.33\columnwidth]{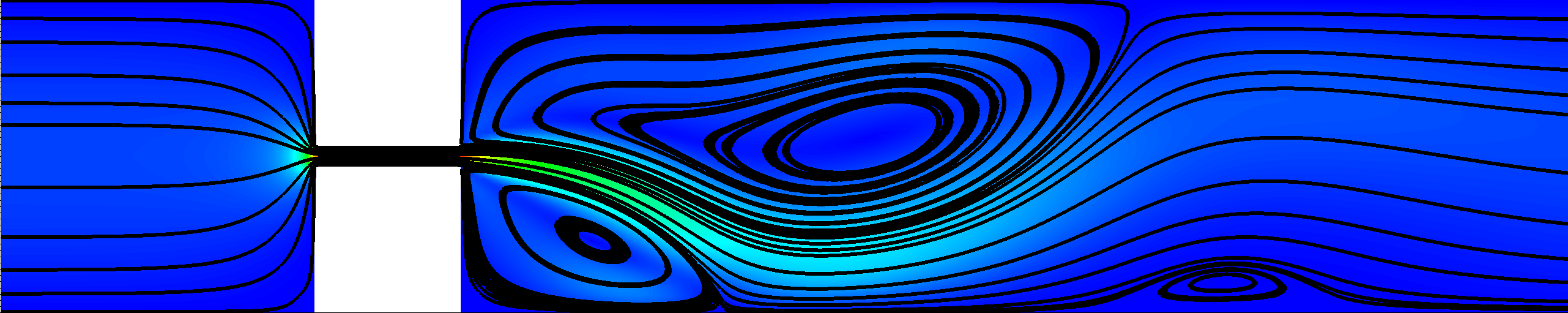}
    }
  \caption{Representative snapshots for the 2D case for $\lambda = 15.4$: velocity magnitude and streamlines for (a) $\re_{2D} = 0.01$, 
  (b) $\re_{2D} = 13.2$, (c) $\re_{2D} = 27.7$ unstable solution, (d) $\re_{2D} = 27.7$ stable solution, and (e) $\re_{2D} =62.2$.}
  \label{fig:snapshots_2d}
\end{figure}

 The online phase for the 2D problem is performed with a RB space of dimension $N=9$. 
In Fig.~\ref{fig:snapshots_2d} we report representative snapshots for the 2D case, corresponding to Reynolds numbers
$\re_{2D} = 0.01, 13.2, 27.7, 62.2$. For very low Reynolds number the solution is characterized 
by symmetry about the horizontal axis and a vertical axis, with a couples of vortices both 
upstream and downstream of the contraction called \emph{Moffatt eddies} \cite{moffatt1}. 
See Fig.~\ref{fig:snapshots_2d}(a). As the inertial effects of fluid become more important (i.e., as $\re_{2D}$ increases), the Moffatt eddies upstream of the contraction gradually diminish in size and two recirculation regions 
of equal size develop downstream of the expansion. 
See Fig.~\ref{fig:snapshots_2d}(b). Symmetry about the vertical axis is lost, but the solution is still symmetric about the horizontal axis.
Past the bifurcation point we can see two solutions: a symmetric one (unstable) and a slightly asymmetric one (stable). 
See Fig.~\ref{fig:snapshots_2d}(c) and (d).
The formation of stable asymmetric vortices in 2D planar expansion is attributed to the \emph{Coanda effect} 
(see \cite{willef1}): an increase in velocity near one wall will lead to a decrease in pressure near that wall and once a pressure difference is established across the channel it will maintain the asymmetry of the flow.
This asymmetric solution remains stable for a certain range of $Re_{2D}$ and asymmetries become stronger with the increasing Reynolds number. See Fig.~\ref{fig:snapshots_2d}(e). 

\begin{figure}[h!]
  \centering
    \includegraphics[width=0.45\columnwidth]{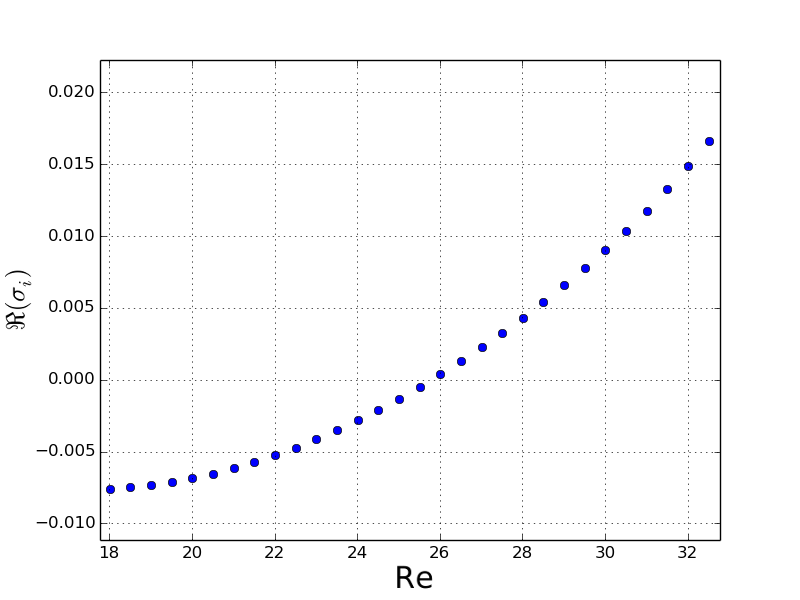}
  \caption{2D case for $\lambda = 15.4$: real part of the eigenvalue
    of matrix $\mathtt{L}$ in (\ref{eq:def_L_matrix}) responsible for the symmetry breaking as a function of the Reynolds number in a neighborhood of a bifurcation point. 
  }
  \label{fig:eig_2d}
\end{figure}

To test the bifurcation detection method described in section~\ref{sec:bif_detection}, we run the online solver parametrized using the 2D basis set with $N=9$ snapshots. In Fig.~\ref{fig:eig_2d}, we plot the real part of the eigenvalue
of matrix $\mathtt{L}$ in (\ref{eq:def_L_matrix}) responsible for the symmetry breaking. We see that the curve crosses the horizontal axis at a Reynolds number of about $\re_{2D,\mathrm{sb}}=26$. This is in good agreement with the critical values for the symmetry breaking reported by \cite{Oliveira,AQpreprint} ($\re_{2D,\mathrm{sb}}=28$) and \cite{mishraj1} ($\re_{2D,\mathrm{sb}}=27$).

\begin{figure}[h!]
  \centering
  \subfigure[eigenvalues of matrix $\mathtt{L}$]{
    \includegraphics[width=0.45\columnwidth]{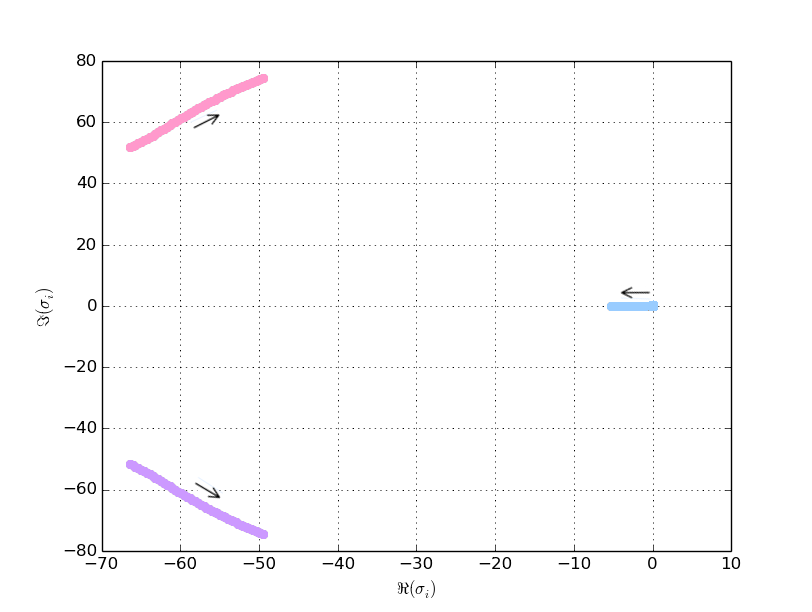}
  }
  \subfigure[zoomed-in view of (a)]{
    \includegraphics[width=0.45\columnwidth]{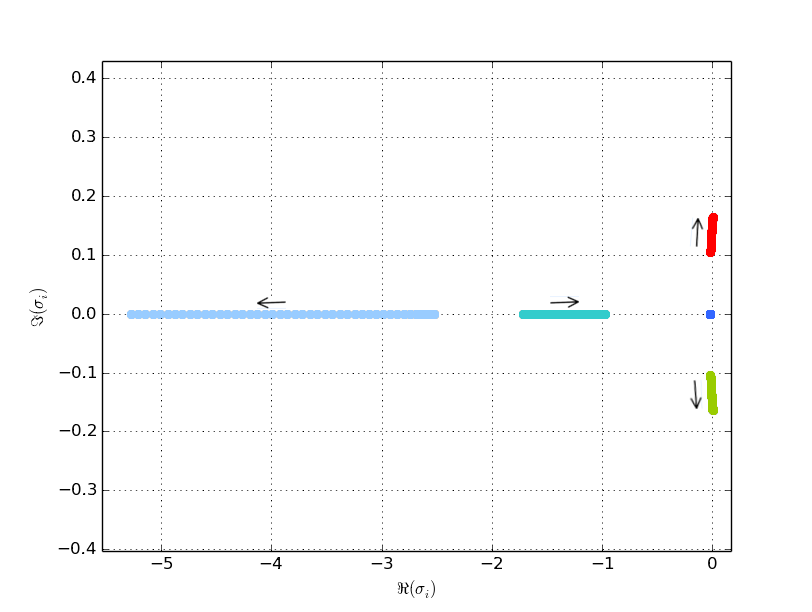}
  }
  \subfigure[zoomed-in view of (b)]{
    \includegraphics[width=0.45\columnwidth]{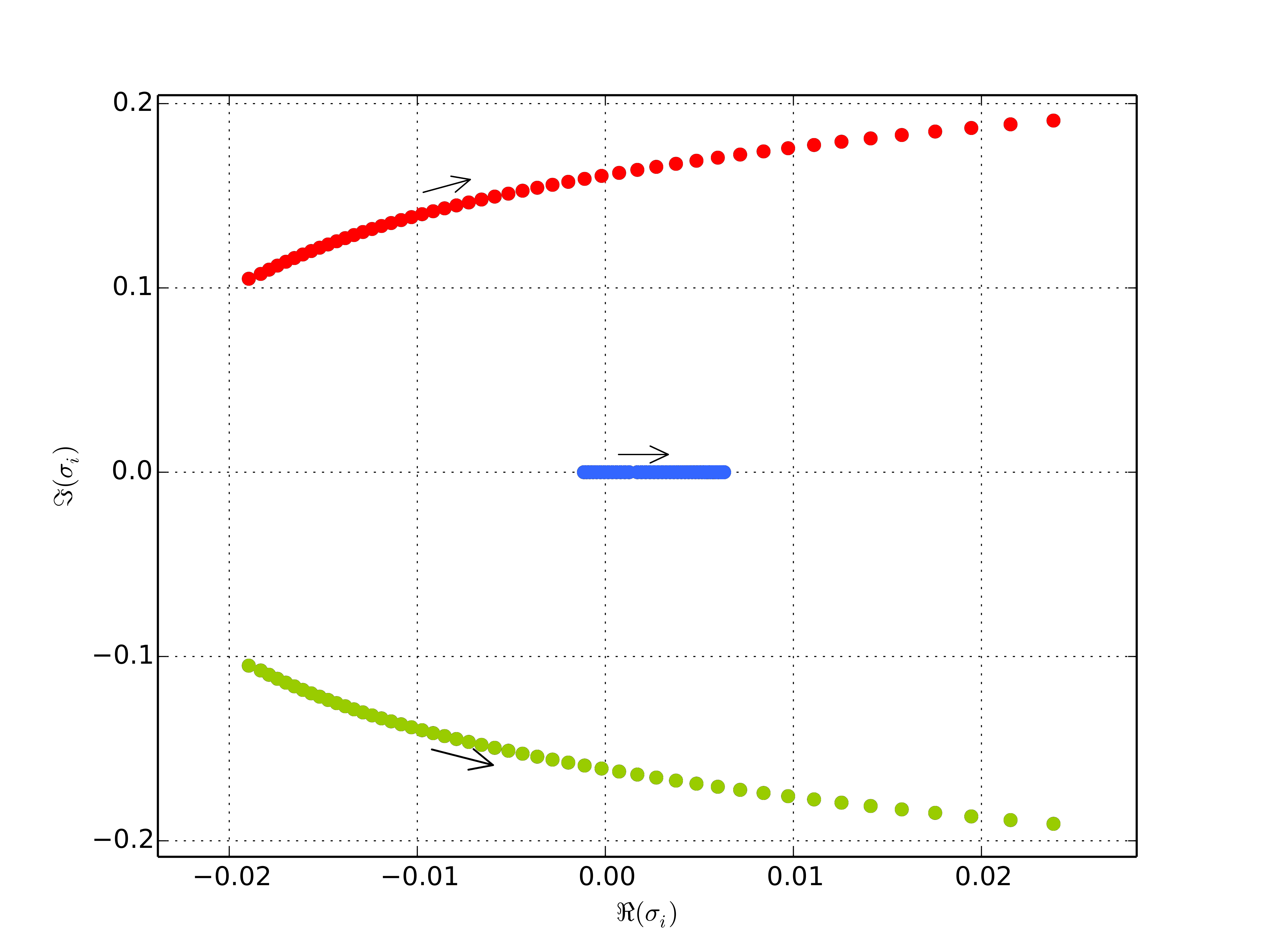}
  }
    \caption{2D case for $\lambda = 15.4$: (a) path of the eigenvalues of matrix $\mathtt{L}$ in (\ref{eq:def_L_matrix})  
    in the complex plane for $\re_{2D}\in[20,55]$.
    Subfigure (b) is a zoomed-in views of subfigure (a) and 
    subfigure (c) is a  zoomed-in views of subfigure (b). The arrows indicate the direction of the increasing Reynolds numbers. 
    The eigenvalue responsible for the bifurcation is the one colored in blue in (c).
     } 
    \label{fig:eig_2d_color}
\end{figure}

Fig.~\ref{fig:eig_2d_color}(a) shows the path of the eigenvalues of matrix $\mathtt{L}$ in (\ref{eq:def_L_matrix})  
in the complex plane for $\re_{2D}\in[20,55]$.
The arrows indicate the direction of the increasing Reynolds numbers. 
Fig.~\ref{fig:eig_2d_color}(b) is a zoomed-in view of Fig.~\ref{fig:eig_2d_color}(a), and
Fig.~\ref{fig:eig_2d_color}(c) is in turn a zoomed-in view of Fig.~\ref{fig:eig_2d_color}(b).
In Fig.~\ref{fig:eig_2d_color}(c) we see the eigenvalue responsible for the bifurcation: 
it is the simple eigenvalue colored in blue that changes sign as the Reynolds number increases.

For this 2D case, the computational savings are significant.
The detection of the bifurcation point using the continuation method required about 80 runs, with a total computational time of around 5 minutes ($0.08h$) on a common desktop computer, which means $3.75$s per online single run. 
The same computations using the full order model described in sec.~\ref{full_order} would have required about 10 CPU-hours per run. 
Hence, adding to the online cost the time required for the RB spaces generation (i.e., the 2 CPU hours required by the POD computations), we can estimate that the  
computational cost for the reduced model is around 11.5\% of the computational cost for the full order model, considering all the operations needed for the bifurcation detection and computation ($N=9$).
\[
\frac{\text{Time to build the RB spaces}+\text{Online time to detect the bifurcation point}}{\text{Time of the equivalent full order computation}}=\frac{9\cdot10h+2h+0.08h}{80\cdot10h}\simeq 11.5\%.
\]
More generally, if only the online runtimes are considered, the computational savings become much more relevant compared with the offline runtimes per single query:
\[
  \frac{\text{RB online query time}}{\text{Equivalent full order single computation}}=\frac{3.75s}{10h}\simeq 0.01\%.
\]

An important quantity to be used as indication if a reduced computational model is  competitive is the \emph{break-even}, comparing all the offline computational times needed to prepare the reduced basis problem ($N=9$) and  an online query with full order model:
\[
  \frac{\text{All full order computations for RB prep.}}{\text{Full order one query  comp. time}}=\frac{9\cdot 10h +2 h}{10h}\simeq 9.2.
\]
suggesting that the use of Reduced Order Methods becomes more and more competitive as the number of queries increases (with 10 or more queries this approach brings already  important computational advantages). Also, this confirms that for the one-parameter scenario, this method could be efficiently adapted to a real-time query tool to be used, e.g., on smartphones or other mobile devices with appropriate \emph{apps}.

\begin{figure}[h!]
  \centering
  \subfigure[$\re_{2D} = 20$, full order]{
    \includegraphics[width=0.45\columnwidth]{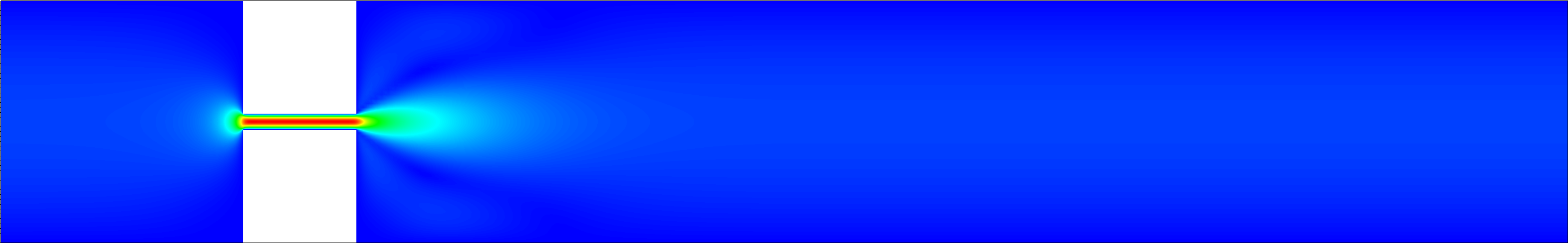}
    }
     \subfigure[$\re_{2D} = 20$, reduced order]{
    \includegraphics[width=0.45\columnwidth]{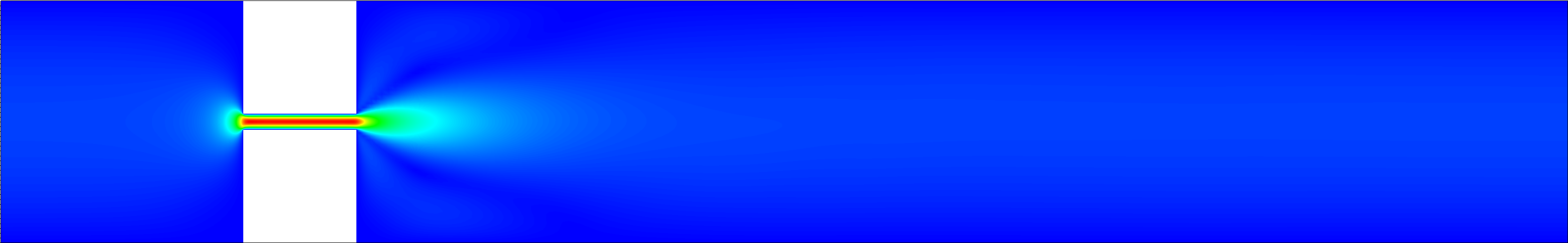}
    }
  \subfigure[$\re_{2D} = 55$, full order]{
    \includegraphics[width=0.45\columnwidth]{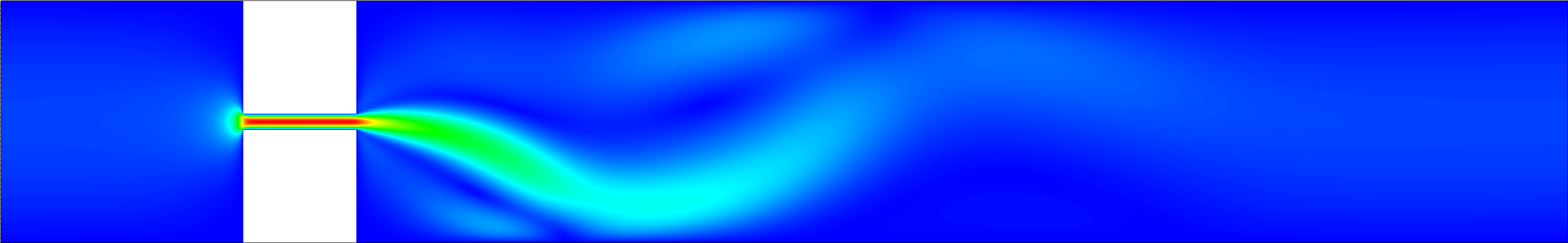}
    }
     \subfigure[$\re_{2D} = 55$, reduced order]{
    \includegraphics[width=0.45\columnwidth]{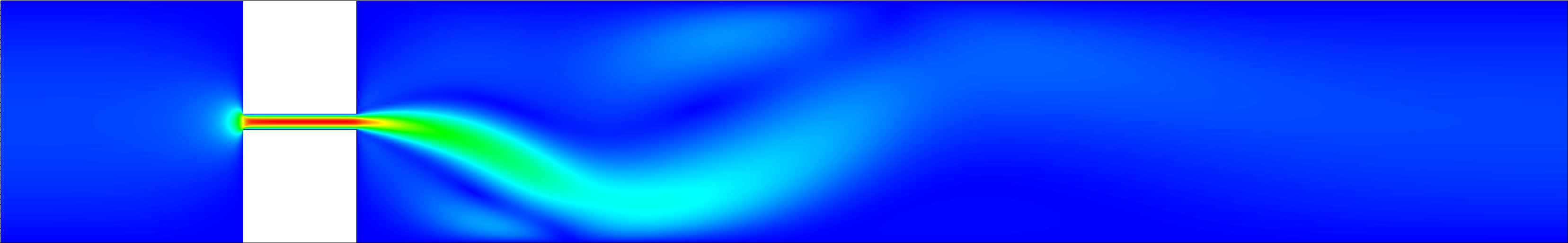}
    }
      \caption{2D case for $\lambda = 15.4$: solutions obtained with the full order (left) and reduced order (right) method for 
      (a) and (b) $\re_{2D} = 20$, (c) and (d) $\re_{2D} = 55$.
      }
  \label{fig:comparison_offline_online}
\end{figure}

We conclude this section by showing that there is no visible qualitative difference between
the solutions obtained with the full order method and with the RB method
for values of the Reynolds number not associated with the snapshots.
See the comparison in Fig.~\ref{fig:comparison_offline_online} 

\subsection{Unstable solution branch}
\label{ramo_instabile}
As already mentioned, for a given expansion ratio $\lambda$ and given $\re_{2D}$ the symmetric flow configuration exists
regardless of whether $\re_{2D}$ is smaller or grater than the critical value $\re_{2D,\mathrm{sb}}$ for the bifurcation.
Indeed, the symmetric branch is the only solution branch existing for Reynolds numbers 
below $\re_{2D,\mathrm{sb}}$, but for Reynolds numbers above $\re_{2D,\mathrm{sb}}$ it becomes unstable.
See, for example, the unstable symmetric flow configuration for $Re_{2D} = 27.7$ in Fig.~\ref{fig:snapshots_2d}(c)
and the corresponding stable asymmetric configuration in Fig.~\ref{fig:snapshots_2d}(d).



The numerical tests have shown that the RB approximation of the unstable branch can be achieved, 
but some care is required with the choice of the trial and test RB spaces.
One way to reconstruct the unstable branch is to use 
only the basis functions coming from the sampling of the unstable branch itself 
for both trial and test spaces. If this strategy is adopted, all the flow configurations 
of the unstable branch will be correctly approximated, but the bifurcation point will not be detected.
On the other hand, if basis functions coming from both the symmetric and the asymmetric 
branch are employed, the bifurcation point can be successfully detected but the fixed point 
scheme fails to converge after the bifurcation point, oscillating without damping between the two solution branches.
In this case, convergence to the symmetric or asymmetric branch after the bifurcation 
point can be achieved through e.g. a predictor-corrector or a pseudo-arclength continuation 
method (see~\cite{Dijkstra}) during the online phase,
with a further programming effort.

On the other hand, if one is interested only in the approximation of the stable solution branches,
there is no need for basis functions coming from the unstable branch
and no need for a continuation method in the online phase.
The reduced basis for the velocity is constructed only with 
basis functions arising from stable branches. This will allow to detect the bifurcation point and 
compute the stable solution for every parameter value. 


\subsection{2D case: two-parameter study}\label{2d_results_2param}

In this section, we still consider a slightly modified 2D channel: the part of the channel upstream of the sudden expansion 
in Fig.~\ref{fig:scheme_2d} is removed, 
since we focus now on the flow downstream of the contraction. The new geometry is thus a rectangle.
We let vary both the Reynolds number and the contraction width, so 
the parameter vector has now two components: $\bs{\mu}=(\re_{2D},\lambda)$.

In this very simple case, the geometry can be parametrized in two different ways:
\begin{itemize}
  \item[-] {\bf Geometric parametrization}: the contraction width is treated as explained in section~\ref{method} and the incompressibility constraint can be imposed through the Piola transformation as explained in section~\ref{sec:online};
  \item [-] {\bf Boundary condition parametrization}: the different aspect ratio of the contraction is imposed by parametrizing the boundary conditions. Indeed, a channel with a contraction of width $\lambda$ will produce in our model a parabolic inner velocity profile dependent on $\lambda$: 
    \begin{equation}
      v_x=
      \begin{cases}
        \displaystyle{-\frac{(y-\lambda L_c)(y+\lambda L_c)}{\lambda^2L_c^2}}\qquad &\text{ if }-\lambda L_c\le y\le \lambda L_c \\
        0 & \text{ otherwise,}
      \end{cases}
      \label{eq:parametrized_inlet}
    \end{equation}
    where the $y$ coordinate has origin on the symmetry axis of the contraction.
\end{itemize}

One advantage of the second strategy is that the RB functions are automatically divergence-free and 
the relatively complex procedure of the Piola transformation does not need to be performed. 
Thus, we choose the boundary condition parametrization.
However, we need to be careful in imposing the inlet velocity profile because 
the mass flow rate constraint as expressed in equations~\eqref{eq:def_average_phi_x} 
and~\eqref{eq:def_lagrange_matrix} is not sufficient to ensure uniqueness of the RB solution.
One possible workaround for this issue is to split the boundary integral~\eqref{eq:def_mass_flow_rate} 
used for the mass flow rate constraint in two parts:
\begin{equation}
\int_{\Ga_\mathrm{in}}v_x\de x=\int_{\Ga_0}v_x\de x+\int_{\Ga_\lambda}v_x\de x,
\label{eq:mass_flow_rate_split}
\end{equation}
where $\Ga_\mathrm{in}$ is the part of $\partial\Om$ where the inlet velocity profile is imposed, 
$\Ga_0$ the part of $\Ga_\mathrm{in}$ where $v_x = 0$ and $\Ga_\lambda$ the part of $\Ga_\mathrm{in}$ where $v_x\ne0$. 
Notice that $\Ga_0\cup\Ga_\lambda=\Ga_\mathrm{in}$ and $\Ga_0\cap\Ga_\lambda=\emptyset$.
We introduce two Lagrange multipliers $\alpha_0$ for $\Ga_0$ and $\alpha_\lambda$ for $\Ga_\lambda$, in order to enforce~\eqref{eq:parametrized_inlet} in integral form as:
    \begin{equation}
      \alpha_0\int_{\Ga_0}v_x\de y=0 \qquad \alpha_\lambda\int_{\Ga_\lambda}v_x\de y=\dot{w_x}
      \label{eq:two_multipliers}
    \end{equation}
    Finally, the two Lagrange multipliers $\alpha_0$ and $\alpha_\lambda$ are treated as additional unknowns, and a linear system analogous to that in equation~\eqref{eq:def_lagrange_matrix} is solved.

The GLC collocation sampling has been carried out on the kinematic viscosity set 
$\nu\in[1.5,5]\cdot10^{-3}$ and on the contraction width set $w_c\in[1/10,1/2]$. 
We obtained 6 values for the kinematic viscosity $\nu=1.5\cdot10^{-3},1.73446\cdot10^{-3},2.375\cdot10^{-3},4.125\cdot10^{-3},4.76554\cdot10^{-3},5\cdot10^{-3}$
and 7 values for the expansion ratio $\lambda =2,3,4,5,6,8,10$, so $N=42$.
Note that the sampling has not been performed directly on the Reynolds number 
due to its dependence on the contraction width.
In table~\ref{tab:cmp_drikakis}, we report the critical Reynolds numbers for the symmetry breaking computed with the RB method for different values of the contraction width.  

\begin{table}[h!]
  \caption{Symmetry breaking Reynolds numbers as a function of the channel contraction width for the 2D case with variable geometry.}
\label{tab:cmp_drikakis}
\centering
\begin{tabular}[tbp]{lccccccc}
\toprule
$\lambda$ & 2 & 3 & 4 & 5 & 6 & 8 & 10 \\
$\re_{2D,sb}$     &  222.5  &  73.5   & 59.125  & 42.75   & 34.5    &  28.5   &  27.0        \\
\bottomrule
\end{tabular}
\end{table}

\begin{figure}[h!]
  \centering
  \includegraphics[width=0.5\textwidth]{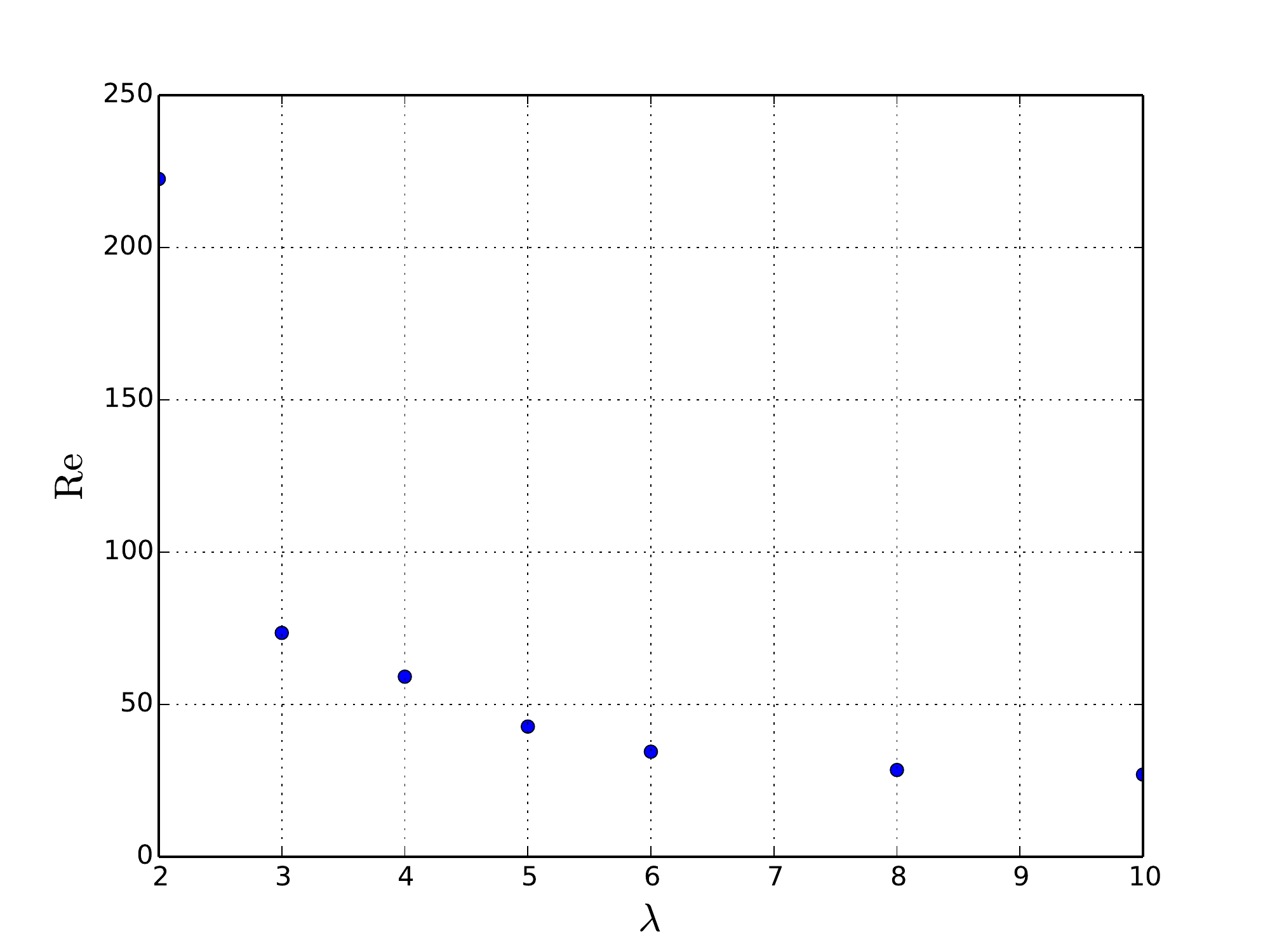}
  \caption{2D case: value of the Reynolds number at the bifurcation point as a function of the expansion ratio $\lambda$. 
}
  \label{fig:lambda-Re}
\end{figure}

The same values in table~\ref{tab:cmp_drikakis} are plotted in figure~\ref{fig:lambda-Re}. 
We clearly see that as the aspect ratio $\lambda$ decreases, the critical Reynolds number for the symmetry breaking increases, 
as observed also in \cite{Drikakis}. We see that $\re_{2D,sb}$ decreases fast for small values of $\lambda$,
while it decreases mildly for $\lambda \geq 6$ (also recall that for $\lambda = 15.4$ we found $\re_{2D,sb} \approx 26$).
We remark that also in this case the results match closely the ones in~\cite{Drikakis}.

Let us analyze the computational time savings allowed by our RB method.
Since we have $N=42$ with a $10h$ cputime needed per single run,
the computational time analysis is given by: 
\[
\frac{\text{Time to build the RB spaces}+\text{Online time to detect the bifurcation point}}{\text{Time of the equivalent full order computation}}=\frac{42\cdot10h+2h+7\cdot0.08h}{7\cdot10h\cdot80}\simeq 7.5\%.
\]
The ratio between a single online reduced order run and a single full order one are  the same as the one considered in the single parameter case (order $10^{-4}$).
The \emph{break-even}, comparing all the offline computational times needed to prepare the reduced basis problem ($N=42$) and  an online query with full order model is:
\[
  \frac{\text{All full order computations for RB prep.}}{\text{Full order one query  comp. time}}=\frac{42\cdot10h+2 h}{10h}\simeq 42.2.
\]
After 43 queries a reduced order computational model brings savings.

\begin{figure}[htbp]
  \centering
    \includegraphics[width=.6\textwidth]{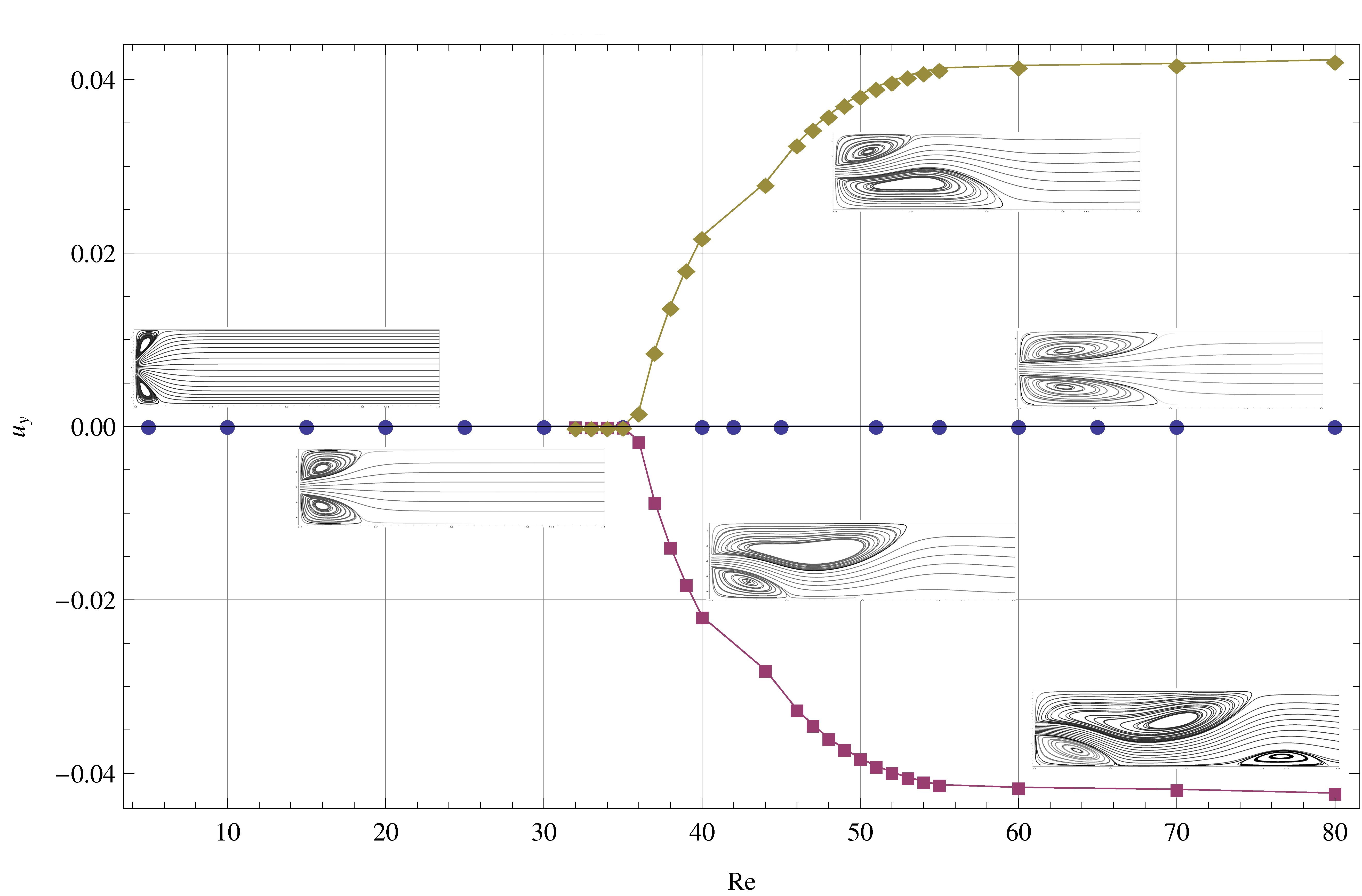}
    \caption{Bifurcation diagram obtained with the Reduced Order Model for $\lambda=6$: vertical component of the velocity $u_y$
    taken on the horizontal axis, at distance 1 from the inlet, versus the Reynolds number.}
  \label{fig:bifurcation_1_6}
\end{figure}

For $\lambda=6$, which is one value among those listed in table~\ref{tab:cmp_drikakis}, 
we plot in figure~\ref{fig:bifurcation_1_6} the vertical component of the velocity is taken on
the horizontal axis, at distance 1 from the inlet, versus the Reynolds number. 
This bifurcation diagram with both the stable and unstable solution branches
 compares very well with the one in~\cite{Drikakis}, 
but it has been obtained at a fraction of the computational time as explained above.

\begin{figure}[h!]
  \centering
     \subfigure[Unstable solution]{
    \includegraphics[width=0.33\columnwidth]{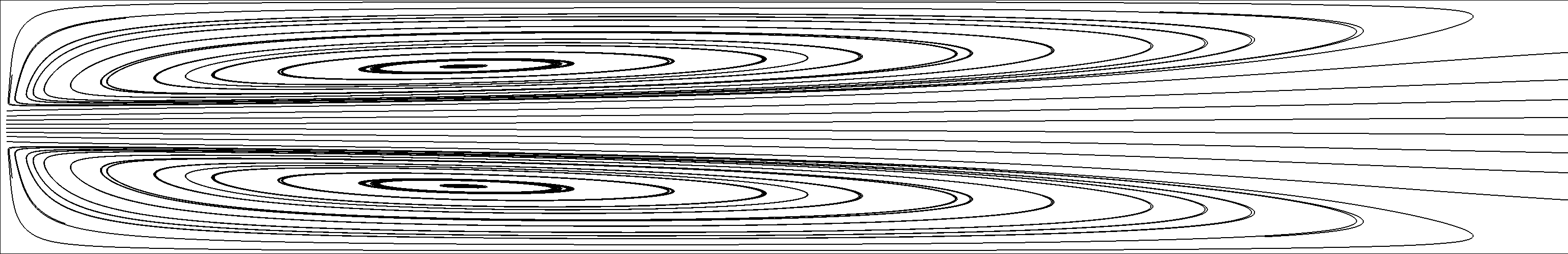}
    }
     \subfigure[Stable solution]{
    \includegraphics[width=0.33\columnwidth]{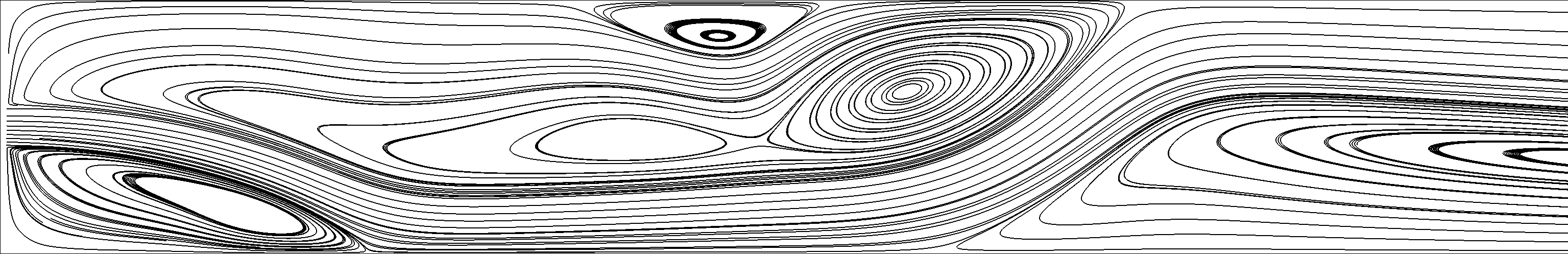}
    }
  \caption{2D case for $\lambda = 6$: streamlines for the (a) stable asymmetric and (b) unstable symmetric 
  solutions at $\re_{2D}=600$.}
  \label{fig:snapshots_2d_600}
\end{figure}

Keeping $\lambda = 6$, we check how the flow evolves as $Re_{2D}$ is pushed to a higher value,
well beyond the parameter range considered in this work. 
Fig.~\ref{fig:snapshots_2d_600} reports the streamlines of both the unstable and stable 
solution at $\re_{2D}=600$. The stable solution in Fig.~\ref{fig:snapshots_2d_600}(b) shows that
the flow structure becomes more complex, with existing recirculations changing shape and growing in size. 
This is consistent with the results presented in \cite{Drikakis,AQpreprint}.

\subsection{3D case}
\label{3d_results}

The three-dimensional channel in Fig.~\ref{fig:scheme_3d} has been obtained by extruding the
two-dimensional geometry in Fig.~\ref{fig:scheme_2d} along the $z$-axis.
Thus, for the 3D case, we would have three parameters: the Reynold number $\re_{3D}$,
the contraction width, and the channel depth.
However, since we have already investigated in Sec.~\ref{2d_results_2param} the influence 
of the expansion ratio (i.e., the contraction width) on the critical Reynolds number
for the symmetry breaking, we fix the contraction width and consider the Reynolds number and the channel depth
as the only parameters.

We set the expansion ratio $\lambda$ to $15.4$, due to the richness of flow patterns described in
Sec.~\ref{2d_results_1param} and reference~\cite{Oliveira}. 
Of course, we expect the vortex structure to be much more complex than in the 2D case.
We are interested in understanding how varying the Reynolds number and the aspect ratio $\AR$ (and thus $\HH$) affects the flow in the expansion channel. 
The goal of this section is to evaluate the effect of the walls on the bifurcating phenomenon. 
Intuitively, when the walls are very far apart (large values of $\HH$), their influence on the central region 
of the channel will be quite small, and the flow pattern can be expected to be close to the 2D case.
On the other hand, when the walls are very close with respect to the channel height (small values of $\HH$), 
a relatively large fraction of the sectional area will be occupied by low velocity fluid. As a result, we can expect 
that the bifurcation will take place at higher Reynolds numbers.
As reported in tables \ref{tab:which_reynolds} and \ref{tab:sampling_data}, we sample sample nine values 
for the Reynolds number in the interval $\re_{2D}\in[0.01,90]$ and eight value of $\HH$. Notice
that the eighth ``value'' of $\HH$ in table \ref{tab:sampling_data} corresponds to the 2D case. 

\begin{figure}[h!]
  \centering
  \subfigure[$\re_{3D} = 0.01$]{
  \includegraphics[width=0.45\columnwidth]{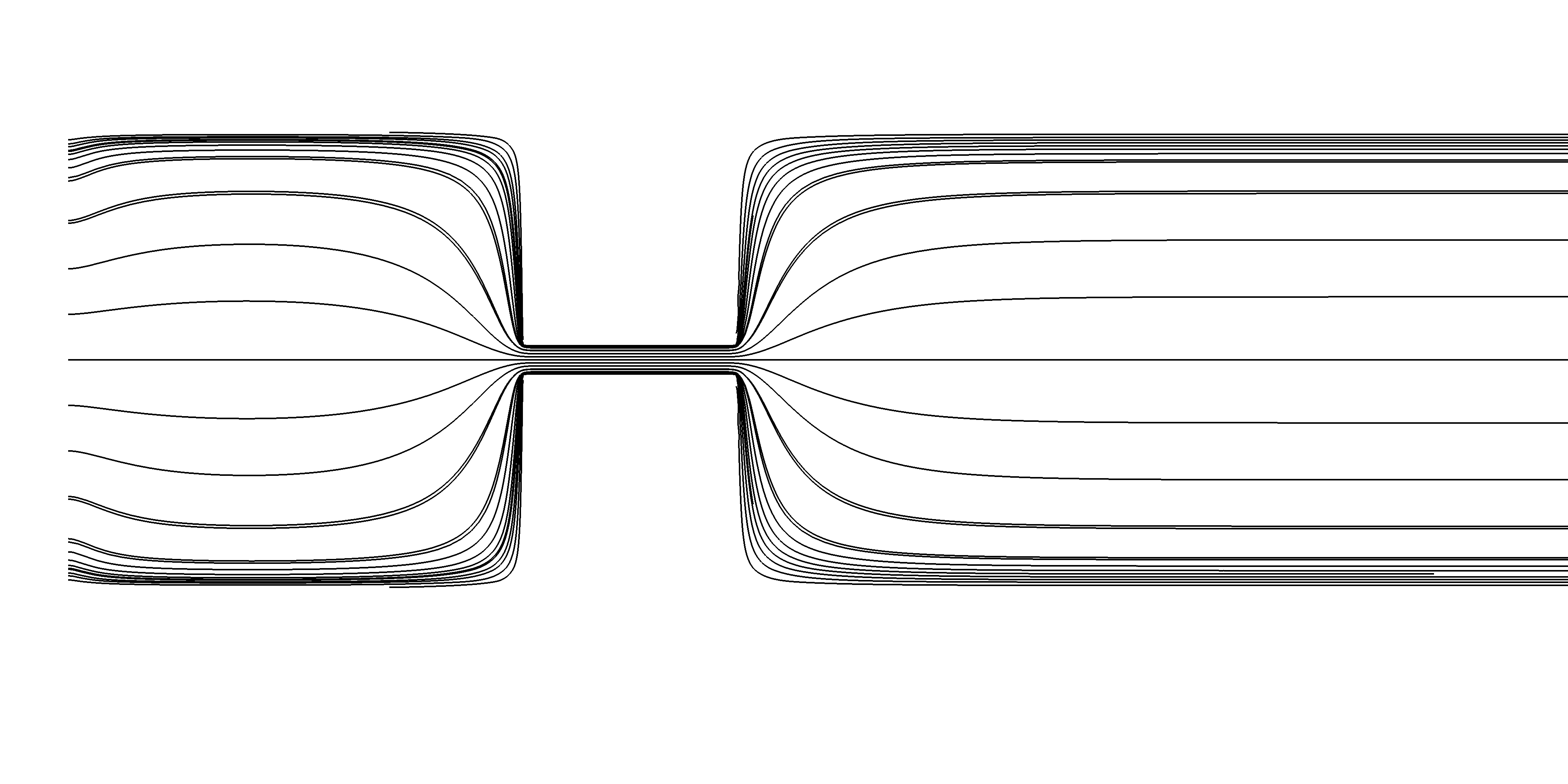}
  }
  \subfigure[$\re_{3D} = 27.786$]{
  \includegraphics[width=0.45\columnwidth]{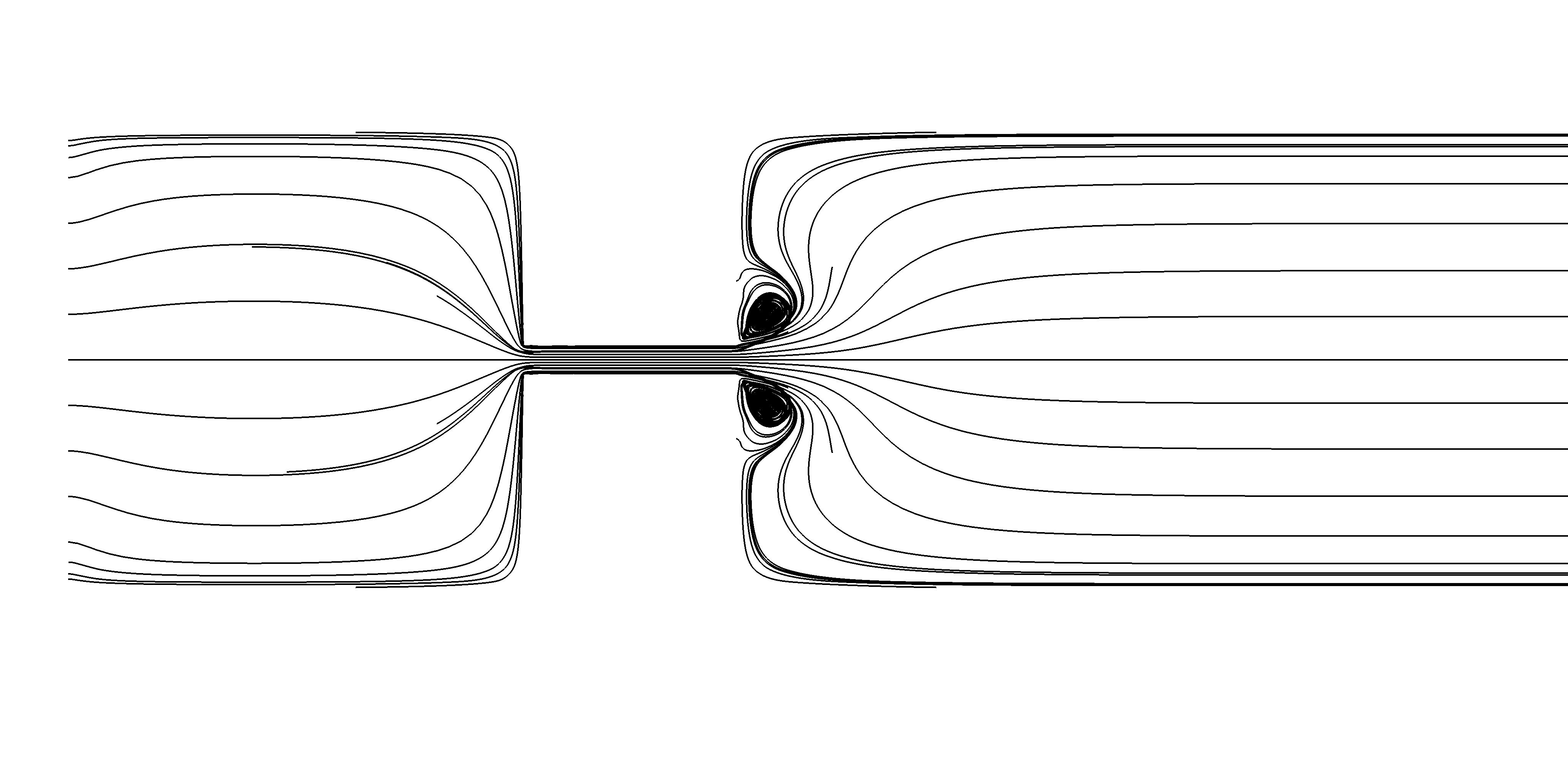}
  }
  \subfigure[$\re_{3D} = 45.005$]{
  \includegraphics[width=0.45\columnwidth]{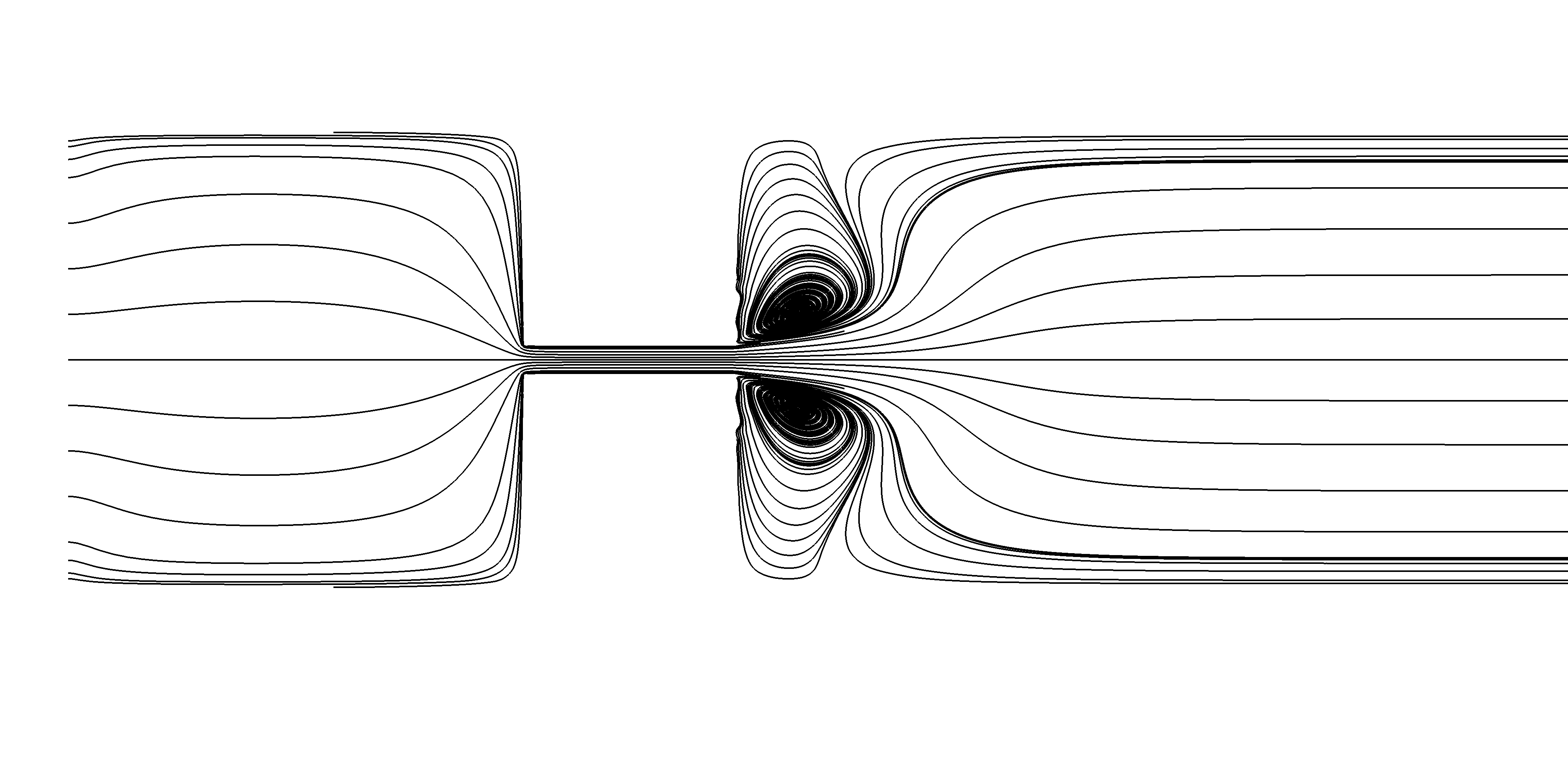}
  }
  \subfigure[$\re_{3D} = 62.224$]{
  \includegraphics[width=0.45\columnwidth]{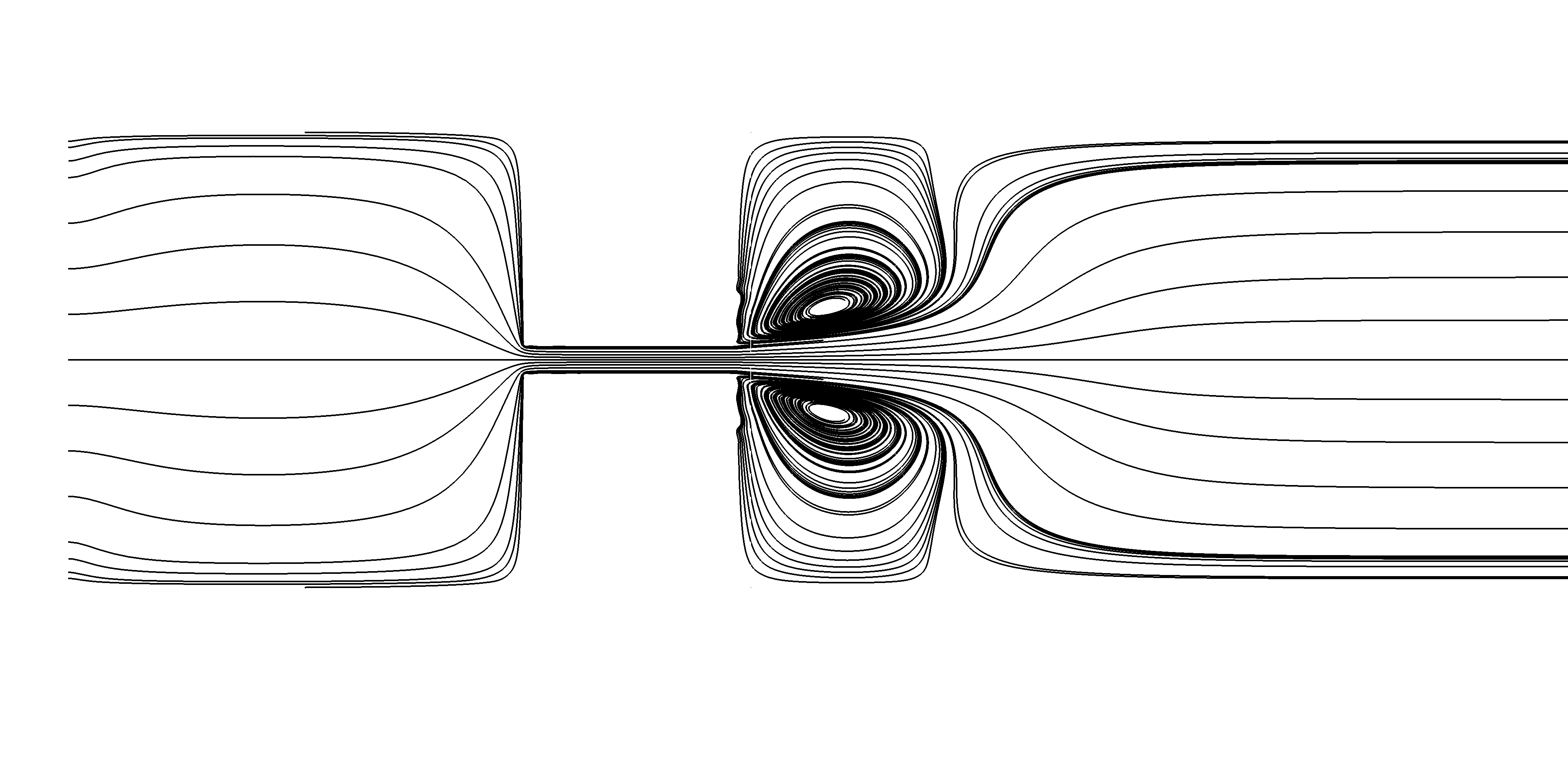}
  }
  \subfigure[$\re_{3D} = 76.821$]{
  \includegraphics[width=0.45\columnwidth]{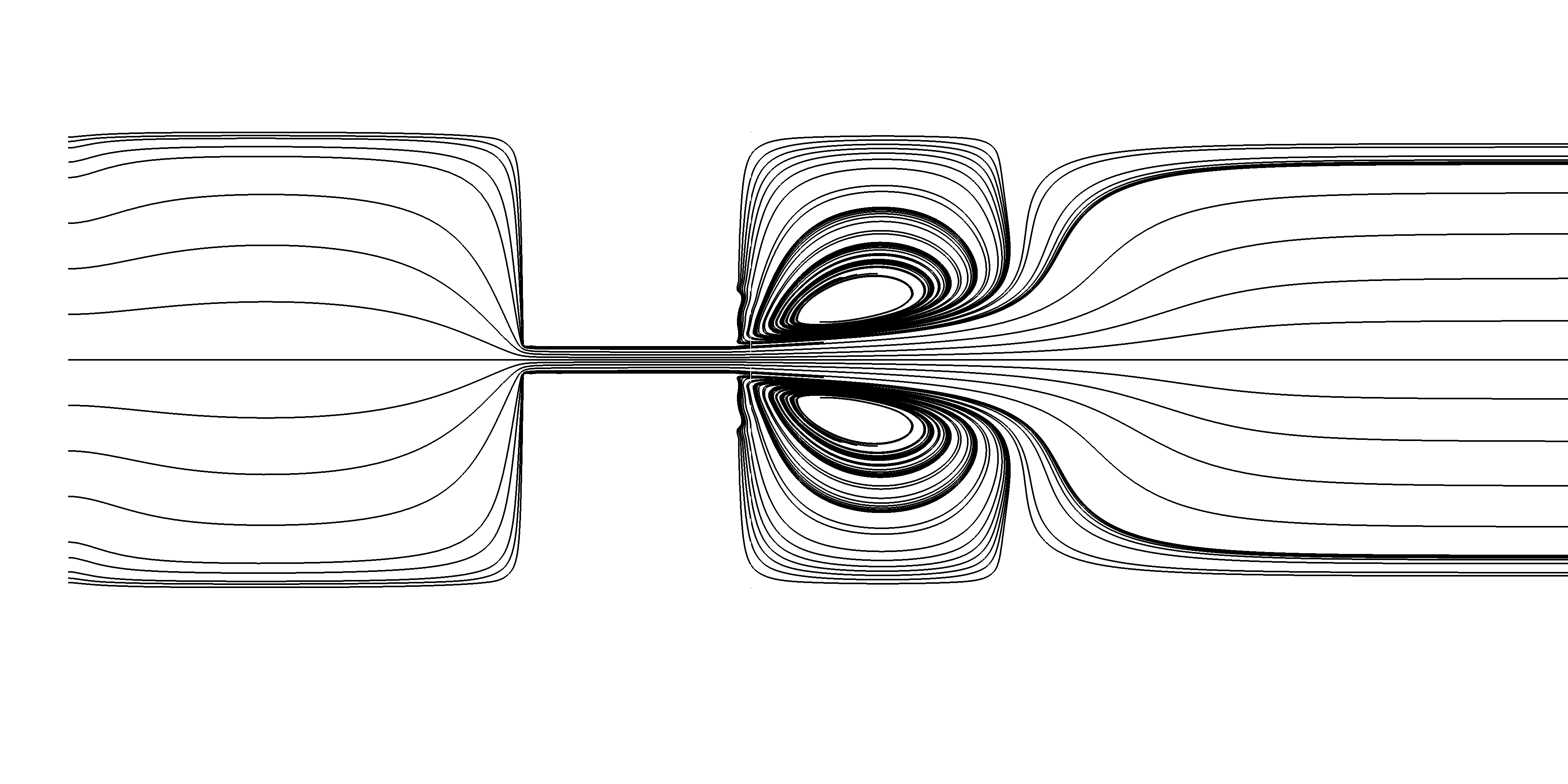}
  }
  \caption{3D case for $\lambda = 15.4$ and $\HH=0.6210$: streamlines on the $xy$-plane (see figure~\ref{fig:scheme_3d}) for 
  (a) $\re_{3D}=0.01$, (b) $\re_{3D}=27.786$, (c) $\re_{3D}=45.005$, (d) $\re_{3D}=62.224$, and (e) $\re_{3D}=76.821$.} 
  \label{fig:3d_n4}
\end{figure}

In order to show the sequence of events as the Reynolds number is increased when 
the aspect ratio is fixed, we set it to 1.6398 which corresponds to $\HH = 0.6210$.
In Fig.~\ref{fig:3d_n4}, we display the streamlines on the $xy$-plane
for different values of the Reynolds number $\re_{3D}$.
At $\re_{3D} = 0.01$, the 3D flow looks similar to the 2D flow: 
(compare Fig.~\ref{fig:3d_n4}(a) with Fig.~\ref{fig:snapshots_2d}(a)) but it features smaller Moffatt eddies. 
As the Reynolds number increases, ``lip vortices'' form, as shown in Fig.~\ref{fig:3d_n4}(b). 
This is in agreement with the observations in \cite{Oliveira} and references therein.
The size of the lip vortices increases as $\re_{3D}$ increases and once they reach the
corner, the vortices continue to grow in the downstream direction, i.e. along the $x$-axis. 
See Fig.~\ref{fig:3d_n4}(c), (d), and (e). By convention, once they expand in the downstream direction
they are called ``corner vortices''. 
Notice that the flow downstream of the expansion is symmetric about the $xz$-plane
up to $\re_{3D} = 76.821$, while asymmetries in 2D (i.e., for $\HH = 1$) arise around
$\re_{2D}=26$. 


\begin{figure}[h!]
  \centering
    \includegraphics[width=0.45\columnwidth]{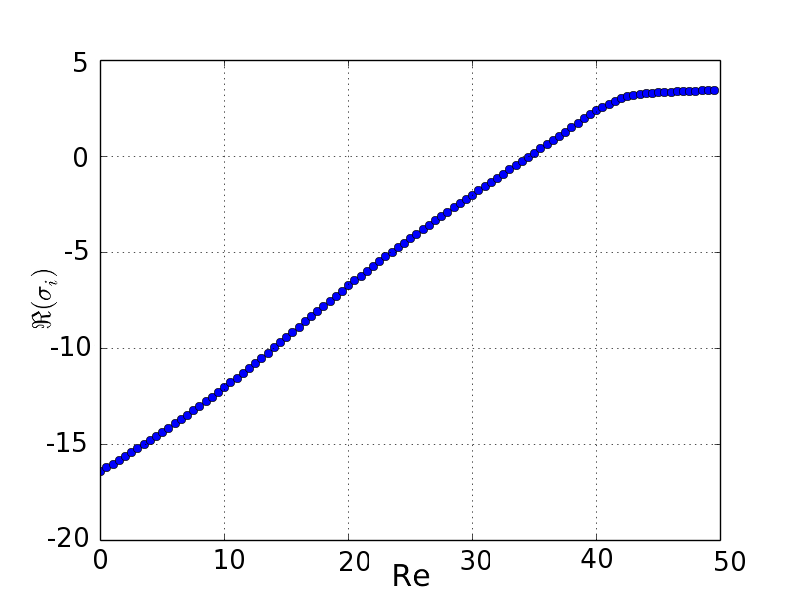}
    \caption{3D case for $\lambda = 15.4$ and $\HH = 0.9517$: 
    real part of the eigenvalue of matrix $\mathtt{L}$ in (\ref{eq:def_L_matrix}) responsible 
    for the symmetry breaking as a function of the Reynolds number in a neighborhood of a bifurcation point. 
    }
  \label{fig:eig_3d}
\end{figure}

\begin{figure}[h!]
  \centering
  \subfigure[eigenvalues of matrix $\mathtt{L}$]{
    \includegraphics[width=0.45\columnwidth]{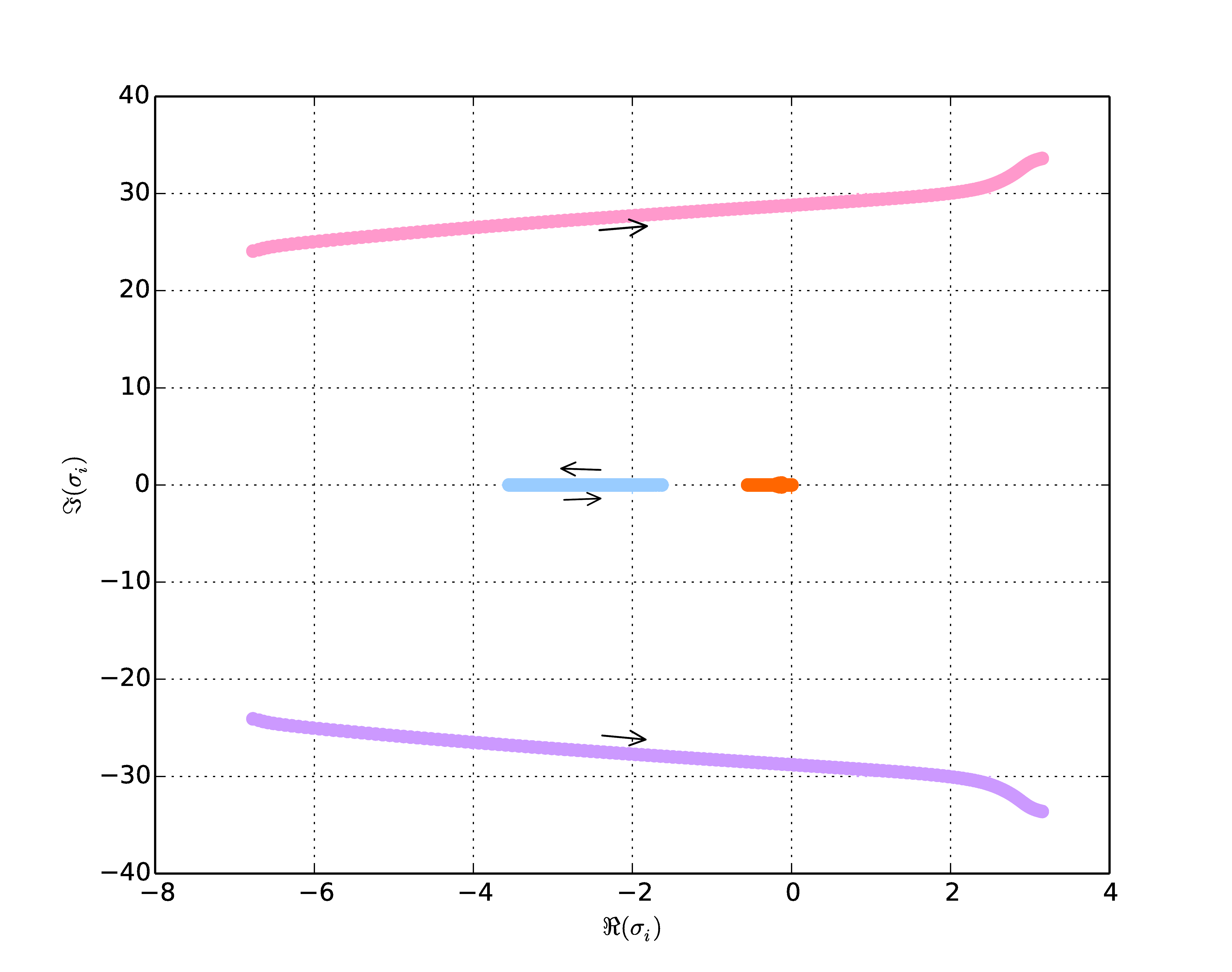}
  }
  \subfigure[zoomed-in view of (a)]{
    \includegraphics[width=0.45\columnwidth]{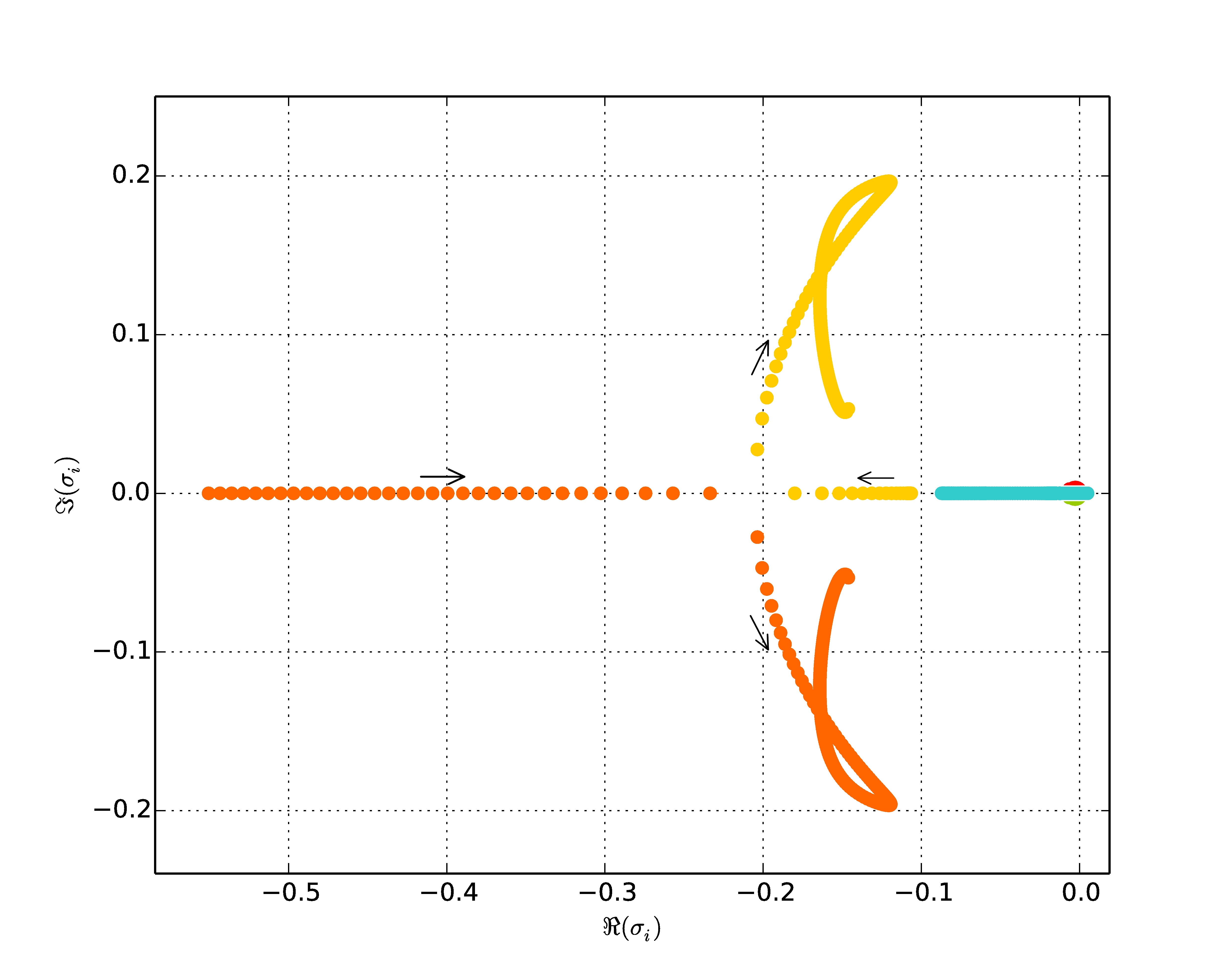}
}
  \subfigure[zoomed-in view of (b)]{
    \includegraphics[width=0.45\columnwidth]{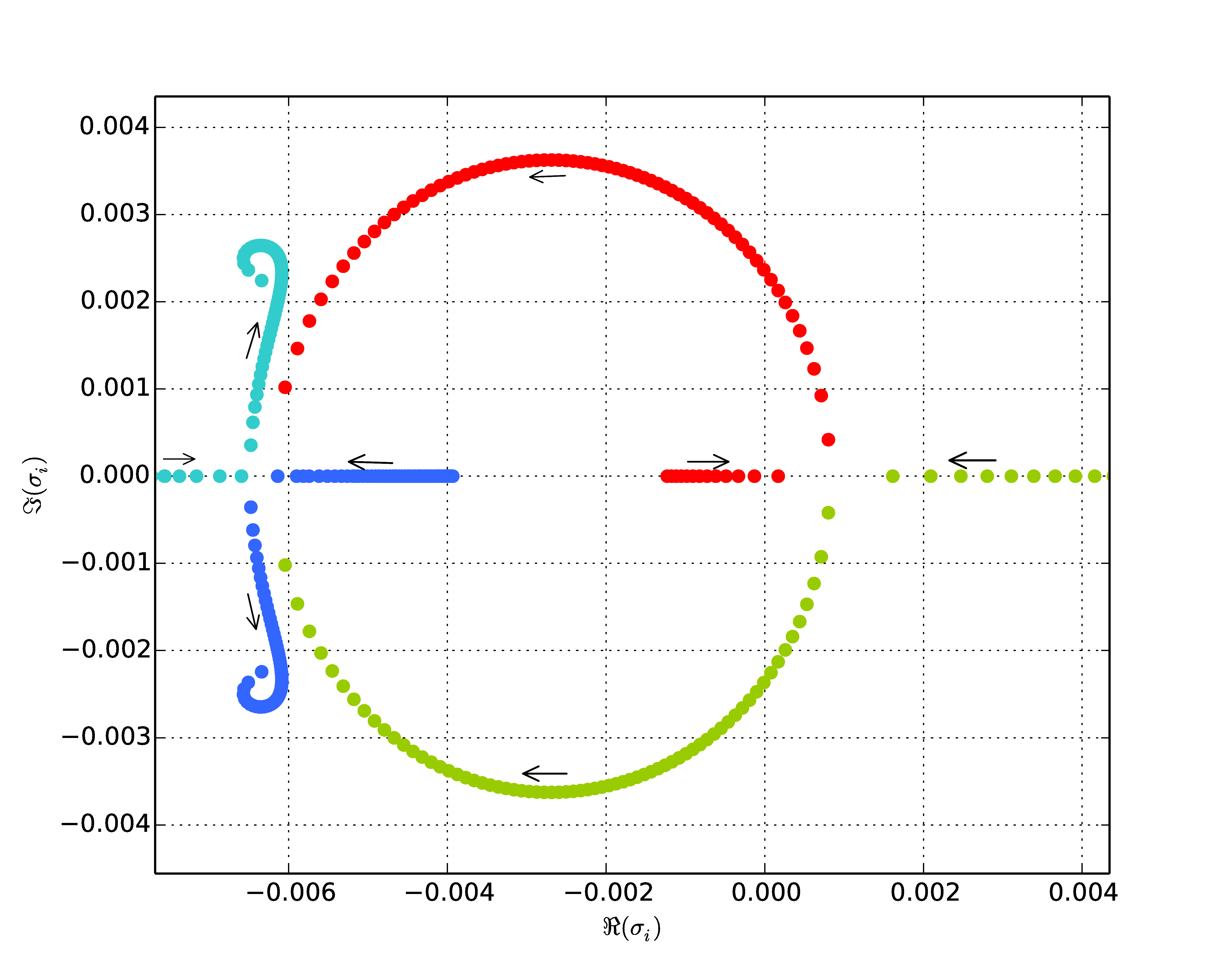}
  }
    \caption{3D case for $\lambda = 15.4$ and $\HH = 0.9517$: (a) path of the eigenvalues of matrix $\mathtt{L}$ in (\ref{eq:def_L_matrix}) 
    in the complex plane for $\re_{3D}\in[0.01,90]$. 
    Subfigure (b) is a zoomed-in views of subfigure (a) and 
    subfigure (c) is a  zoomed-in views of subfigure (b). 
    The arrows indicate the direction of the increasing Reynolds numbers.
    The eigenvalue in red on the real axis in (c) is responsible for the bifurcation point.} 
    \label{fig:eig_3d_color}
\end{figure}

Let us consider the geometry with $\HH=0.9517$, which corresponds to
the largest aspect ratio among those in Table \ref{tab:sampling_data} for which we have 
an actual 3D geometry. 
We proceed with the computation of the symmetry breaking bifurcation point 
using the bifurcation detection method described in section~\ref{sec:bif_detection}. 
Since $\HH$ is fixed, we consider a total of 9 basis functions for the online computation, 
corresponding to the different values of Reynolds number reported in table~\ref{tab:sampling_data}. 
In Fig.~\ref{fig:eig_3d}, we plot the real part of the eigenvalue of matrix $\mathtt{L}$ in (\ref{eq:def_L_matrix})
responsible for the symmetry breaking.
We see that the curve crosses the horizontal axis at a Reynolds number of about $35$. This 
coincides with the critical value for the symmetry breaking reported by \cite{Oliveira}. 
For the sake of completeness, in figure~\ref{fig:eig_3d_color} we report the path of all the eigenvalues in the complex plane.

For $\lambda = 15.4$, the critical Reynolds number for the symmetry breaking in the 2D geometry (i.e., $\HH = 1$) 
found in Sec.~\ref{2d_results_1param} is $\re_{2D,\mathrm{sb}}=26$. See Fig.~\ref{fig:eig_2d}.
When $\HH$ is decreased to $0.9517$, the critical Reynolds number for the symmetry breaking
increases to $\re_{3D,\mathrm{sb}}=35$, as shown in Fig.~\ref{fig:eig_3d}. 
If $\HH$ is further decreased to $0.6210$, we saw in Fig.~\ref{fig:3d_n4} that the flow remains 
symmetric up to $\re_{3D} = 76.821$. As expected, at low values of $\HH$ the proximity of vertical 
walls make the flow fully three-dimensional (instead of quasi-2D) inhibiting the symmetry breaking.
Thus, as $\HH$ decreases $\re_{3D,\mathrm{sb}}$ becomes larger and larger.

\begin{figure}[h!]
  \centering
  \subfigure[$\re_{3D} = 27.79$, $\HH = 0.8165$, $xy$-plane]{
  \includegraphics[width=0.45\columnwidth]{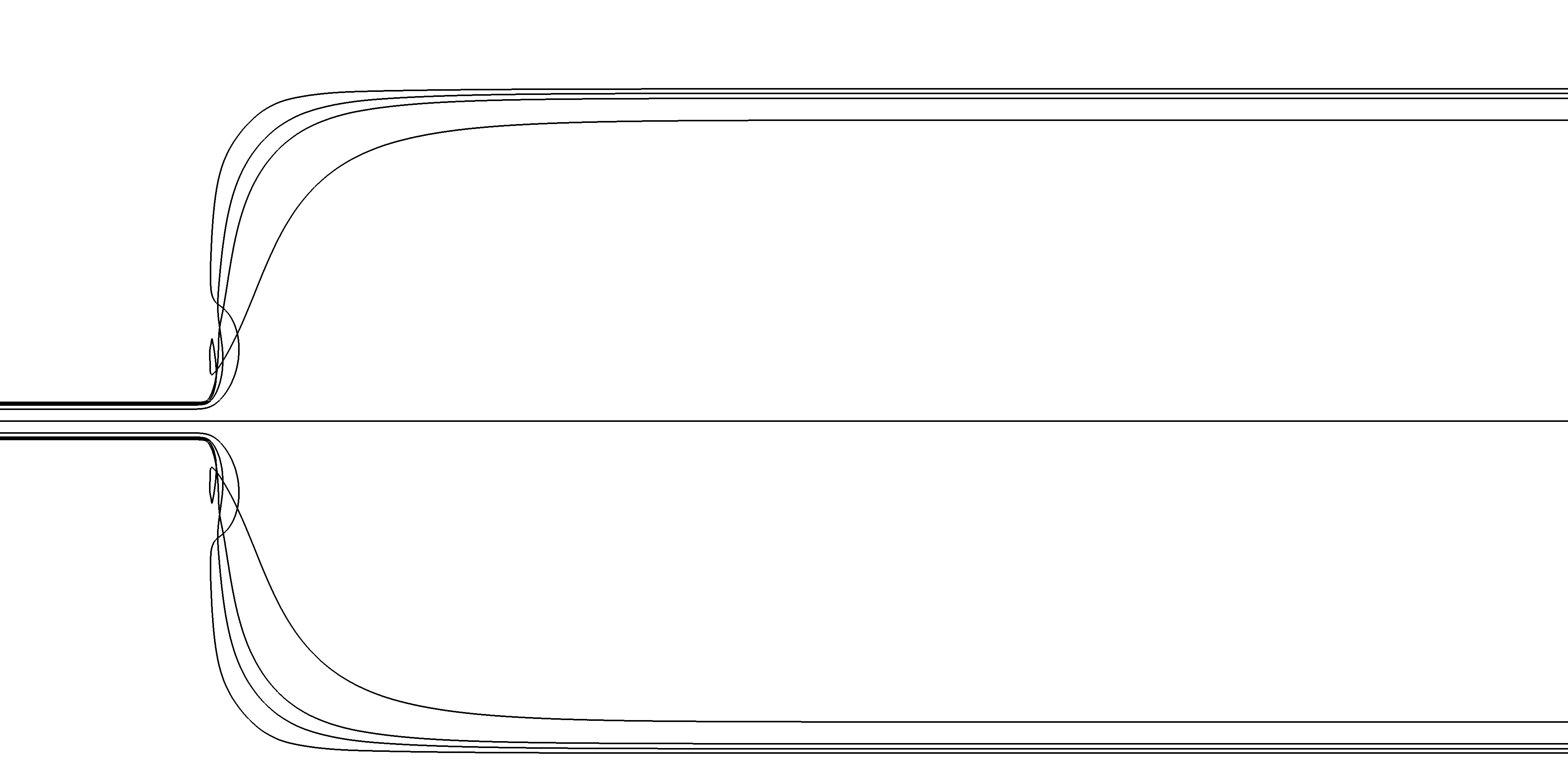}
  }
  \subfigure[$\re_{3D} = 27.79$, $\HH = 0.8165$, $yz$-plane]{
  \includegraphics[width=0.45\columnwidth]{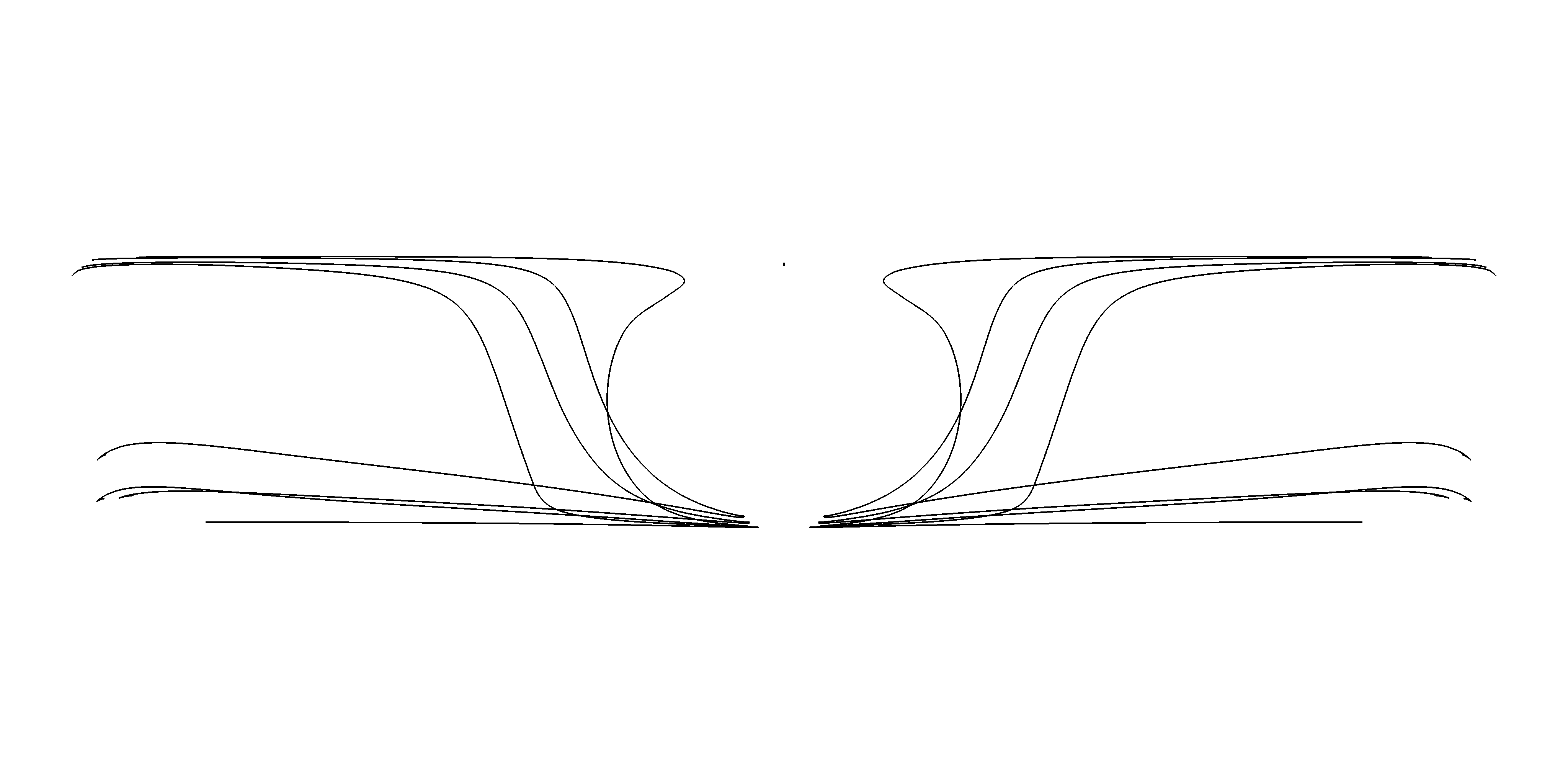}
  }
   \subfigure[$\re_{3D} = 27.79$, $\HH = 0.9517$, $xy$-plane]{
  \includegraphics[width=0.45\columnwidth]{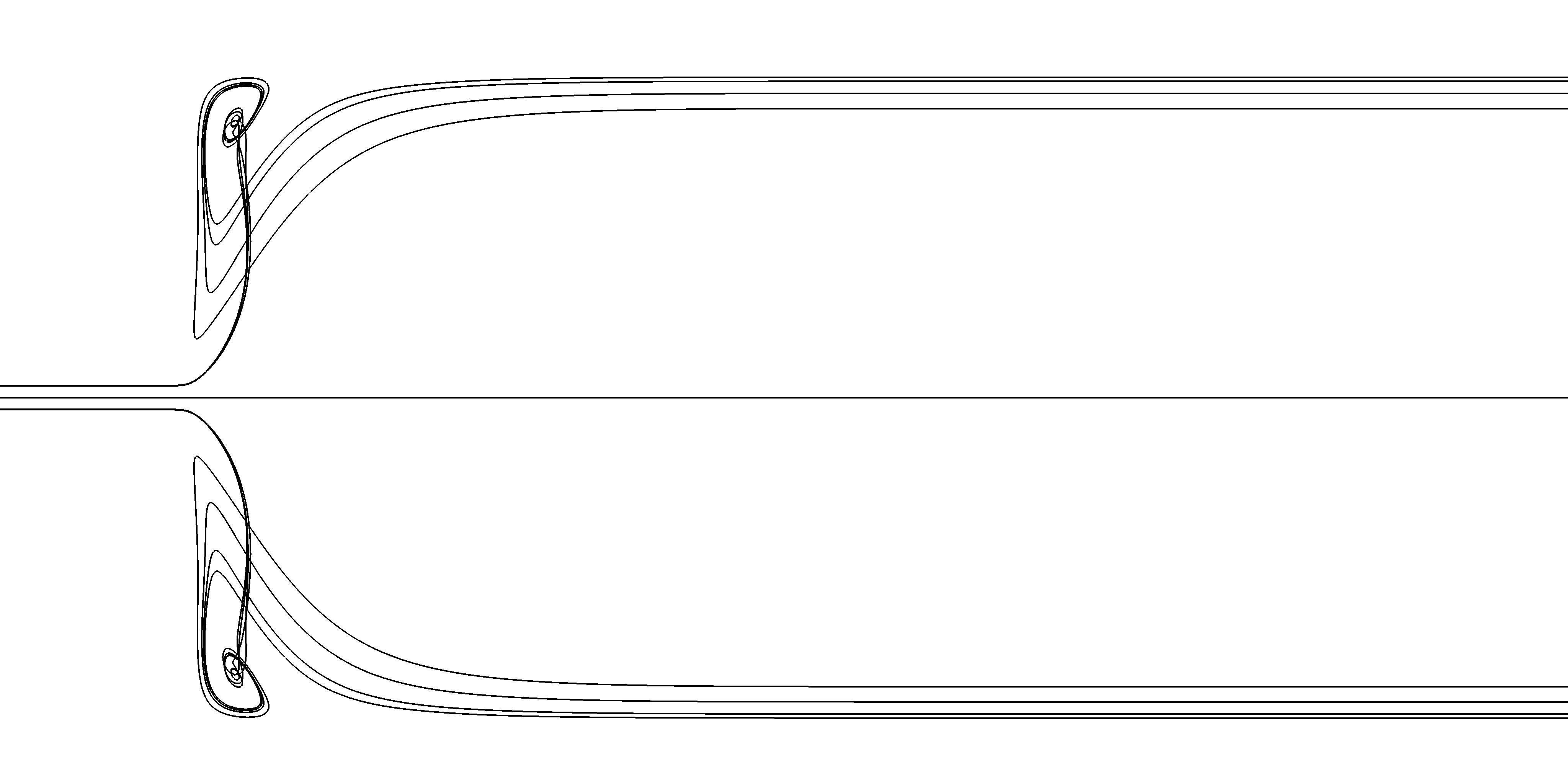}
  }
  \subfigure[$\re_{3D} = 27.79$, $\HH = 0.9517$, $yz$-plane]{
  \includegraphics[width=0.45\columnwidth]{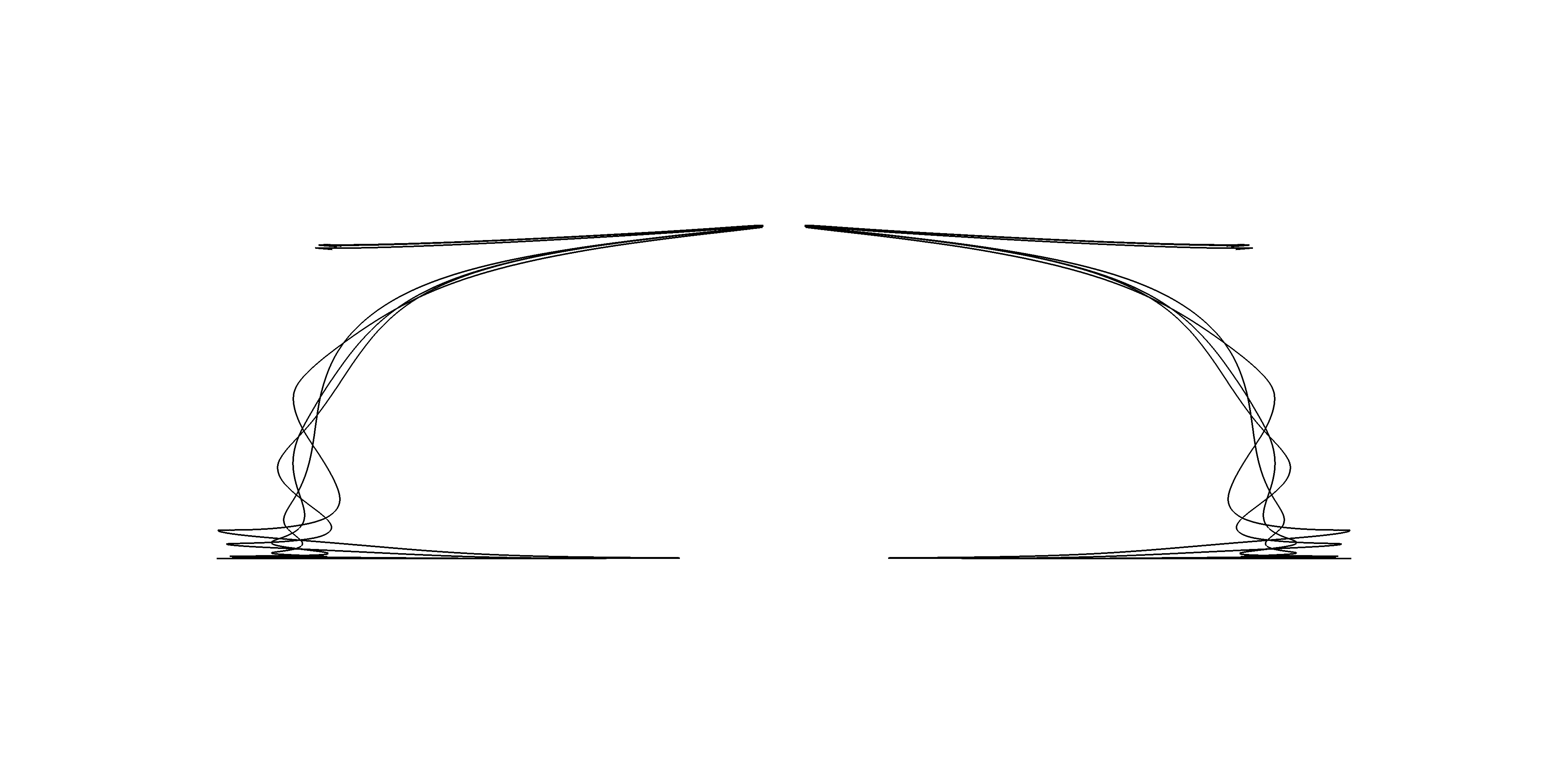}
  }
   \subfigure[$\re_{3D} = 45.01$, $\HH = 0.9517$, $xy$-plane]{
  \includegraphics[width=0.45\columnwidth]{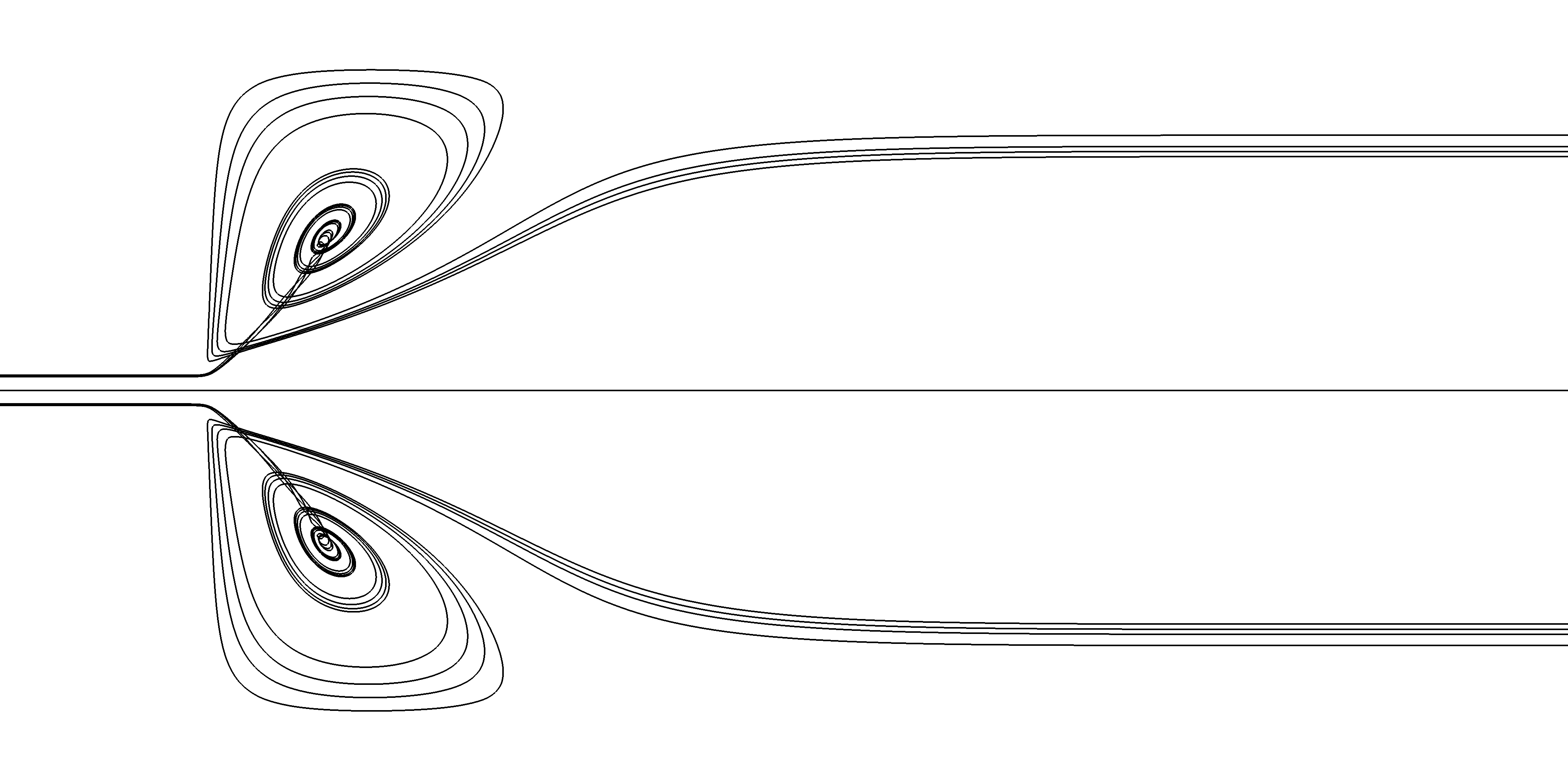}
  }
  \subfigure[$\re_{3D} = 45.01$, $\HH = 0.9517$, $yz$-plane]{
  \includegraphics[width=0.45\columnwidth]{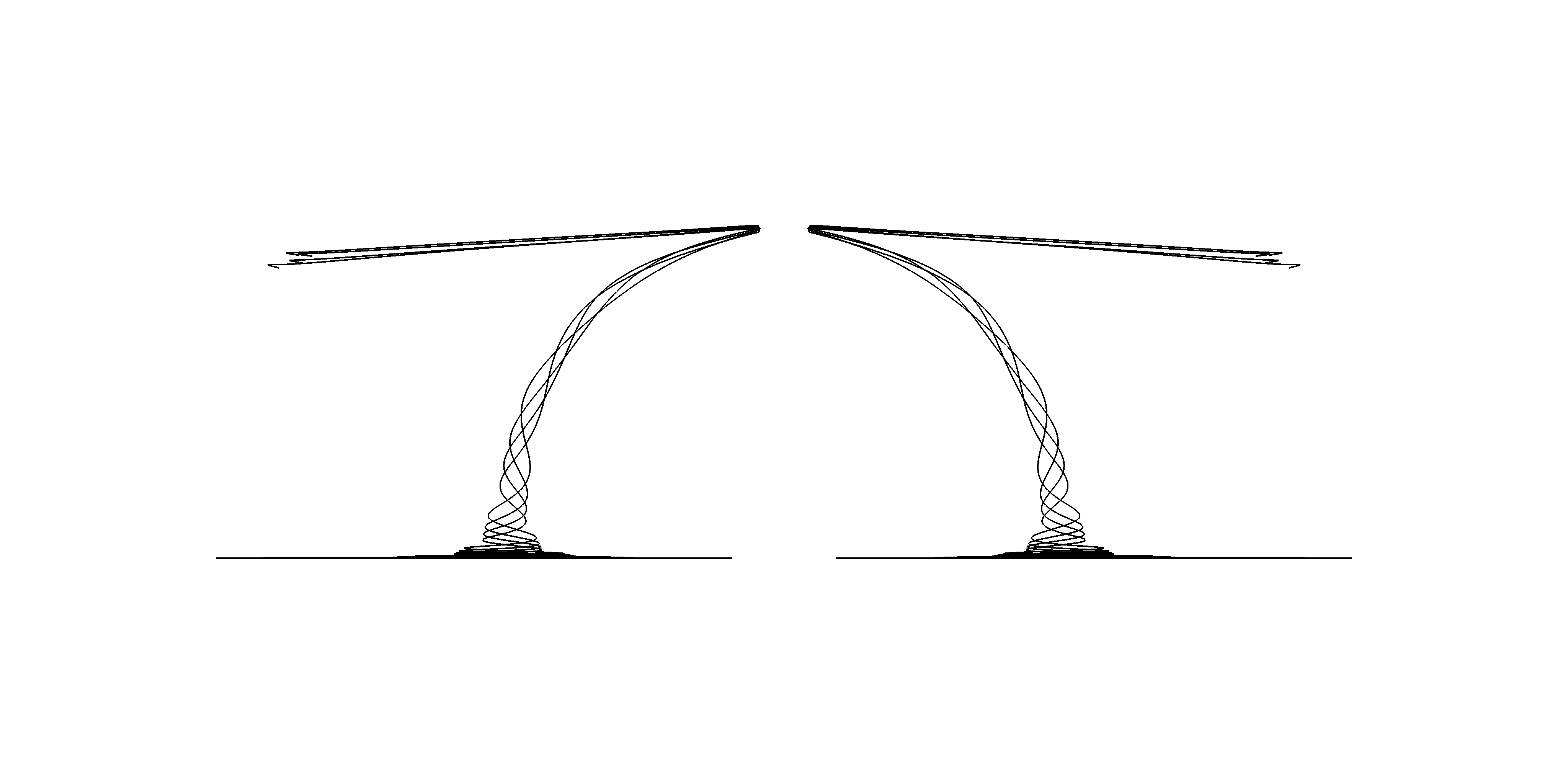}
  }
   \subfigure[$\re_{3D} = 62.22$, $\HH = 0.8165$, $xy$-plane]{
  \includegraphics[width=0.45\columnwidth]{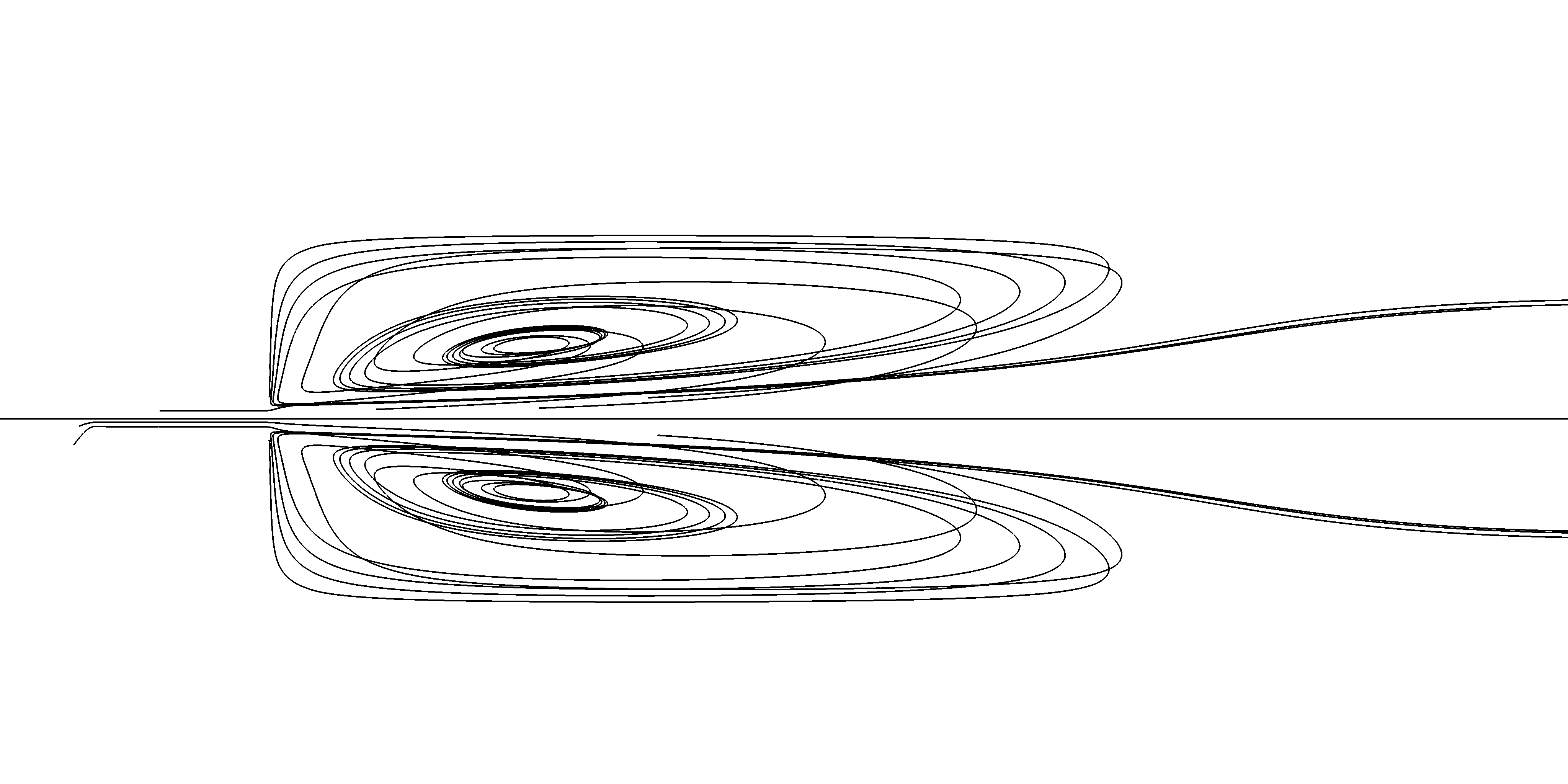}
  }
  \subfigure[$\re_{3D} = 62.22$, $\HH = 0.8165$, $yz$-plane]{
  \includegraphics[width=0.45\columnwidth]{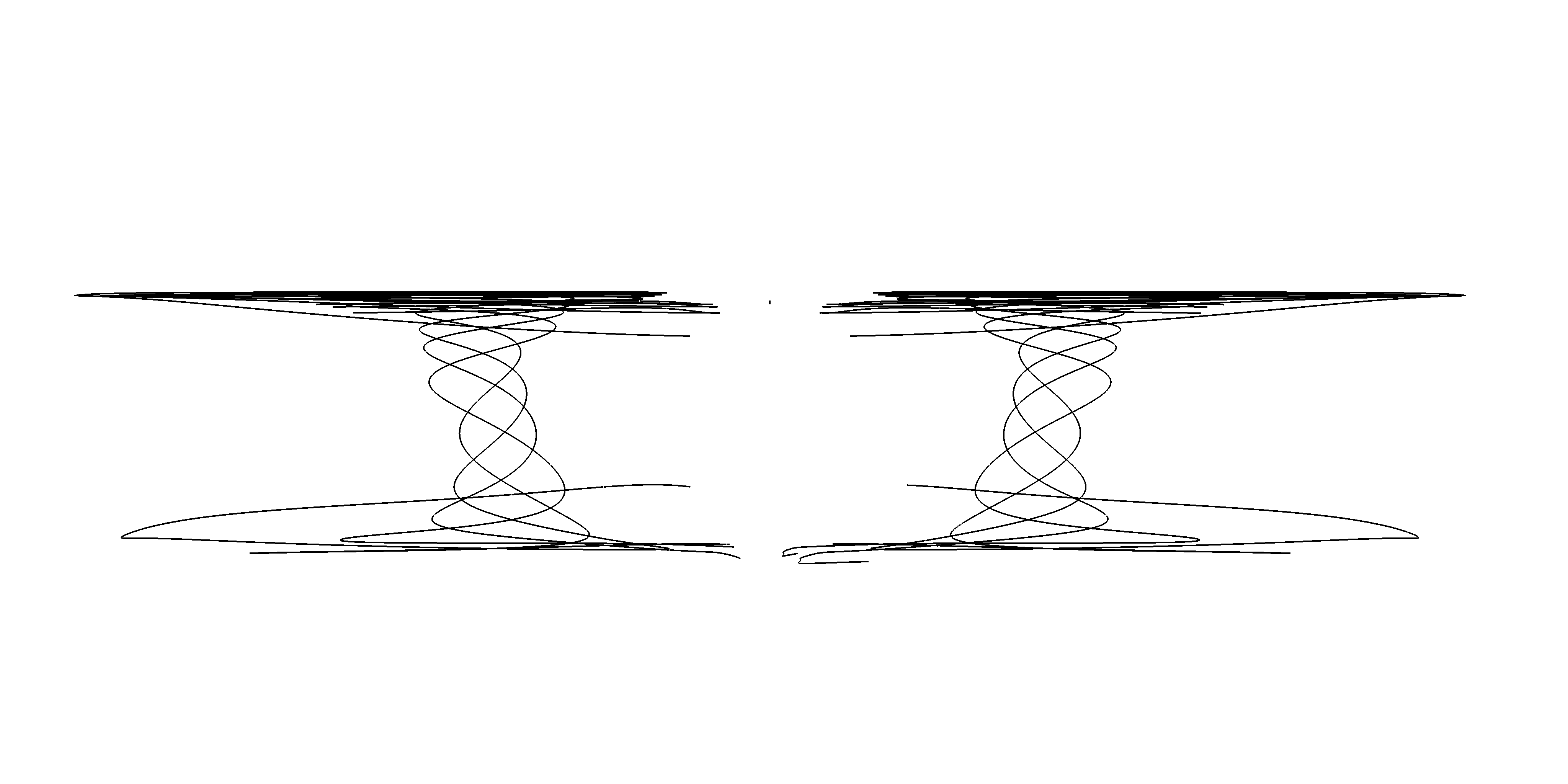}
}
  \caption{3D case for $\lambda = 15.4$: streamlines on the $xy$-plane (left) and $yz$-plane (right)  
  and (a) and (b)  $\HH=0.8165$, $\re_{3D} = 27.79$, (c) and (d)  $\HH=0.9517$, $\re_{3D} = 27.79$, 
  (e) and (f)  $\HH=0.9517$, $\re_{3D} = 45.01$, (g) and (h)  $\HH=0.8165$, $\re_{3D} = 62.22$.
The projection on the $yz$-plane for symmetry reasons shows only half of the geometry.}
  \label{fig:cmp_vortices}
\end{figure}

Next, we let both the geometric parameter $\HH$ and the Reynolds number vary. 
We display in figure~\ref{fig:cmp_vortices} the streamlines on the $xy$-plane (left) and $yz$-plane
for representative values of the two parameters.
For low values of $\HH$ and $\re_{3D}$, the flow develops
without forming vortices, with the streamlines deviating only slightly out of plane. 
See Fig.~\ref{fig:cmp_vortices}(a) and (b), which have been obtained for $\re_{3D} = 27.79$ and $\HH = 0.8165$.
As the channel increases in width, the streamlines gradually become 
fully three-dimensional, especially in the vortex region. 
See Fig.~\ref{fig:cmp_vortices}(c) through (f).
Notice how Fig.~\ref{fig:cmp_vortices}(c) and (d), obtained for $\HH = 0.9517$ and $\re_{3D} = 27.79$,
differ from Fig.~\ref{fig:cmp_vortices}(a) and (b), obtained for the same Reynolds
number but a smaller $\HH$.
The corner vortices in Fig.~\ref{fig:cmp_vortices}(e) looks similar to the recirculations observed in 2D (see, e.g., Fig.~\ref{fig:snapshots_2d}(b)).
However, in a 3D geometry the presence of a top and bottom bounding wall leads to
complex 3D spiraling recirculation structures \cite{Chiang2000,Tsai2006}, 
as shown in Fig.~\ref{fig:cmp_vortices}(f). See also Fig.~\ref{fig:cmp_vortices}(g) and (h).

\begin{figure}[h!]
  \centering
  \subfigure[$xy$-plane, unstable solution]{
  \includegraphics[width=0.45\columnwidth]{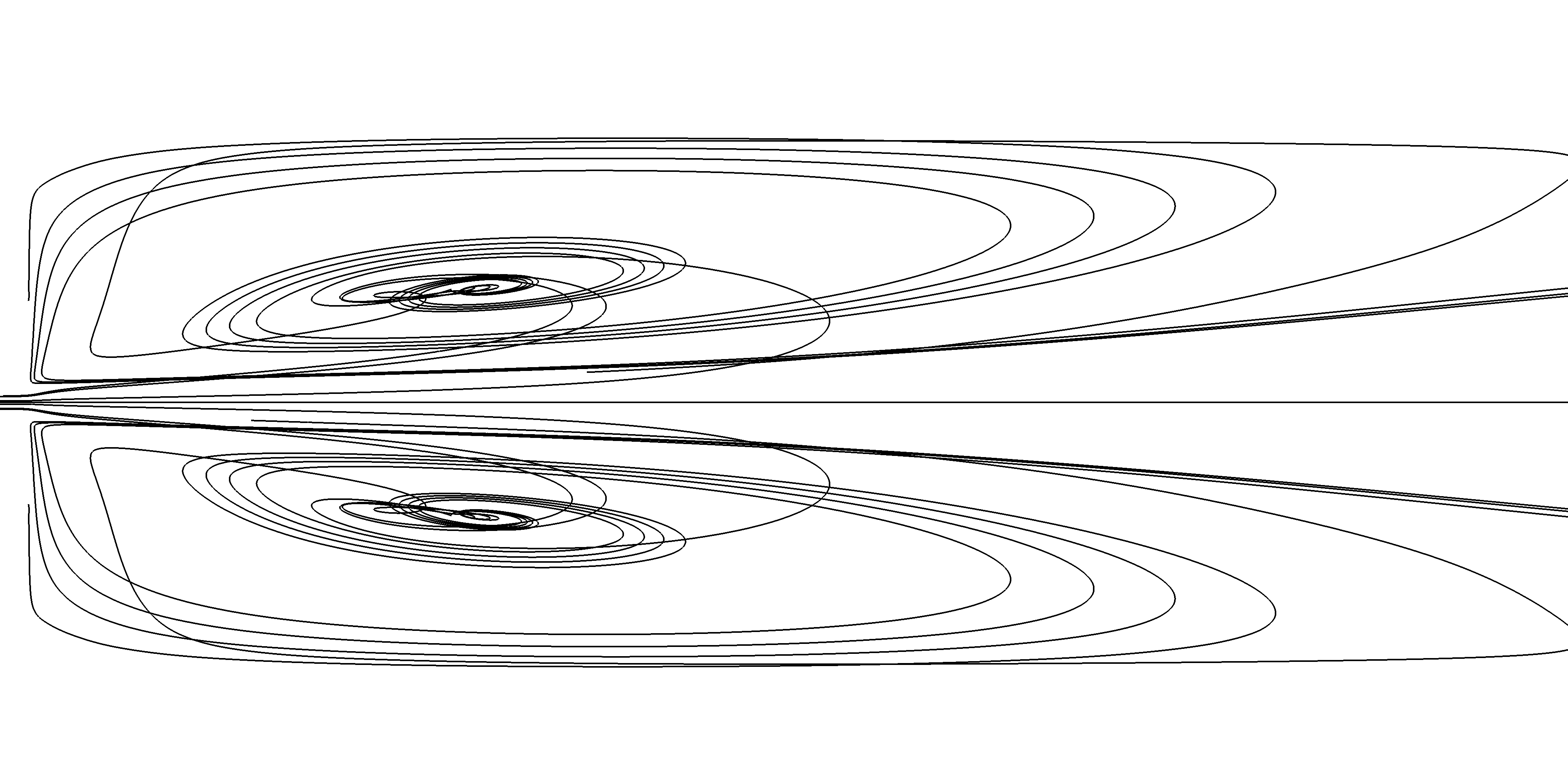}
  }
  \subfigure[$yz$-plane, unstable solution]{
  \includegraphics[width=0.45\columnwidth]{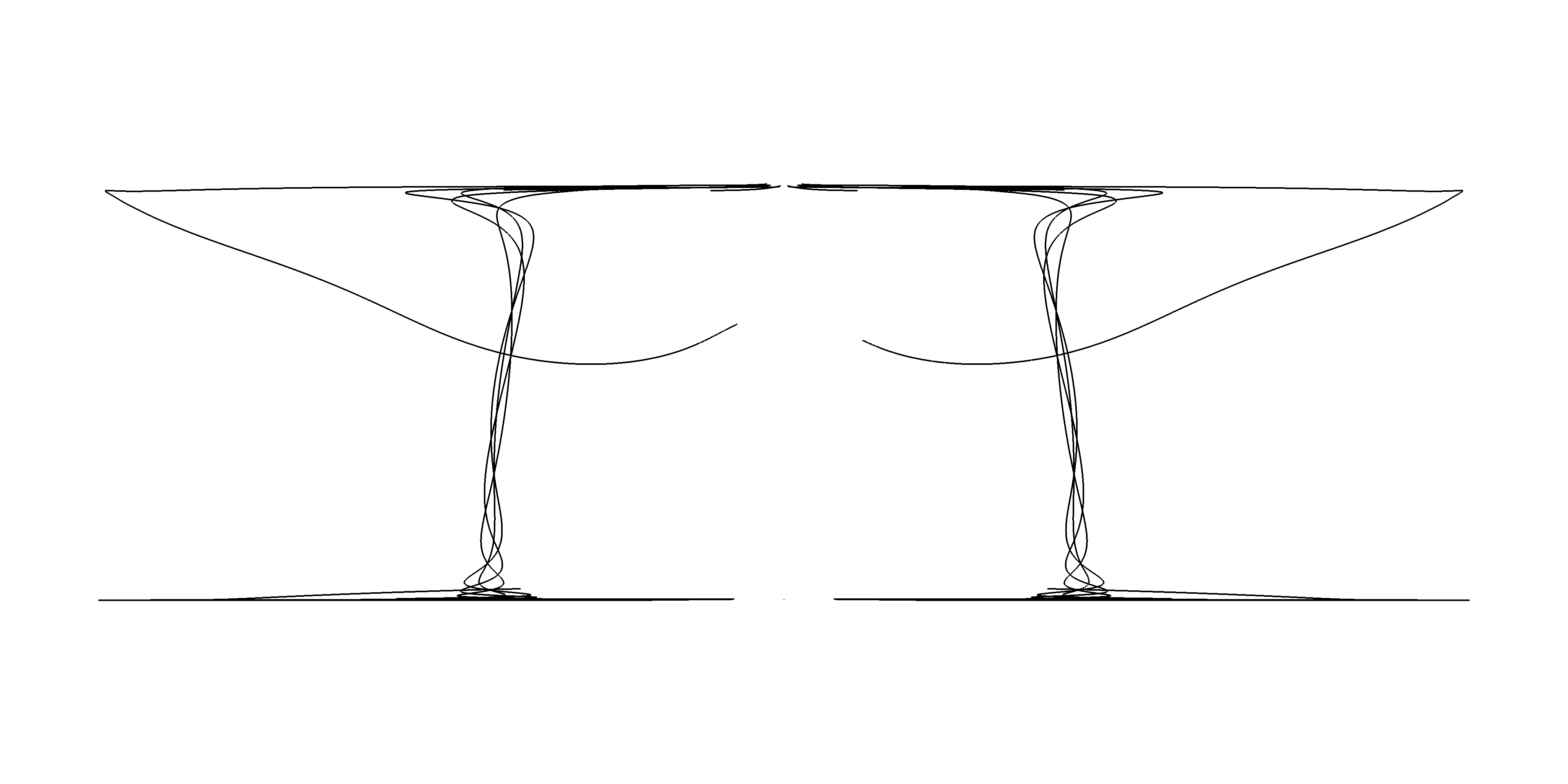}
  }
  \subfigure[$xy$-plane, stable solution]{
  \includegraphics[width=0.45\columnwidth]{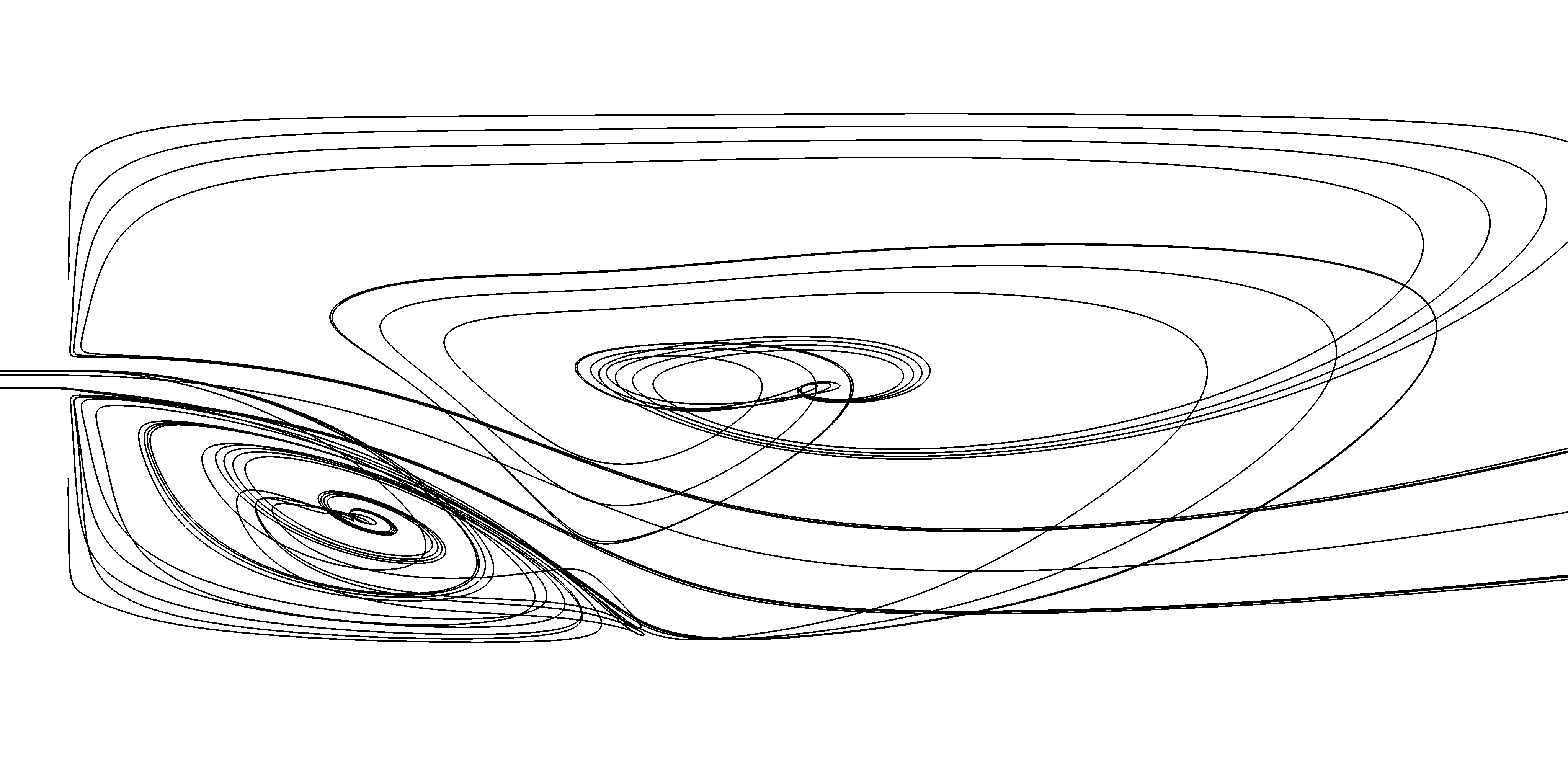}
  }
  \subfigure[$yz$-plane, stable solution]{
  \includegraphics[width=0.45\columnwidth]{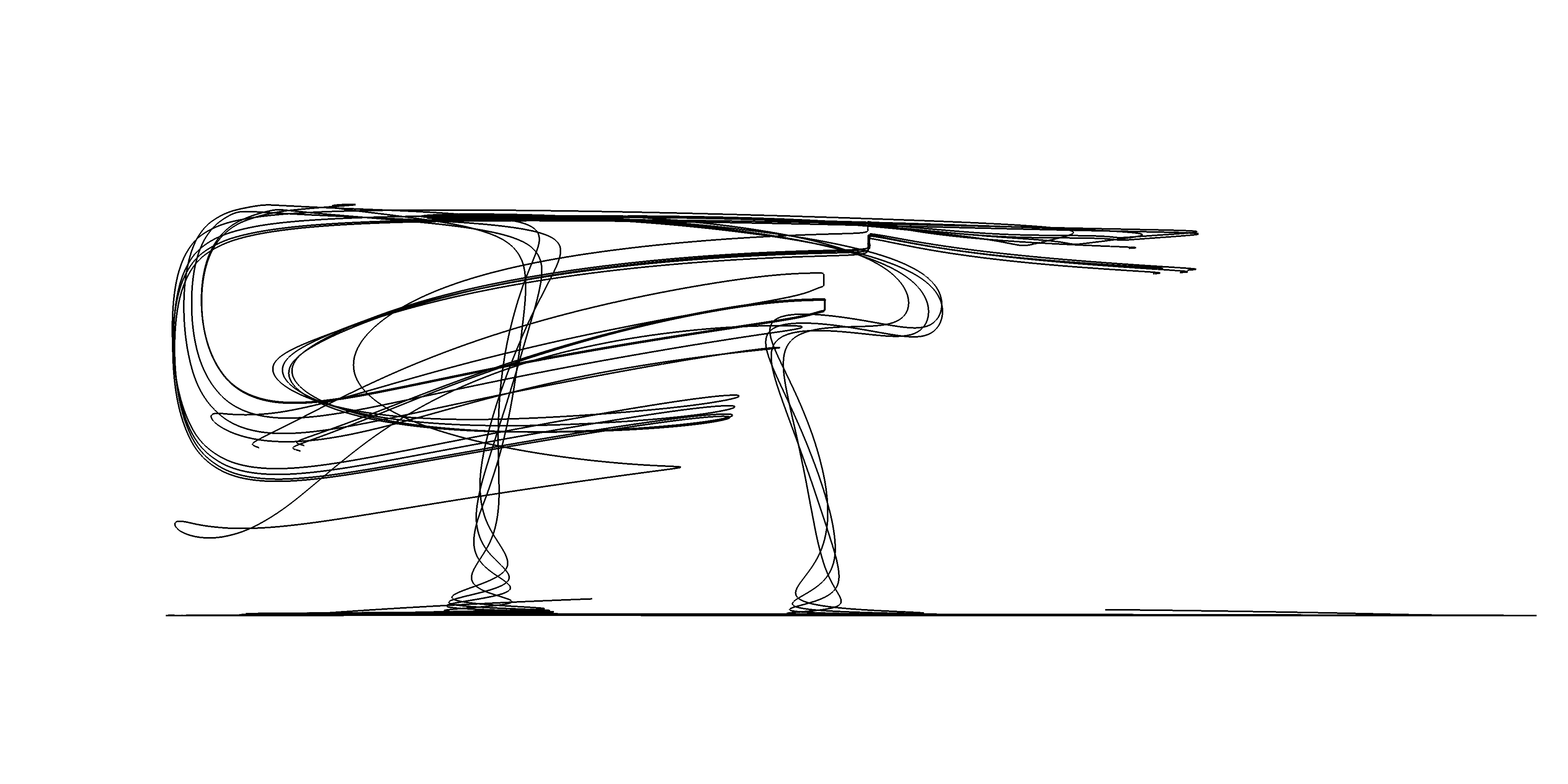}
  }
  \caption{3D case for $\lambda = 15.4$, $\HH=0.9517$, and $\re_{3D} = 76.82$: 
  streamlines on the $xy$-plane (left) and $yz$-plane (right) 
  for (a) and (b) unstable solution, (c) and (d)  stable solution.
The projection on the $yz$-plane for symmetry reasons shows only half of the geometry.}
  \label{fig:cmp_stable_unstable}
\end{figure}

Fig.~\ref{fig:cmp_stable_unstable} shows that streamlines on the $xy$-plane and $yz$-plane
of both the unstable (symmetric) and the stable (asymmetric) solution
for a value of the Reynolds number ($\re_{3D} = 76.82$) past the bifurcation point,
the usual expansion ratio $\lambda = 15.4$, and $\HH=0.9517$.
The vortex pattern becomes even more intricate after the bifurcation point, with the vortices promoting 
the mixing between distant regions of the channel. Due to the symmetry of the 
geometry and the boundary conditions, there is no flow crossing the midline $xy$ plane. 

Each 3D full order computation requires about $240h$ of CPU time, and the preprocessing time is about $40h$.
The computational time savings estimate for the two parameter ($\re$ and $\HH$) case
is given by:
\[
\frac{\text{time to build the RB spaces}+\text{online time to detect the bifurcation point}}{\text{time of the equivalent full order computation}}=\frac{56\cdot240h+40h+0.05h}{7\cdot7\cdot10\cdot240h}\simeq 11.4\%.
\]
where, based on the experience acquired with the 2D case, we suppose that 7 runs per each 
parameter are required to have a reasonable tracking of the bifurcation points in the parameter space. 
With 2 parameters this amounts to 49 runs, each run requiring on average 10 full simulations.
Thus, in the 3D case the \emph{break-even} is given by:
\[
  \frac{\text{All full order computations for RB prep.}}{\text{Full order one query  comp. time}}=\frac{56\cdot240h+40 h}{240h}\simeq 56.2.
\]
The interpretation of this result is that a reduced order model can be expected to bring savings if more than 56 runs are planned.

To test our method, we select a geometric aspect ratio 
not considered in the sampling phase, and we try to recover some characterizing 
flow features as a function of the Reynolds number. We consider $\AR=2.12$ 
(corresponding to $\HH=0.679$) and we reconstruct the profile of the normalized axial velocity:
\begin{equation}
  \frac{v_x}{\langle v_x\rangle_c}=v_x\frac{\int_{\Om\cap\mathit{\Pi}_c}v_x\de\bs{x}}{| \Om\cap\mathit{\Pi}_c|},
  \label{eq:def_axial_vx}
\end{equation}
where $\mathit{\Pi}_c$ is any plane crossing the contraction section and orthogonal to the channel axis 
and $|\Om\cap\mathit{\Pi}_c|$ is the measure of the intersection between the plane $\mathit{\Pi}_c$ and the domain $\Om$.
We also consider the normalized axial velocity gradient:
\begin{equation}
  \frac{\partial_x v_x}{\langle v_x\rangle_c}w_c.
  \label{eq:def_axial_grad}
\end{equation}

\begin{figure}[h!]
  \centering
  \includegraphics[width=0.45\columnwidth]{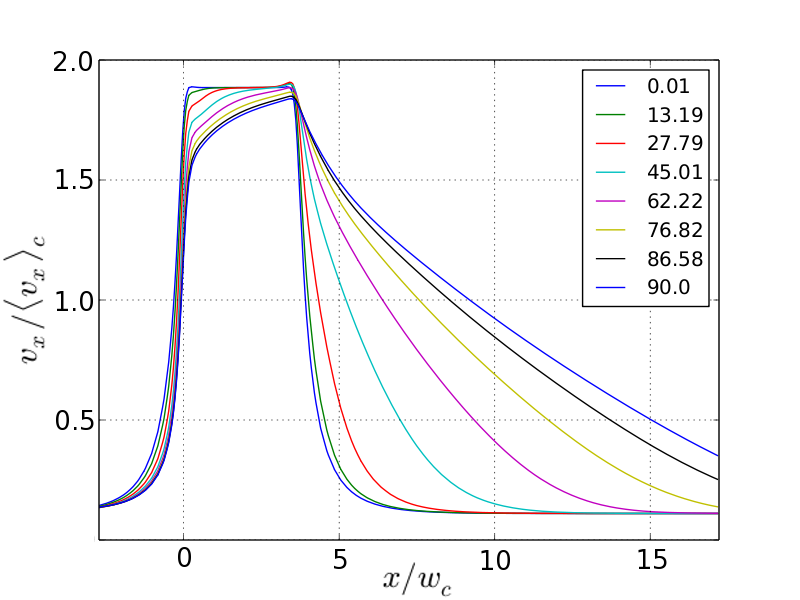}
  \includegraphics[width=0.45\columnwidth]{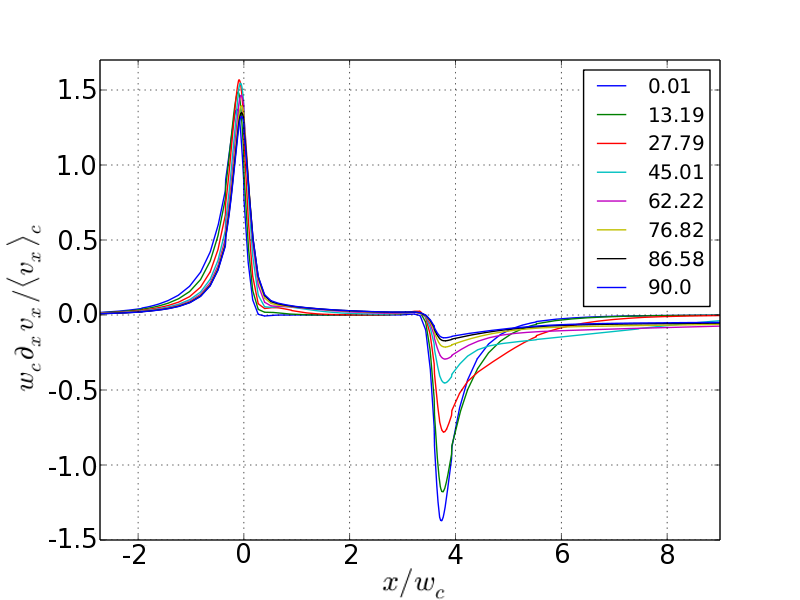}
  \caption{3D case for $\lambda = 15.4$ and $\HH=0.9517$: normalized axial velocity 
  (\ref{eq:def_axial_vx}) and normalized axial gradient (\ref{eq:def_axial_grad}) as 
  a function of the normalized distance from the contraction inlet. 
  The curves in the two figures are computed for values of the Reynolds 
  number between $0.01$ and $90$. The different curves refer to the values of 
  Reynolds number reported in the legend.} 
  \label{fig:axial}
\end{figure}

We plot the normalized axial velocity (\ref{eq:def_axial_vx}) and normalized 
axial gradient (\ref{eq:def_axial_grad}) along the center line for different values of the Reynolds number
in figure~\ref{fig:axial}(a) and (b), respectively. 
The results are in good qualitative agreement with those reported in \cite{Oliveira}.
Concerning the normalized axial velocity, for small Reynolds numbers the curve is almost a 
symmetric step function, since the viscosity is sufficiently high to avoid large velocity gradients 
both inside the cross-section and along the channel length. 
As the Reynolds number is increased, the curve becomes more and more asymmetric, 
and the averaging effect of the viscosity takes longer to smooth out the velocity gradients. 
This is visible from the long tail of the curves with higher Reynolds number.
The viscosity has also a clear effect on the normalized axial gradient in Fig.~\ref{fig:axial}(b):
the two spikes show that the velocity gradients in proximity of the variations in channel width 
increase as the Reynolds number increases.
We remark that the graphs in figure~\ref{fig:axial} can be easily drawn 
by saving the normalized axial velocity and normalized axial gradient for 
the RB functions and using these as to interpolate the desired output in real time. 
This feature is particularly interesting in the real-time query case, since it does not need 
to search a large database during the postprocessing phase.

\begin{figure}[h!]
  \centering
  \includegraphics[width=0.33\columnwidth]{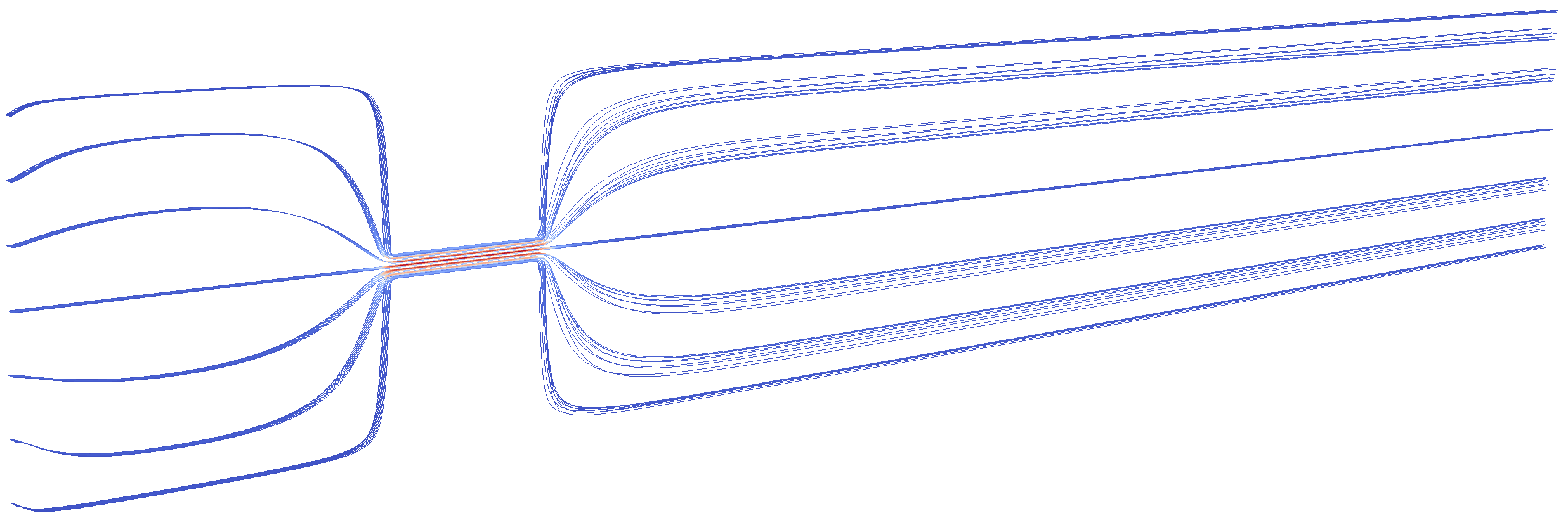}
  \includegraphics[width=0.33\columnwidth]{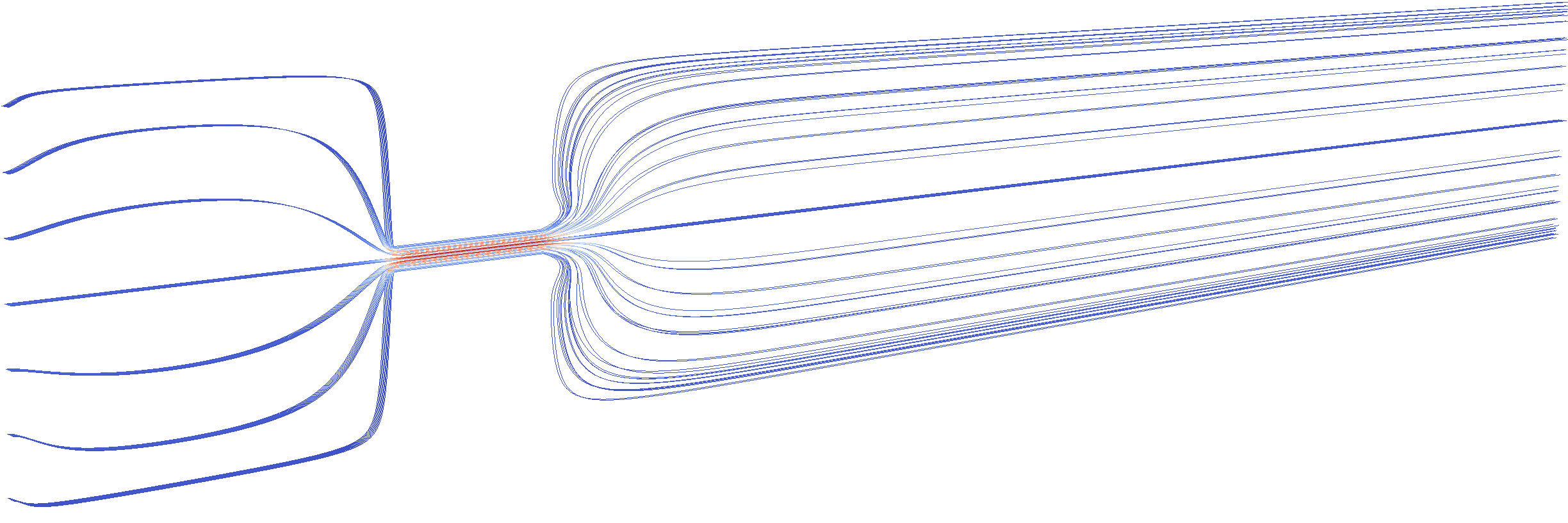}
   \includegraphics[width=0.33\columnwidth]{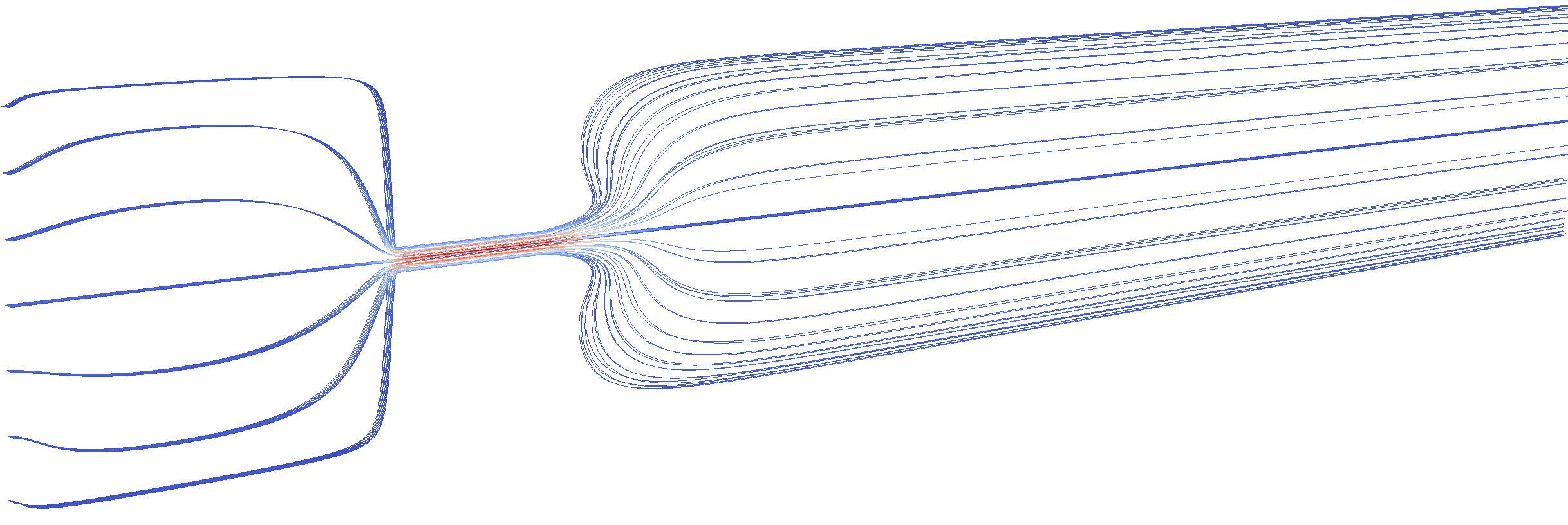}
  \caption{3D case for $\lambda = 15.4$ and $\HH=0.2085$ streamlines for $\re_{3D}=0.01$ (left), $\re_{3D}=23$ (center), and $\re_{3D}=90$ (right).}
  \label{fig:3D_1}
\end{figure}

\begin{figure}[h!]
  \centering
  \includegraphics[width=0.33\columnwidth]{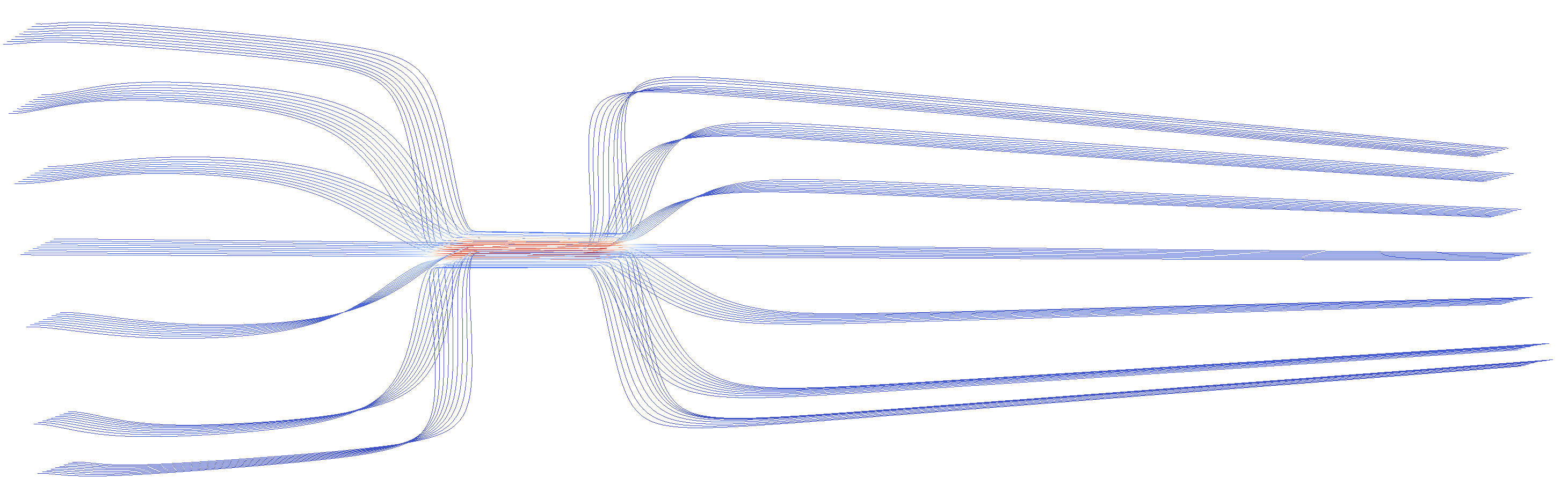}
   \includegraphics[width=0.33\columnwidth]{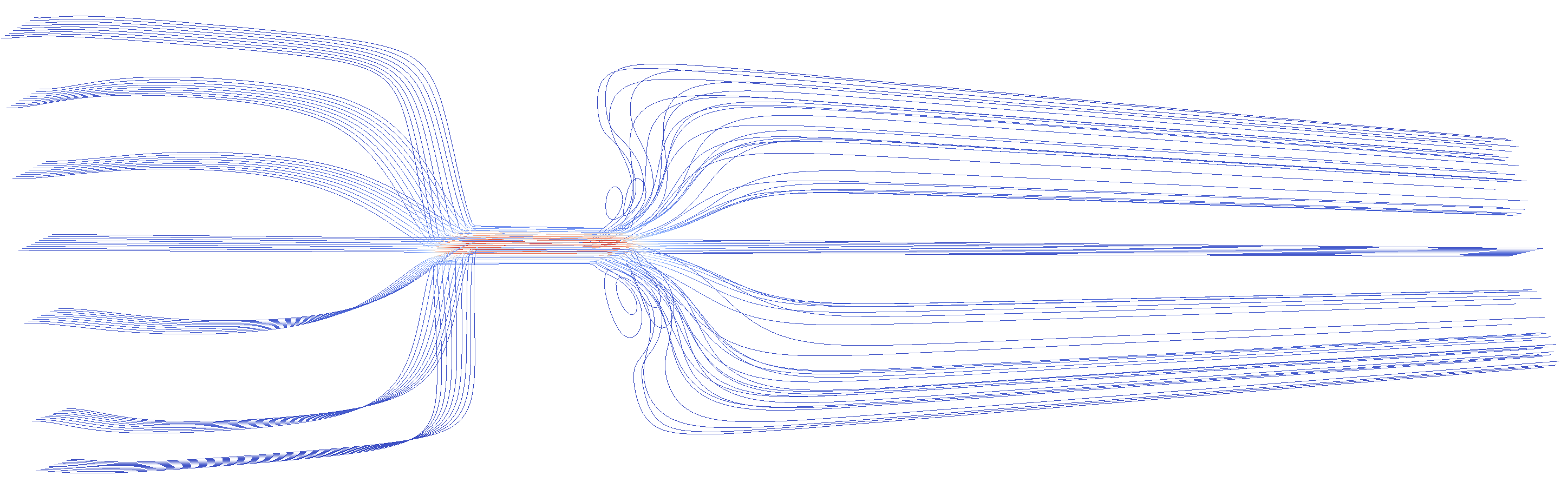}
    \includegraphics[width=0.33\columnwidth]{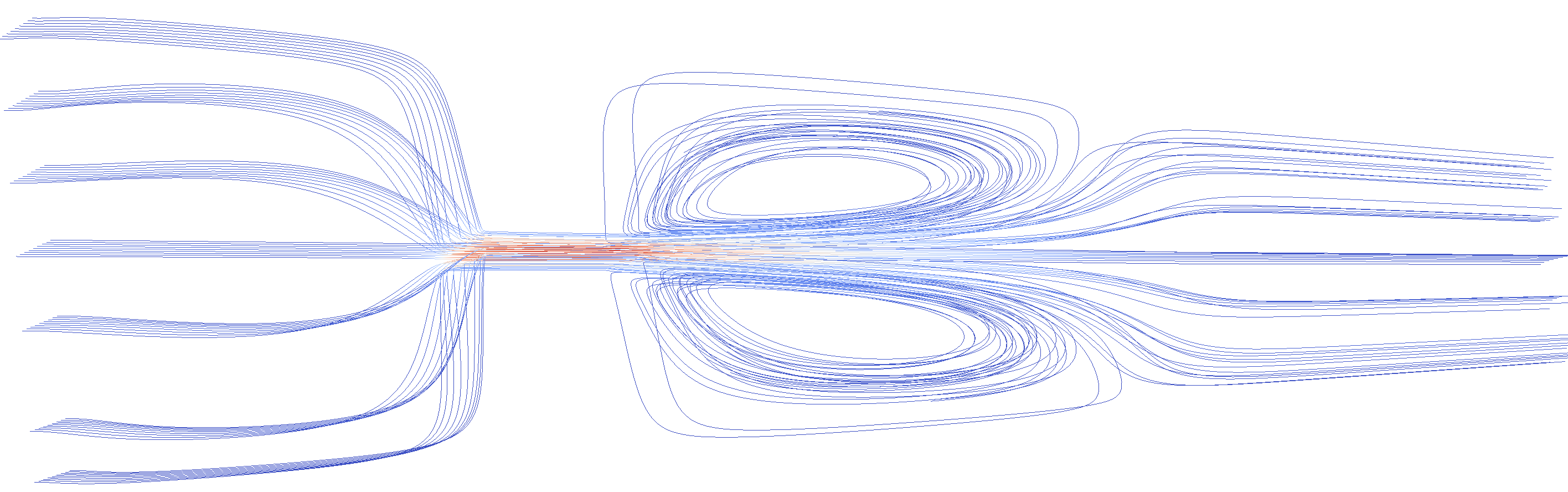}
  \caption{3D case for $\lambda = 15.4$ and $\HH=0.6210$: streamlines for $\re_{3D}=0.01$ (left), (b) $\re_{3D}=13$ (center), and $\re_{3D}=90$ (right).}
  \label{fig:3D_2}
\end{figure}

\begin{figure}[h!]
  \centering
  \includegraphics[width=0.33\columnwidth]{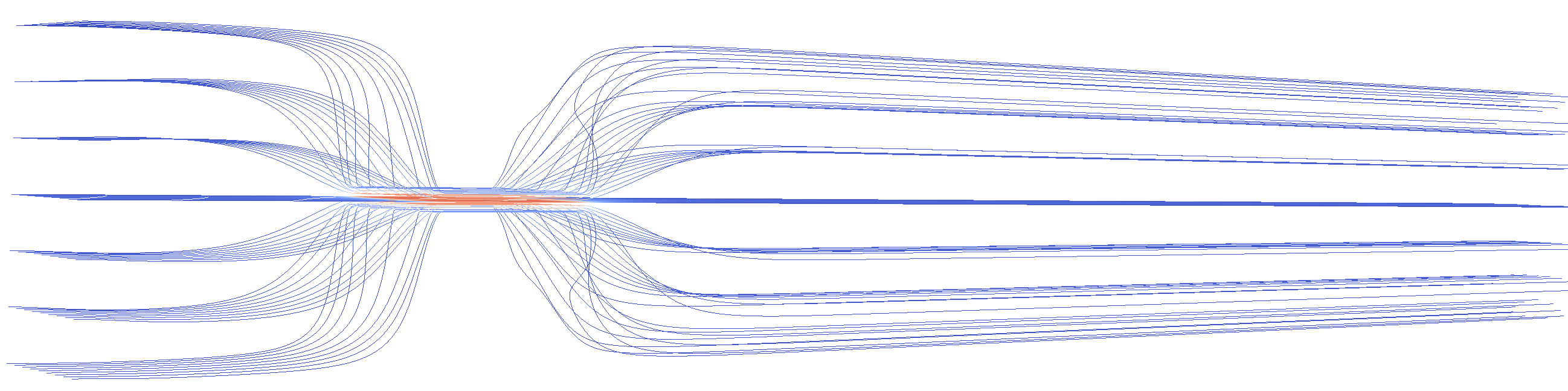}
  \includegraphics[width=0.33\columnwidth]{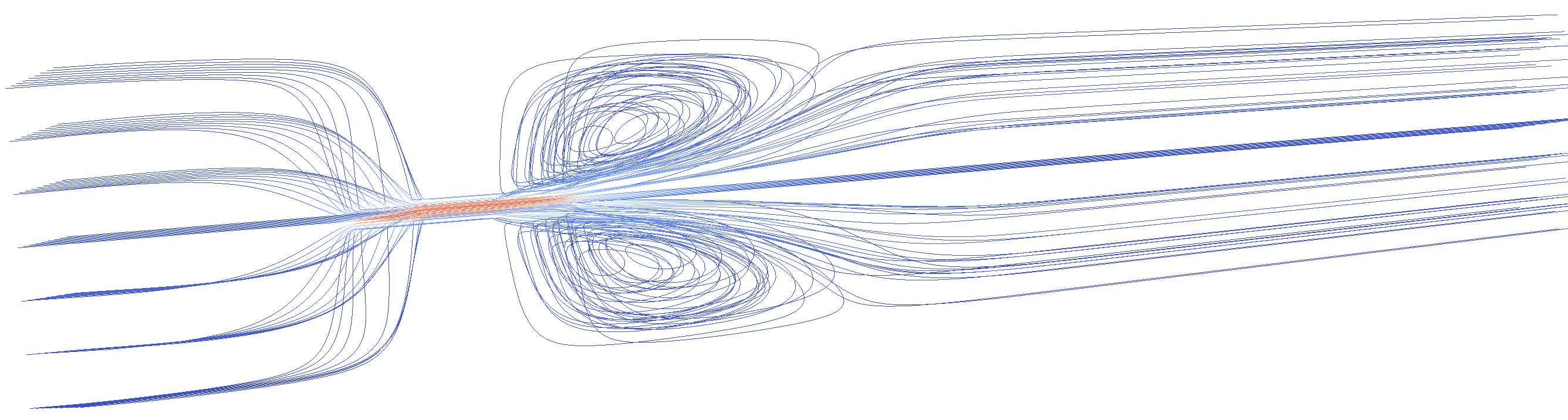}
     \includegraphics[width=0.33\columnwidth]{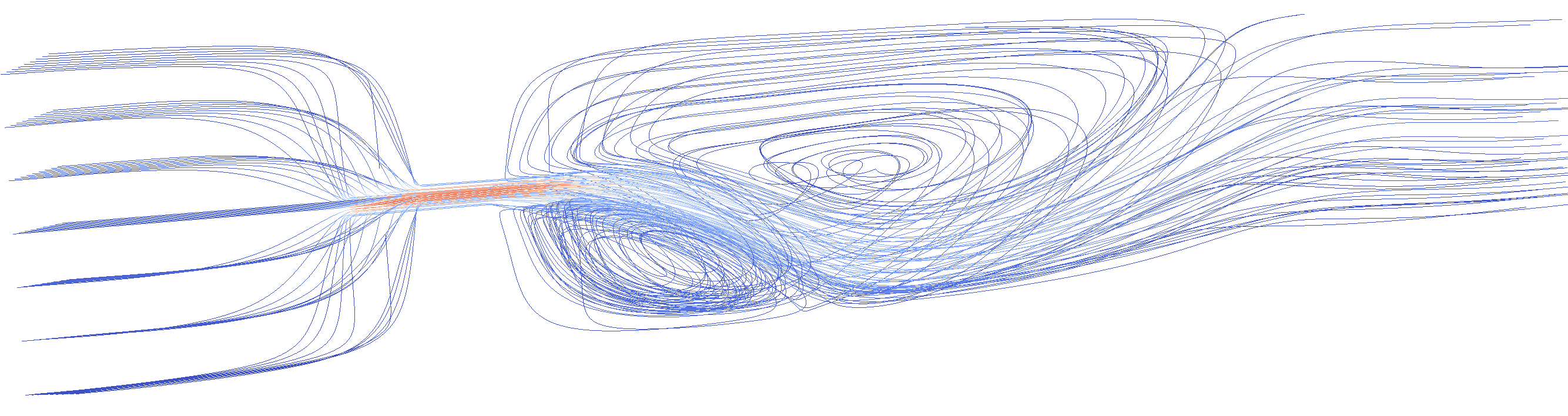}
     \caption{3D case for $\lambda = 15.4$ and $\HH=0.9517$: streamlines for $\re_{3D}=4.466$ (left), $\re_{3D}=27$ (center), and $\re_{3D}=90$ (right).
     }
  \label{fig:3D_4}
\end{figure}

We conclude the section with the streamlines for the flow associated to 
$\HH = 0.2085$ (in Fig.~\ref{fig:3D_1}), $\HH=0.6210$ (in Fig.~\ref{fig:3D_2}), and $\HH=0.9517$ (in Fig.~\ref{fig:3D_4})
for a small value, a medium value, and a large of $\re_{3D} \in [0.01, 90]$.
In particular, compare the solutions for $\re_{3D} = 90$ (leftmost panel in Fig.~\ref{fig:3D_1},~\ref{fig:3D_2}, and \ref{fig:3D_4}).
They clearly show that at low values of $\HH$ the symmetry breaking bifurcation
is pushed to higher values of $\re_{3D}$ due the vertical walls.



%

Based on the results presented in Sec.~\ref{2d_results_2param} and \ref{3d_results},
we conclude that eccentric mitral regurgitant jets are produced by long (large $\HH$) and narrow
(large $\lambda$) orifices. In fact, such slender orifices
associated with eccentric jets, seem to resemble the coaptation geometry of the mitral
valve. Coaptation is the region where the two leaflets of the mitral valve meet
Our hypothesis is that Coanda effect occurs in mitral valves in which the leakage, i.e.,
regurgitation, occurs along a large section of the coaptation zone, rather than at an isolated
point, leading to a possibly significant regurgitant volume. This is corroborated by clinical
observations indicating that eccentric regurgitant jets are, indeed, prevalent in patients
with severe MR \cite{littlei3,chandra2011three,shanks2010quantitative}.

\begin{figure}[h!]
  \centering
  \subfigure[3D printed plate with orifice]{
   \includegraphics[width=0.24\columnwidth]{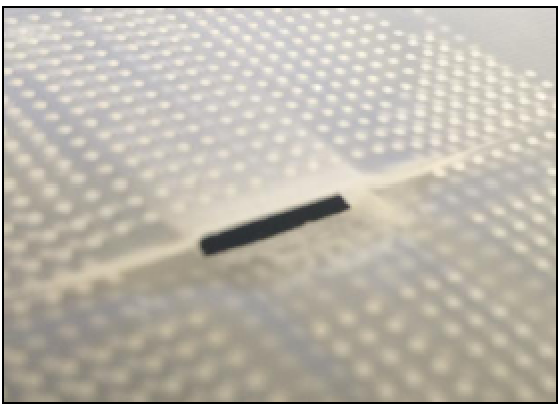}
  }
\subfigure[Mock heart chamber]{ 
  \includegraphics[width=0.37\columnwidth]{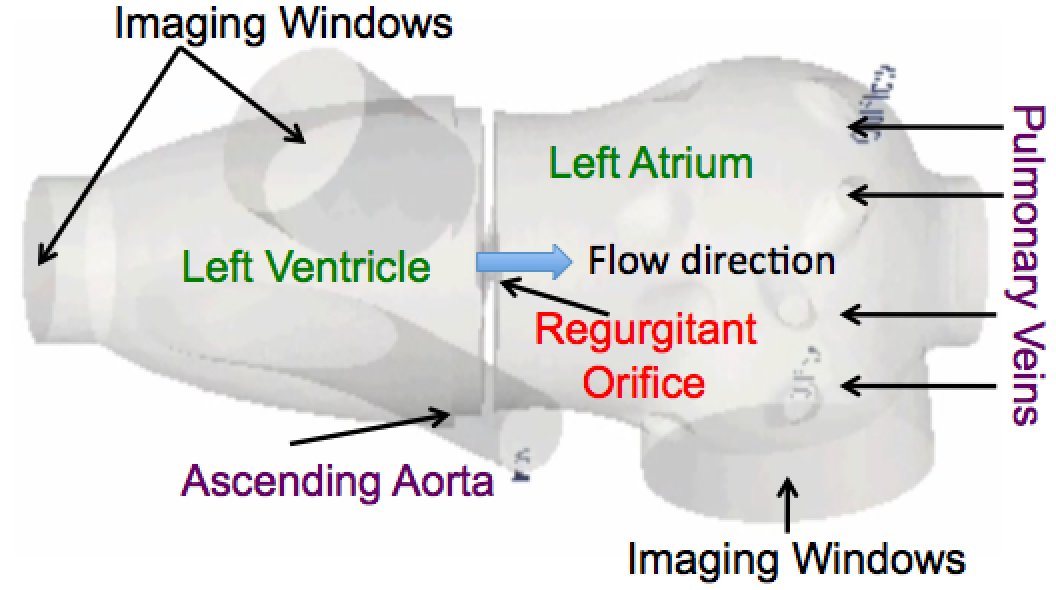}
  }
  \subfigure[Mock heart chamber]{ 
  \includegraphics[width=0.22\columnwidth]{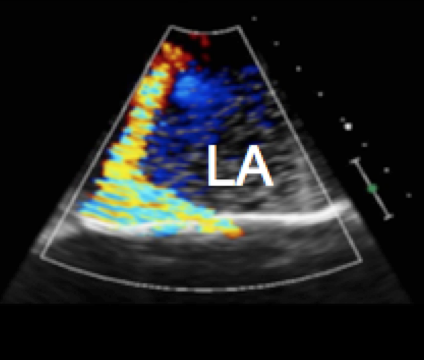}
  }
     \caption{
     (a) Close-up view of the 3D printed divider plate with a long and narrow orifice to mimic a closed leaky mitral valve,
     (b) geometry of the mock heart chamber with the divider plate between mock left ventricle and mock left atrium (LA), and
     (c) 2D Doppler echocardiographic image of the regurgitant jet in the mock heart chamber. Conic distortion in (c) occurs due to the use of convex array transducer.
     }
  \label{fig:exp_setup}
\end{figure}

Before the study presented in this manuscript, our collaborators at the Houston Methodist DeBakey Heart \& Vascular Center 
had never succeeded in reproducing the Coanda effect in vitro. Following our results, they
designed a long and narrow orifice in a divider plate that mimics a closed leaky mitral valve. A close-up view
of the orifice is in Fig.~\ref{fig:exp_setup}(a). The divider plate was mounted
on an anatomically correct mock (left) heart chamber developed to study the use of 2D and 
3D color Doppler techniques in imaging the clinically relevant intra-cardiac flow events 
associated with regurgitant jets \cite{littlei1,littlei2}. See Figure~\ref{fig:exp_setup}(b). 
The chamber is connected to a pulsatile flow loop. 
The fluid in the mock heart chamber is water with 30\% glycerin added to mimic blood
viscosity. Notice that this is consistent with modeling blood as a Newtonian fluid in Sec.~\ref{setting}.
From the 2D Doppler echocardiographic image
in Fig.~\ref{fig:exp_setup}(c) we see that indeed the slender orifice in
Fig.~\ref{fig:exp_setup}(a) generates a regurgitant jet that hugs the wall.
See also \cite{WangQuaini}. We expect also that these studies could enhance in the near future in vivo studies and applications.

%% file: tex/conclusions.tex
\section{Conclusions and perspectives}\label{conclusions}

The symmetry breaking bifurcation (Coanda effect) has been studied in parametric flows, representing a simplified test case
for regurgitant mitral valve flows.
Our preliminary work shows that standard reduced order methods (e.g., Reduced Basis and /or Proper Orthogonal Decomposition)
allow to capture complex physical and mathematical phenomena, such as bifurcations in the parametrized
Navier-Stokes equations, at a fraction of the computational cost required by full order order methods.
In order to detect the bifurcation points, the reduced parametric Navier-Stokes equations have been 
supplemented with a generalized eigenvalue problem, also cast into the reduced order setting.
This work is also an example of computational collaboration between high performance computing and reduced order methods:
thanks to the computational gains with the same resources we can treat more complex problems. 
This computational collaboration has demonstrated the ability to provide reliable and accurate results with significant reduction of computational times. Results have been validated both with the full-order model and by comparison with parametric studies available in literature for both 2D and 3D cases.

Research perspectives in this field include the development of proper error bounds for the detection of the bifurcation points and the verification of the accuracy. At the state of the art this aspect is carried out by supplementing the state equation with a generalized eigenvalue problem, solved with the same reduced order method proposed for the state equation.
Moreover, we plan on taking into account the interaction of the fluid with elastic walls (i.e., elastic valve leaflet) \cite{FLD:FLD4252}.
This would lead to important improvements in the study of this complex multiphysics nonlinear problem and a better understanding 
of how the Coanda effect is influenced by the valve elasticity.


%% file: tex/perspectives.tex
\section{Acknowledgements}
The authors want to thank Prof. S. Canic, Prof. R. Glowinski (University of Houston) and S. Little MD (The Methodist Hospital, Houston) 
for the fruitful discussions. 
The research in this work has been partially supported by the National Science Foundation under grants DMS-1620384, DMS-1263572 and  DMS-1109189 (Quaini),  INDAM-GNCS 2015 project ``Computational Reduction Strategies for CFD and Fluid-Structure Interaction Problems'', by the INDAM-GNCS 2016 projects ``Tecniche di riduzione della complessit\`{a} computazionale per le scienze applicate'', by PRIN project ``Mathematical and numerical modeling of the cardiovascular system, and their clinical applications'', and by European Union Funding for Research and Innovation -- Horizon 2020 Program -- in the framework of European Research Council Executive Agency: H2020 ERC CoG 2015 AROMA-CFD project 681447 ``Advanced Reduced Order Methods with Applications in Computational Fluid Dynamics''. Computations have been performed on the SISSA cluster \emph{Ulysses} and on the CINECA clusters (COGESTRA project 2015).